\documentclass[review]{elsarticle}
\usepackage{setspace}
\usepackage{lineno,hyperref}
\usepackage{amssymb,amsfonts,amsthm}
\usepackage{amsmath}
\usepackage{graphicx}
\usepackage{float}
\usepackage{subcaption}
\usepackage{parselines}
\modulolinenumbers[5]

\usepackage[margin=1.0in]{geometry}
\begin{document}

\title{Neural Networks with Inputs Based on Domain of Dependence and A Converging Sequence for Solving Conservation Laws, Part I: 1D Riemann Problems}
\author{Haoxiang Huang\footnotemark[1],\; Yingjie Liu\footnotemark[2] \; and Vigor Yang\footnotemark[3]}

\begin{abstract}
Recent research works for solving partial differential equations (PDEs) with deep neural networks (DNNs) have demonstrated that spatiotemporal function approximators defined by auto-differentiation are effective for approximating nonlinear problems, e.g. the Burger’s equation, heat conduction equations, Allen-Cahn and other reaction-diffusion equations, and Navier-Stokes equation. Meanwhile, researchers apply automatic differentiation in physics-informed neural network (PINN) to solve nonlinear hyperbolic systems based on conservation laws with highly discontinuous transition, such as Riemann problem, by inverse problem formulation in data-driven approach. 
However, it remains a challenge for forward methods using DNNs without knowing part of the solution to resolve discontinuities in nonlinear conservation laws. In this study, we incorporate 1st order numerical schemes into DNNs to set up the loss functional approximator instead of auto-differentiation from traditional deep learning framework, e.g. TensorFlow package, which improves the effectiveness of capturing discontinuities in Riemann problems. 
In particular, the 2-Coarse-Grid neural network (2CGNN) and 2-Diffusion-Coefficient neural network (2DCNN) are introduced in this work. We use $2$ solutions of a conservation law from a converging sequence, computed from a low-cost numerical scheme, and in a domain of dependence of a space-time grid point as the input for a neural network to predict its high-fidelity solution at the grid point. Despite smeared input solutions, they output sharp approximations to solutions containing shocks and contacts and are efficient to use once trained.
\end{abstract}

\maketitle

\renewcommand{\thefootnote}{\fnsymbol{footnote}}
\footnotetext{1. {\tt E-mail: hcwong@gatech.edu}. Woodruff School of Mechanical Engineering, Georgia Institute of Technology, Atlanta, GA 30332, USA.}
\footnotetext{2. {\tt E-mail: yingjie@math.gatech.edu}. School of Mathematics, Georgia Institute of Technology, Atlanta, GA 30332, USA. Research is supported in part by NSF grant DMS-1522585.}
\footnotetext{3. {\tt Email: vigor.yang@aerospace.gatech.edu}. Daniel Guggenheim School of Aerospace Engineering, Georgia Institute of Technology, Atlanta, GA 30332, USA.}
\renewcommand{\thefootnote}{\arabic{footnote}}

\section{Introduction}
There has been a lot of research works on numerical methods for solving conservation laws whose solutions may contain shocks and contact discontinuities, such as the Godunov scheme~\cite{Godunov59}, MUSCL scheme~\cite{vanLeer73,vanLeer79}, ENO~\cite{ENO87,ShuOsh88} and WENO~\cite{WENO,JiShu96} schemes, hierarchical reconstruction~\cite{LiuShu07,LiuShu07b,XuLiuShu09}, and many others. Numerical techniques are important for studying high-speed aerodynamic flows which plays a substantial role in aircraft designs, combustion problems and astronomy physics ~\cite{YANG00, XWang17, Umesh21, NumStarForm19}. However, the development of machine learning techniques for solving hyperbolic conservation laws are still in early stage.  

In the past decade, data driven modeling in machine learning has been widely developed for multiple scientific disciplines, including image processing, biomedical applications and engineering design optimization ~\cite{MLMedicine16, PINNrev21, siamRevOpt18, CNN12, CKSPOD21, KSPOD19, EmulationLES18}. Especially, advances in computational resources, e.g. graphics processing unit (GPU) and tensor processing unit (TPU), accelerate training speed in deep learning frameworks, e.g. TensorFlow, PyTorch, for computer vision, natural language processing and other important scientific disciplines~\cite{CVPR19, NIPS2017_3f5ee243, chen2020learning}. 

Recently, \textit{Raissi et al.} ~\cite{RAISSI2019686} employed deep neural networks with physics informed conditions,  e.g. initial conditions, boundary conditions and functional forms of partial differential equations (PDEs) by applying chain rules for differentiating compositions of functions utilizing automatic differentiation ~\cite{BaydinPearlmutterRadulSiskind17}. This method has achieved many successes in data driven methods for identifying nonlinear PDEs or predicting turbulent mixing, vortex induced vibration (VIV) with given governing PDEs, e.g., the Navier-Stokes equations ~\cite{JIN21, clarkZakiMeneveau_2021, RaissiBabaeeGivi19,RAISSI2019686,  raissiJFM19, Raissi1026}. Meanwhile, PINN has demonstrated that it can solve more physical problems from reality, e.g., Lattice Boltzmann Equation with the Bhatnagar-Gross-Krook collision, hypersonic flow, heat transer, nano-optics, electro-convection and so on ~\cite{lou2020physicsinformed, Chen20, MFNN_CompxFluids21, HT21, mao2020deepmmnet, MFNN_liu19}. Besides PINN developed based on DNN, recent researchers have discovered that physics-informed convolutional neural network (CNN) are capable of resolving physical problems with irregular boundary values ~\cite{PysGeoNet21}. Meanwhile, to explore the possibility of learning temporal stages in the future while data may not be available, \textit{Psaros et al.} ~\cite{psaros21metalearning} extended PINN to the meta-learning framework. In \textit{Nguyen} and \textit{Tsai et al.} \cite{Nguyen2020ASP}, a neural network is developed for solving second order wave equations by using the solution and its gradient computed from a low order scheme as input to predict the higher order solution.

To examine the possibility of applying PINN to solving the Euler equations, \textit{Mao et al.} ~\cite{MAO20} designed forward and inverse methods by selecting clustered training data points in the discontinuity region along with randomly selected data points in other spatio-temporal positions, which enhance the accuracy of predictions for the Lax and Sod problems.  Another study by \textit{Michoski et al.} ~\cite{MICHOSKI20} compares solutions computed by finite volume methods and results from predictions of PINN. In their PINN algorithm design, diffusion terms are introduced in functional forms of Euler equations inside the neural network to smooth out discontinuities, though there is no detailed discussion of the functional forms. 
We are still wondering whether automatic differentiation along works in defining functional forms of Euler equations inside the neural network for solving the equations, without knowing any part of the solution. Is there any efficient methodology using neural networks to compute the Euler equations with results comparable to those from high resolution non-oscillatory schemes? We would like to refer to prior works on locating shocks or contacts from a numerical solution or smeared solution \cite{Harten89,ShenPark07}.   

In our study, we first apply automatic differentiation to approximate Euler equations (with or without artificial diffusion) in the cost functional form of PINN which doesn’t seem to produce good approximation of the solutions. To resolve this issue, 1st order numerical schemes are used to approximate the Euler equations in the cost functional form of PINN. They produce smeared approximations to the solutions comparable to those obtained from corresponding schemes, and are referred to as numerical PINN (N-PINN) throughout the paper. To improve the resolution of N-PINN, scattered data points are randomly selected as the additional training data in the vicinities of discontinuities of the solution in the training process. This greatly increases the resolution of N-PINN. Finally, we introduce our best neural network methods which are efficient once trained and are able to produce sharp results competitive to those from high resolution numerical schemes. In order for the neural network to predict the solution at a space-time grid point, we first compute two coarse grid 1st order numerical solutions and then use them in a carefully selected small domain of dependence of the grid point as input.  Alternatively, instead of using two grids, we can use a low-cost scheme on one grid
to compute a conservation law perturbed with two diffusion coefficients, and use the solutions in a carefully selected small domain of dependence of the grid point as input. We find that the 2-coarse-grid neural network and 2-Diffusion-Coefficient neural network can reduce the error of prediction much better than that from N-PINN even with selected partial data and at the same time reduce training time substantially, which may save time for engineering design optimization.

\section{Problem Setup}
Multi-layer neural networks ~\cite{NN89} has found their application in solving PDEs, e.g., through PINN which uses automatic differentiation (auto-differentiation) to set up functional forms of PDEs \cite{RAISSI2019686}. However, the Euler equations remain a challenge for this approach because of the complexity of the system and 
discontinuities occurring in the solution. We first set up PDE functional forms inside the neural network using 1st order numerical schemes, including the Rusanov scheme, Leapfrog scheme with diffusion terms, and Leapfrog scheme with nonlinearly weighted diffusion terms to approximate the solution.

In this study, we start with 1-D conservation laws for our proposed methods. A scalar conservation law is defined as

\begin{equation}
\label{cons-law}
   \frac{\partial U}{\partial t} + \frac{\partial f(U)}{\partial x} = 0,  x\in\Omega\subset \mathcal{R},  t\in[0, T],
\end{equation}
\\
and the 1-D Euler equations for an ideal gas are    

\begin{equation}
    \frac{\partial}{\partial t}\left(\begin{array}{c} \rho \\ \rho v \\E \end{array}\right) +\frac{\partial }{\partial x} \left(\begin{array}{c} \rho v\\ \rho v^2 \\ v( E + p)\end{array}\right) = 0,\; x\in\Omega\subset \mathcal{R}, \; t\in[0, T]~,
\end{equation}
where $\rho$, $u$ and $p$ are density, velocity
and pressure respectively,  $\Omega$ is an interval,
\begin{equation}
    E = \frac{p}{\gamma - 1}+\frac{1}{2}\rho v^2~,
\end{equation}
$\gamma$ is the specific heat ratio and its value is $1.4$ throughout all the cases studied in the paper.  

To approximate the solution of the Euler equations by DNN, we use physics informed dataset, e.g. initial conditions, boundary conditions and training data from $\rho$, $v$, and $p$ to allow DNN to infer all states of interest in the spatiotemporal domain $(x,t)$.
 We design a DNN with a spatiotemporal point as its input and the prime variables ($\rho$, $v$, and $p$) of the predicted solution at the point as its output. In order to define the loss function of the DNN, we first partition the space-time domain with a rectangular grid, then consider a $1^{st}$ order numerical scheme expressed in terms of the prime variables at grid points (through the output of DNN) approximating each equation of the Euler equations.  The loss function is the summation of two parts. The first part is the loss of the physics-informed dataset from initial conditions, boundary conditions and training data (if available), while the second part is the sum of square of the residue error from the 1st order scheme at each grid point for each equation of the Euler system. 
The weights and biases of the DNN can be learned from minimizing the loss function in this PINN. The loss function is defined as follows.
\\
\begin{equation}
\label{total-loss}
    Err_{total} = Err_1 + Err_2,
\end{equation}
\\
where $Err_1$ is the loss of physics-informed data defined as
\\
\begin{equation}
\label{loss1}
    Err_1 = \frac{1}{N_1}\sum^{N_1}_{i=1}[u(t^i, x^i)-u^i]^2,
\end{equation}
\\
superscript $i$ indicates a spatiotemporal point where the value of a prime variable $u$ is available as $u^i$, $u(t^i, x^i)$ is the corresponding DNN-predicted value at the point, $Err_2$ measures the residue error (denoted as the $1^{st}$ order scheme functional form) defined as
\begin{equation}
\label{loss2}
    Err_2 = \frac{1}{N_2}\sum^{N_2}_{j=1}|L(t^j, x^j)|^2,
\end{equation}
superscript $j$ indicates a spatiotemporal grid point, $L(t^j, x^j)$ is the residue of the scheme at the grid point (in which the solution values inside the scheme are substituted by corresponding DNN-predicted values) and $|\cdot|$ is the Euclidean norm.
We denote $\{t^i,x^i,u^i\}_{i=1}^{N_1}$ as the physics-informed dataset which contains initial and boundary data (known as the forward problem) and possibly additional training data (known as the inverse problem.)

\section{Numerical Physics-Informed Neural Network (N-PINN)}

The original PINN defines the functional form by auto-differentiation. However, auto-differentiation seems to have difficulty in capturing discontinuous solutions for a highly nonlinear system such as the Euler system. We use PINN to predict the solution of the Euler system without any training data point in the interior of the spatiotemporal domain. The schematic below depicts the training procedure for PINN. 

\begin{figure}[H] \centering 
\includegraphics[width=1.0\linewidth]{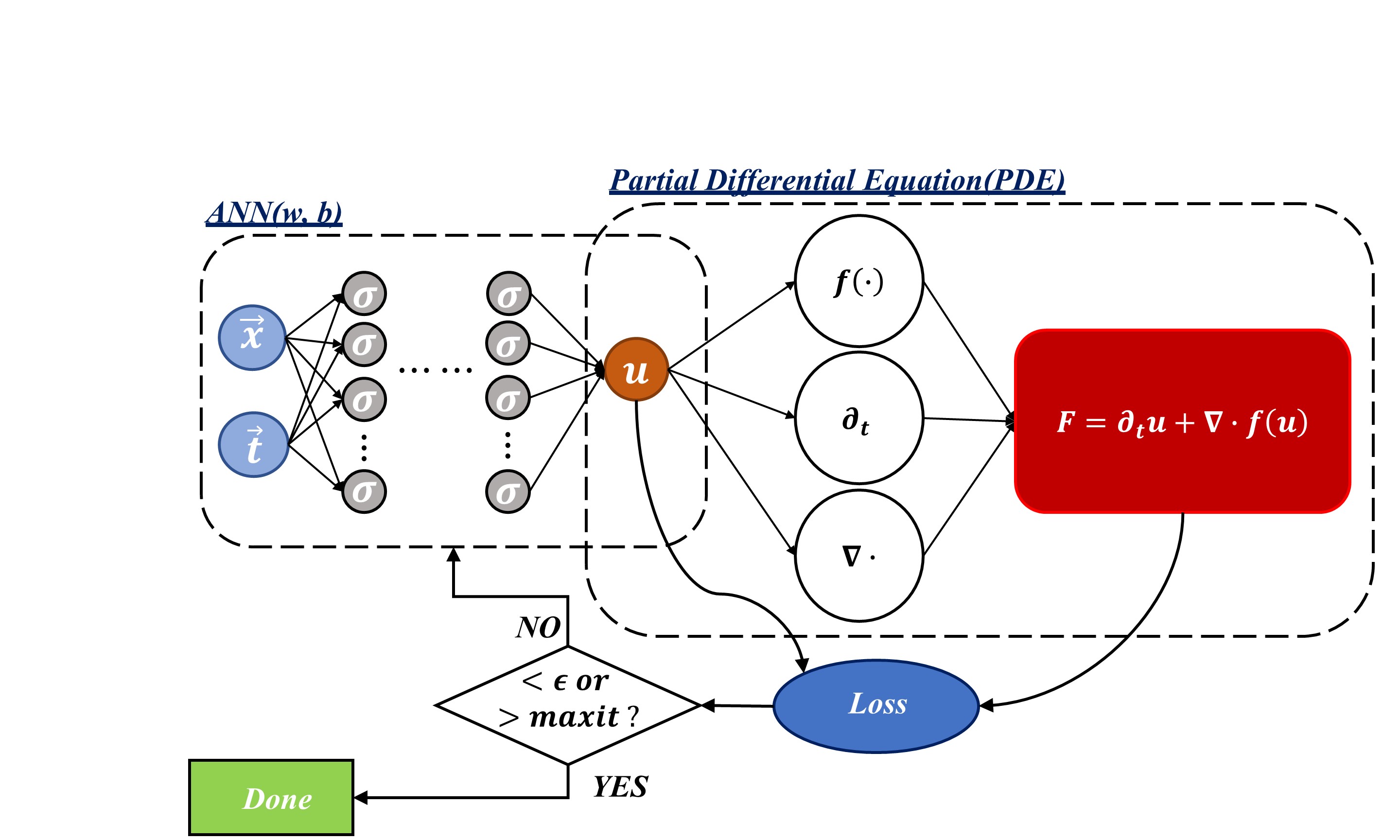}
\caption{Schematic of physics informed neural network (PINN) with PDE functional form (residue) approximated by auto-differentiation.}
\label{fig 1:PINN}
\end{figure}

In this section, We define another PDE functional form for the Euler equations by using a $1^{st}$ order numerical scheme, and the schematic describing the training procedure is modified as followed:

\begin{figure}[H]\centering
\includegraphics[width=1.0\linewidth]{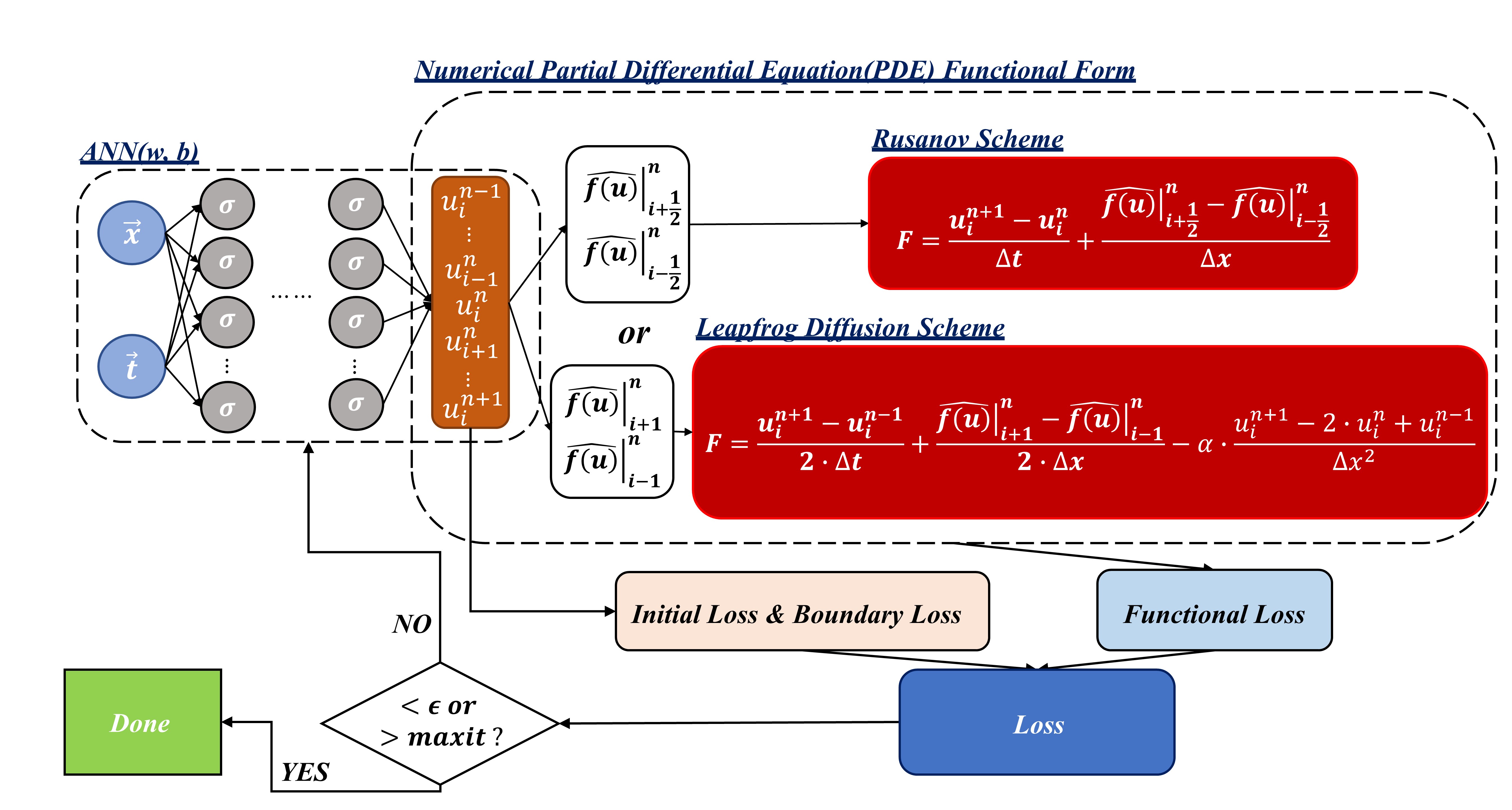}
\caption{Schematic of N-PINN with PDE functional form (residue) approximated by a $1^{st}$ order numerical scheme.}
\label{fig2:NumPINN}
\end{figure}

The numerical scheme used in the functional form of N-PINN is either the Rusanov scheme~\cite{Rusanov61} or a Leapfrog-type scheme (applied to the Euler equations with arificial diffusion terms, see Fig.~\ref{fig2:NumPINN}.) 
Note that the Leapfrog scheme is not stable for computing a diffusive evolution equation. But when using its residue as the loss functional in N-PINN, there is no stability issue. This is a nice feature of N-PINN which allows us to use schemes without worrying about their conventional stability.
The following subsections will discuss the outcomes of predictions by N-PINN.

\subsection{N-PINN with the Rusanov scheme}
\label{Rusanov N-PINN}
The first numerical scheme employed for constructing the functional form in neural network is the Rusanov scheme which is a $1^{st}$ order conservative scheme for conservation laws.   Equation (\ref{cons-law}) can be discretized by the scheme as 

\begin{equation}
    \label{Rusanov-scheme}
    \frac{U^{n+1}_i - U^n_i}{\Delta t} +\frac{\widehat{f(U)}|^{n}_{i+\frac{1}{2}}-\widehat{f(U)}|^{n}_{i-\frac{1}{2}}}{\Delta x} = 0~,
\end{equation}
where  
\begin{equation}
    \widehat{f(U)}|^{n}_{i+\frac{1}{2}} = \frac{1}{2}\{{f(U^{n}_{i+1})} + {f(U^{n}_{i})}\} + \frac{\alpha}{2}\{U^n_i - U^n_{i+1}\}~,
\end{equation}
$\alpha$ is the largest characteristic speed defined by
\begin{equation}
    \label{alpha}
    \alpha = \max_i \{|v_i|+c_i\}~,
\end{equation}
\\
and $c_i=\sqrt{\gamma \frac{p_i}{\rho_i}}$ is the sound speed. In this N-PINN we fix $\alpha$ to a number determined 
by (\ref{alpha}) using the given initial value, multiplied by a factor, e.g., $2$.
We present three examples predicted by the N-PINN.
\\
\\
\textbf{Example 1. Burgers' Equation}

Consider the Burgers’ equation 

\begin{equation}
\label{Burgers-eq}
    \frac{\partial U}{\partial t} + \frac{\partial (\frac12 U^2)}{\partial x} = 0, x\in[-1, 1], t\in[0, 1],
\end{equation}
\\
with the initial condition given as
\begin{gather*}
    U(0,x) = -\sin(\pi x)~.
\end{gather*}

The viscous version of (\ref{Burgers-eq}) is computed by neural network in \textit{Raissi et al.}~\cite{RAISSI2019686}. The original solution with periodic boundary conditions and specified initial condition is simulated from Chebfun package~\cite{Driscoll2014}. The loss function is minimized by the L-BFGS optimizer (in TensorFlow ~\cite{abadi2016tensorflow, BaydinPearlmutterRadulSiskind17}) during the training process. The total loss is defined in eqn. (\ref{total-loss}). Instead of defining the functional form by auto-differentiation in TensorFlow, it is constructed in ({\ref{loss2}}) with residue $L$ defined to be the left-hand-side of (\ref{Rusanov-scheme}).  

\begin{figure}[H] \centering 
\includegraphics[width=0.8\linewidth]{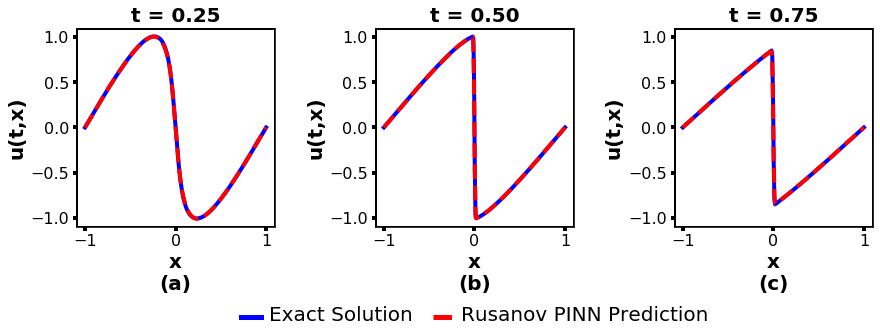}
\caption{Burgers’ equation: (a) to (c) exhibit the prediction and exact solution snapshots from chosen temporal locations, $t = 0.25$, $t = 0.50$, and $t = 0.75$; blue solid line represents exact solution and red dashed line represent prediction results.}
\label{Burgers-Equation-Rusanov}
\end{figure}

Fig.~\ref{Burgers-Equation-Rusanov} summarizes the prediction of solution for Burgers’ equation by N-PINN. To proceed the training process, the initial and boundary conditions ($0$ in this example) are applied in (\ref{loss1}) where the initial condition is given at a set of $256$ data points and boundary conditions are given at $200$ data points ($100$ data points for each boundary). The associated DNN consists of $10$ hidden layers and $40$ neurons in each layer. The resulting error measured by relative $L_2$ norm is $1.328e-2$. Compared with the original PINN, the error is $2$ orders of magnitude higher because Rusanov scheme is only first order. However, the error is lower than that of the framework using the numerical gaussian process to solve time-dependent and nonlinear partial differential equations ~\cite{RAISSI_17,Raissi_gaussian18}.
\\
\\
\textbf{Exampe 2. Lax Problem}

The Lax problem is a Riemann problem for the Euler equations with 
initial value 
\begin{gather*}
    U(0,x) = (0.445, 0.311, 8.928)^T, \text{if}\ 0\leq x\leq 0.5~,\\
    U(0,x) = (0.5, 0.0, 1.4275)^T, \text{if}\ 0.5<x\leq 1.0~,
\end{gather*}
where  $U(t,x) = (\rho, \rho v, E)^T$.

\begin{figure}[H]\centering
\begin{subfigure}[b]{.48\textwidth}
  \centering
  \includegraphics[width=1.0\linewidth]{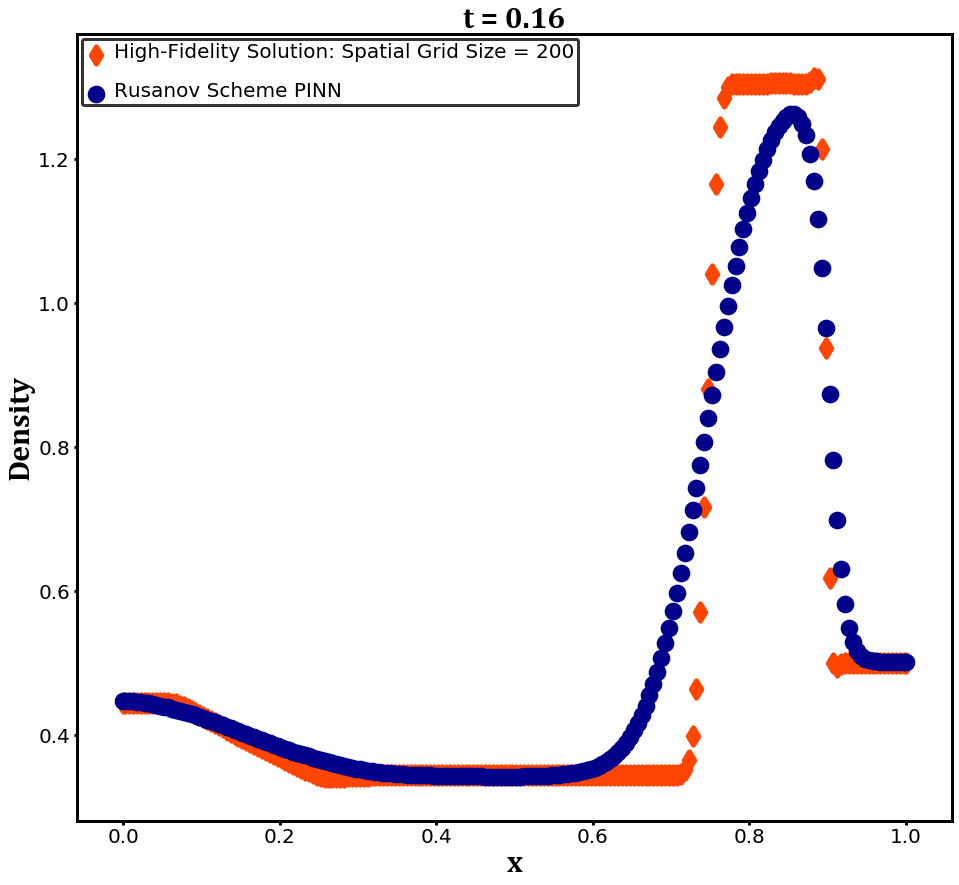}
\end{subfigure}
\begin{subfigure}[b]{.48\textwidth}
  \centering
  \includegraphics[width=1.0\linewidth]{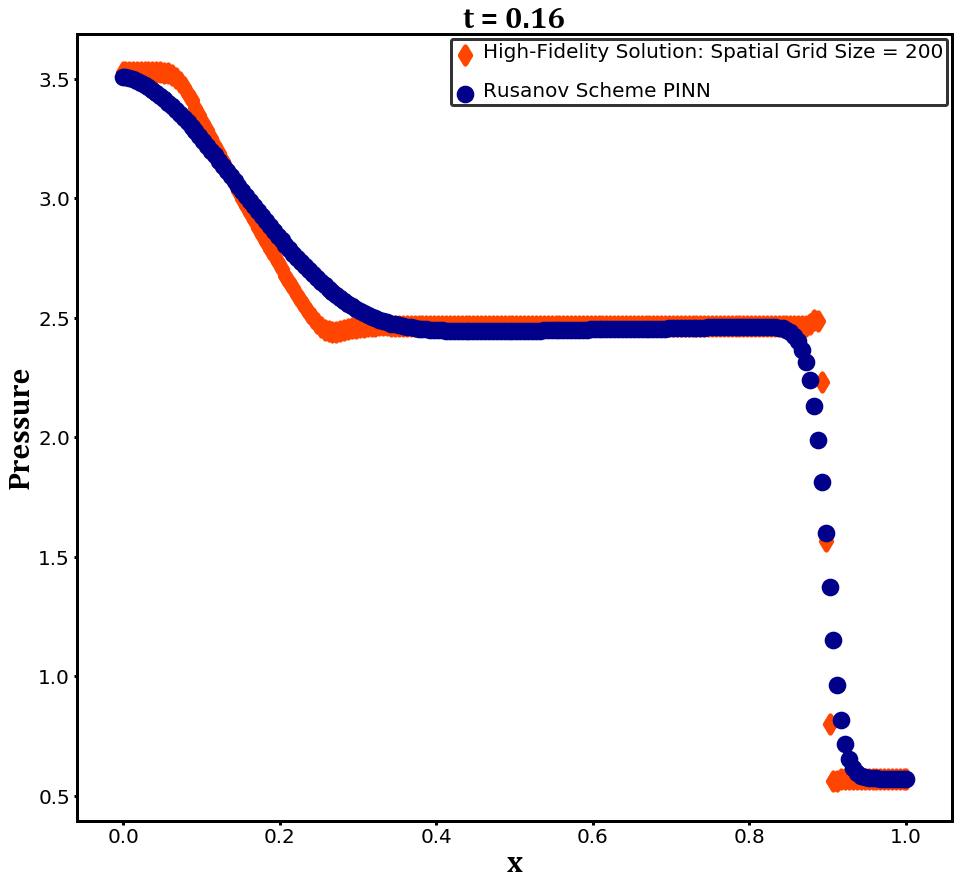}
\end{subfigure}
\begin{subfigure}[b]{.48\textwidth}
  \centering
  \includegraphics[width=1.0\linewidth]{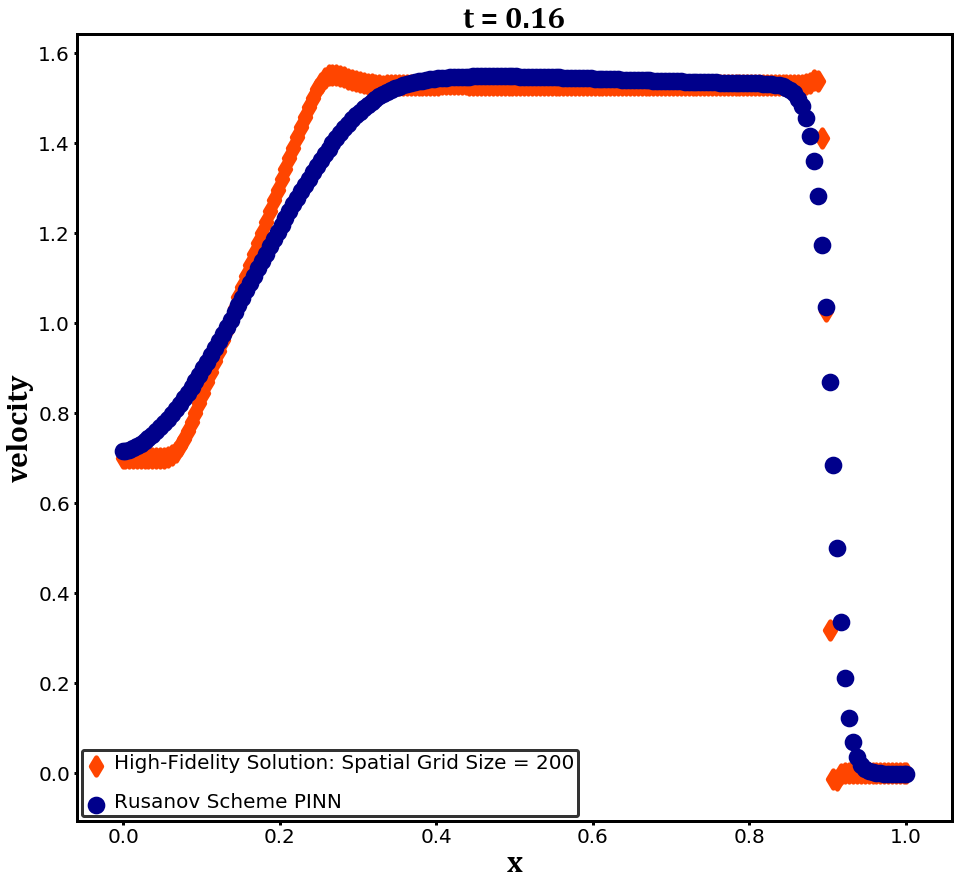}
\end{subfigure}
\caption{Final time prediction $(t=0.16)$ of the Lax problem  by N-PINN based on Rusanov scheme (blue) in comparison with the reference solution (red).}
\label{Rusanov scheme: Final time of Lax problem}
\end{figure}

Fig.~\ref{Rusanov scheme: Final time of Lax problem} presents the comparative results of the Lax problem predicted by N-PINN (based on the Rusanov scheme), and a high resolution reference solution computed by a 3rd order finite volume scheme using non-oscillatory hierarchical reconstruction (HR) limiting ~\cite{LiuShu07b} and partial neighboring cells~\cite{XuLiuShu09} in HR, which will be used for computing reference solutions of the Euler equations throughout the paper.  Even though N-PINN is able to capture the discontinuous solution reasonably well, it has a large numerical diffusion around the discontinuities due to the underlying
first order scheme used in the loss function. Since the number of unknowns in the computational domain (with $200$ cells in $[0,1]$) of the Lax problem is greater than that of the Burgers’ equation, we increase the number of neurons accordingly. The associated DNN consists of 6 hidden layers and 180 neurons in each layer. We have tested multiple fully-connected neural network structures to search for the optimal one for training, and found that the structure applied in this case is the most robust one. The weight for the functional loss ($Err_2$ in (\ref{total-loss})) is 0.7 while the weights for initial and boundary value loss ($Err_1$) are $0.15$ each. The greater weight of $Err_2$ is to enhance the resolution around shocks and contacts.
\\
\\
\textbf{Example 3. Sod Problem}

The Sod problem is a Riemann problem for the Euler equations with initial value 
\begin{gather*}
    U(0,x) = (1.0, 0.0, 2.5)^T, \text{if}\ 0\leq x\leq 0.5~,\\
    U(0,x) = (0.125, 0.0, 0.25)^T, \text{if}\ 0.5<x\leq 1.0~.
\end{gather*}

\begin{figure}[H]\centering
\begin{subfigure}[b]{.48\textwidth}
  \centering
  \includegraphics[width=1.0\linewidth]{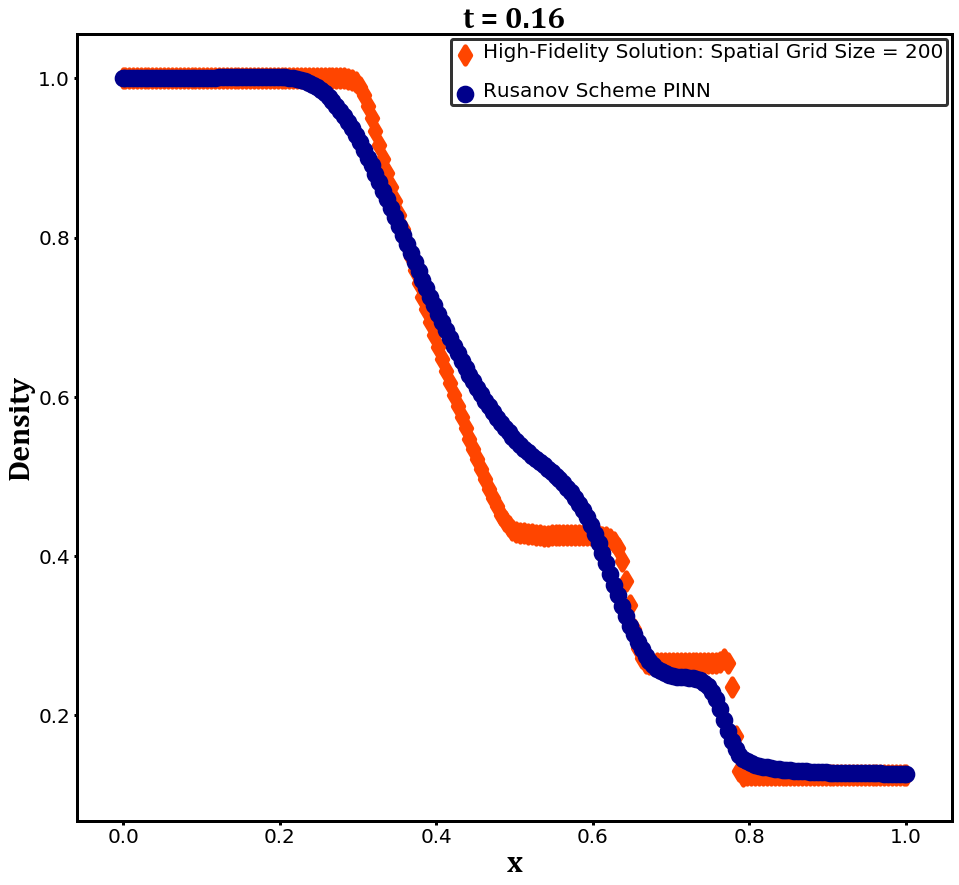}
\end{subfigure}
\begin{subfigure}[b]{.48\textwidth}
  \centering
  \includegraphics[width=1.0\linewidth]{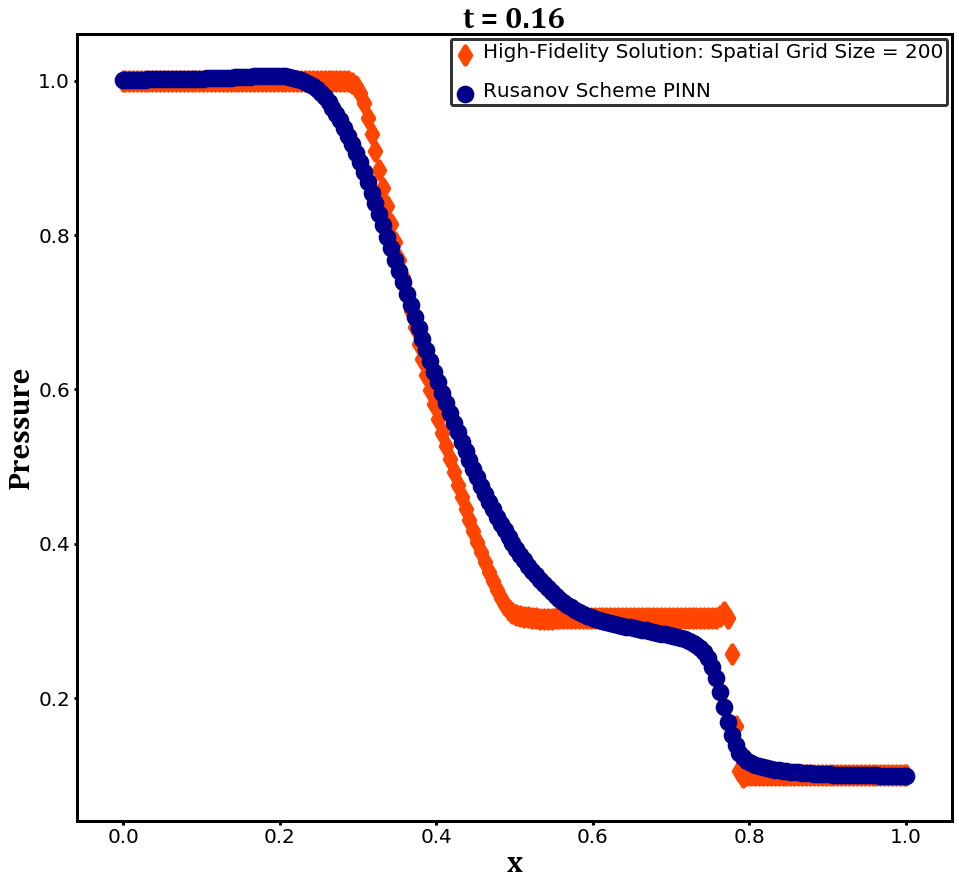}
\end{subfigure}
\begin{subfigure}[b]{.48\textwidth}
  \centering
  \includegraphics[width=1.0\linewidth]{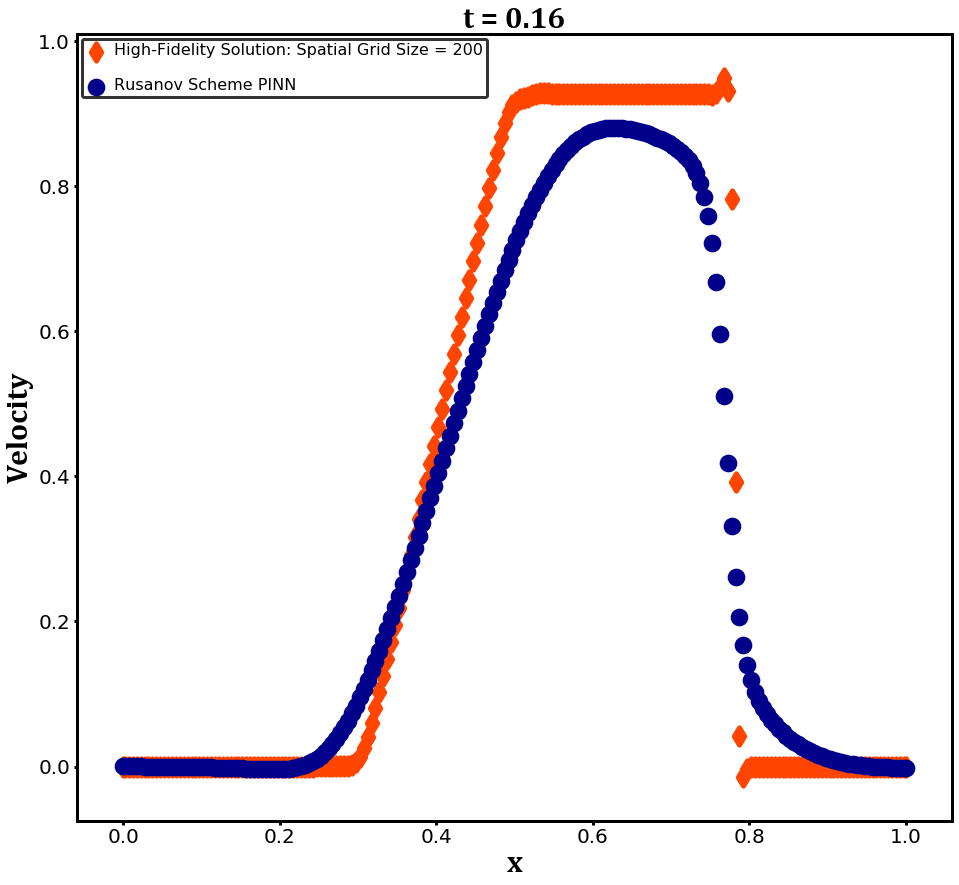}
\end{subfigure}
\caption{Final time prediction $(t=0.16)$ of the Sod problem by N-PINN based on Rusanov scheme (blue) in comparison with the reference solution (red).}
\label{Rusanov scheme: Final time of Sod problem}
\end{figure}

Fig.~\ref{Rusanov scheme: Final time of Sod problem} shows the comparison of the final-time solution of N-PINN (based on the Rusanov scheme) and the reference solution. The numerical diffusion is large due to the first order underlying scheme used in the loss function. The setup of the associated DNN as well as the weights for parts of the loss function are the same as those in the Lax problem.

\subsection{Physics Informed Neural Network with Leapfrog Scheme and Artificial Diffusion }
\begin{figure}[H]\centering
\begin{subfigure}[b]{.48\textwidth}
  \centering
  \includegraphics[width=1.0\linewidth]{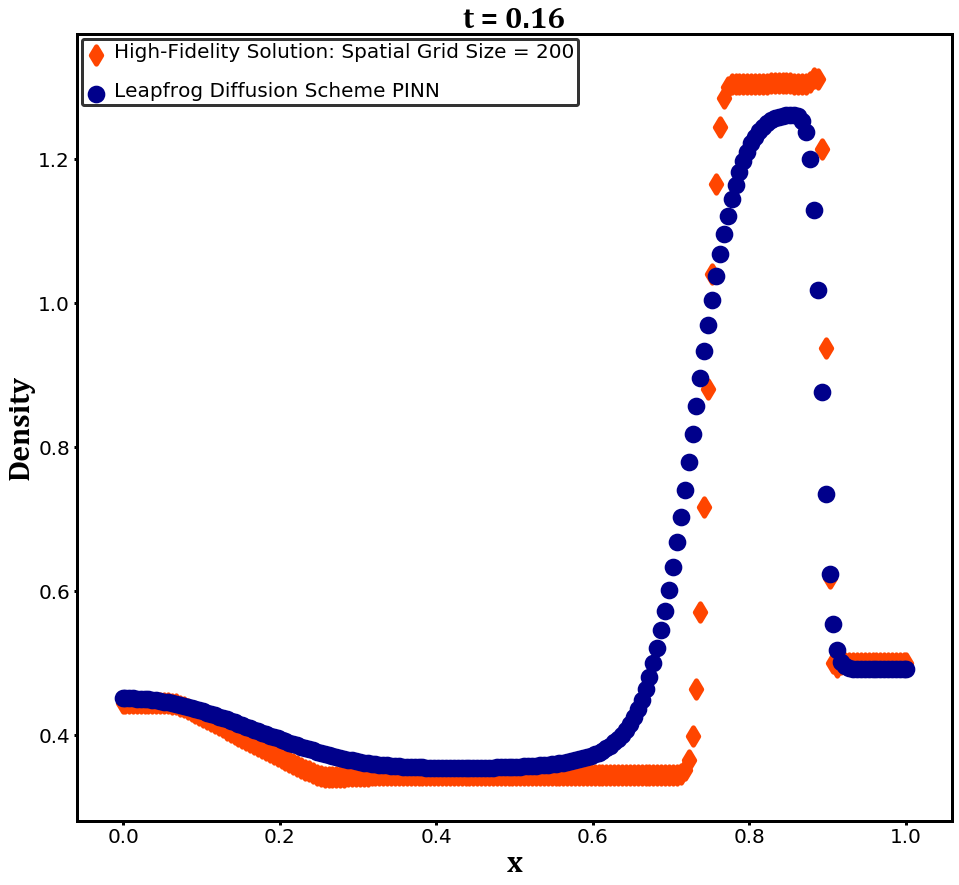}
\end{subfigure}
\begin{subfigure}[b]{.48\textwidth}
  \centering
  \includegraphics[width=1.0\linewidth]{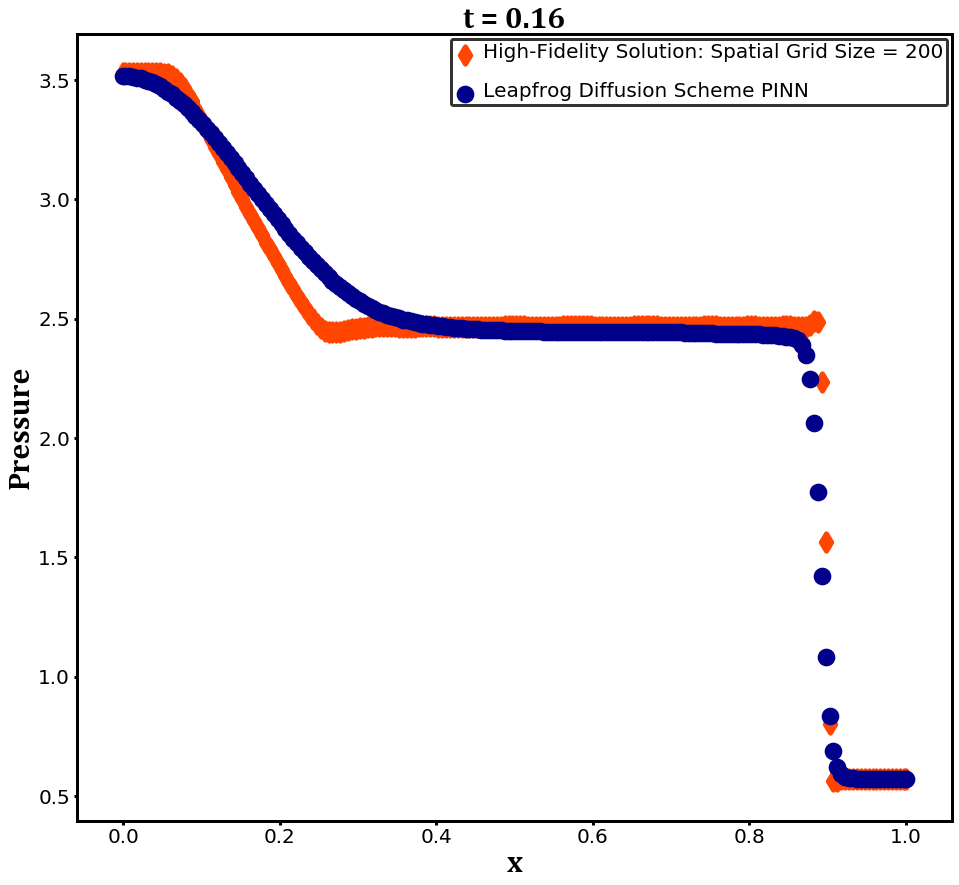}
\end{subfigure}
\begin{subfigure}[b]{.48\textwidth}
  \centering
  \includegraphics[width=1.0\linewidth]{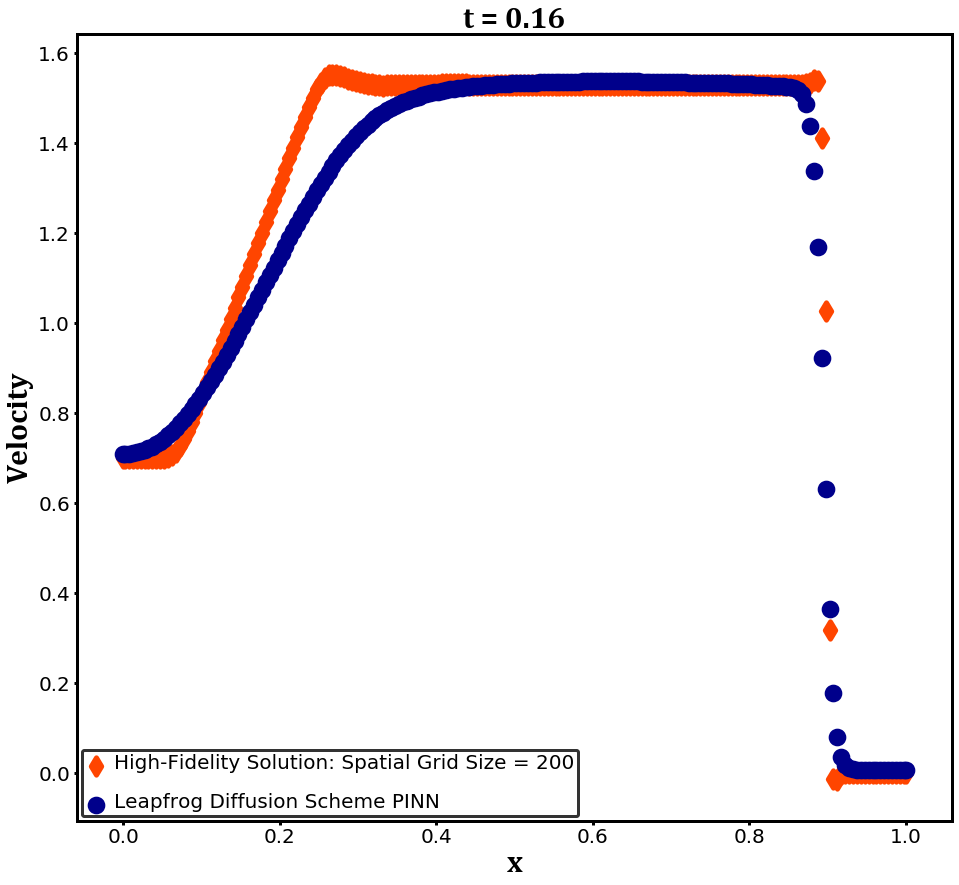}
\end{subfigure}
\caption{Final time prediction $(t=0.16)$ of the Lax problem by N-PINN based on Leapfrog and diffusion (blue) in comparison with the reference solution (red).}
\label{Leapfrog Diffusion scheme: Final time of Lax problem}
\end{figure}

Leapfrog scheme is a commonly used second order scheme for solving first order hyperbolic PDEs. However,  artificial diffusion terms of size $O(\Delta x)$ need to be added to the equations in order to damp the numerical oscillations around discontinuities generated by a linear second order scheme. The central discretization of the diffusion terms 
yields the following (unstable) scheme
\\
\begin{equation}
\label{leapfrog-diffusion}
    \frac{U^{n+1}_i - U^{n-1}_i}{2\Delta t} +\frac{{f(U)}|^{n}_{i+1}-{f(U)}|^{n}_{i-1}}{2\Delta x} - \alpha\frac{U^n_{i+1} - 2\cdot U^n_i + U^n_{i-1}}{\Delta x^2}=0~,
\end{equation}
\\
where $\alpha = \Delta x$.
Nevertheless, the residue of the scheme, {\it i.e.}, the left-hand-side of (\ref{leapfrog-diffusion}) can be used as residue $L$ in (\ref{loss2}), a nice property of the neural network. The resulted N-PINN (N-PINN based on Leapfrog and diffusion) performs similarly to that based on the Rusanov scheme. Fig.~\ref{Leapfrog Diffusion scheme: Final time of Lax problem} displays the final-time predicted solution of N-PINN based on Leafrog and diffusion compared to the reference solution, which has slightly higher resolution than that computed by N-PINN based on the Rusanov scheme. We use the same neural network structures and weights for the loss function during the training procedure as those for N-PINN based on the Rusanov scheme. 

\section{Numerical Physics Informed Neural Network (N-PINN) with Additional Training Data}
\label{Sec: partial-data}
\subsection{Methodology}
N-PINN based on first order schemes has relatively poor resolution around shocks and contacts because of the large numerical diffusion of the underlying scheme. We have tried using the residue of the second order MUSCL scheme~\cite{vanLeer79} in the loss function of N-PINN but the minimization procedure has difficulty in reducing it to the desired level. In previous studies \cite{RAISSI2019686,PINNrev21,MAO20}, sparse training data from the exact solution is used in PINN to provide more information to the loss functional and enhance the accuracy of prediction. We provide in the first part of the loss functional an additional (\ref{loss1}) $21.9\%$ of total data points ($7437$ out of  $202\times 168$ space-time grid points for the Lax problem; $3003$ out of $202\times 79$ grid points for the Sod problem) allocated around shocks and contacts, where at the same time we decrease the weight (e.g., from $1$ to $0.1$) of local residue error $L$ in (\ref{loss2}). 

\subsection{Results and Discussions}
\begin{figure}[H]\centering
\begin{subfigure}[b]{.48\textwidth}
  \centering
  \includegraphics[width=1.0\linewidth]{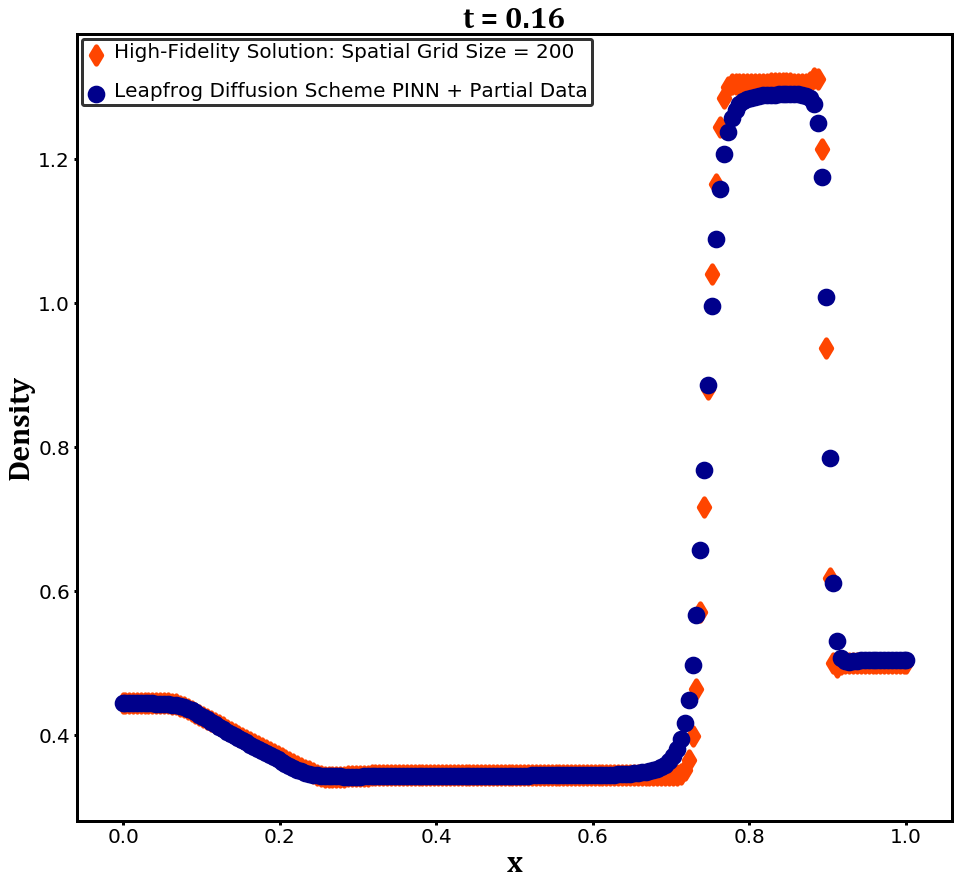}
\end{subfigure}
\begin{subfigure}[b]{.48\textwidth}
  \centering
  \includegraphics[width=1.0\linewidth]{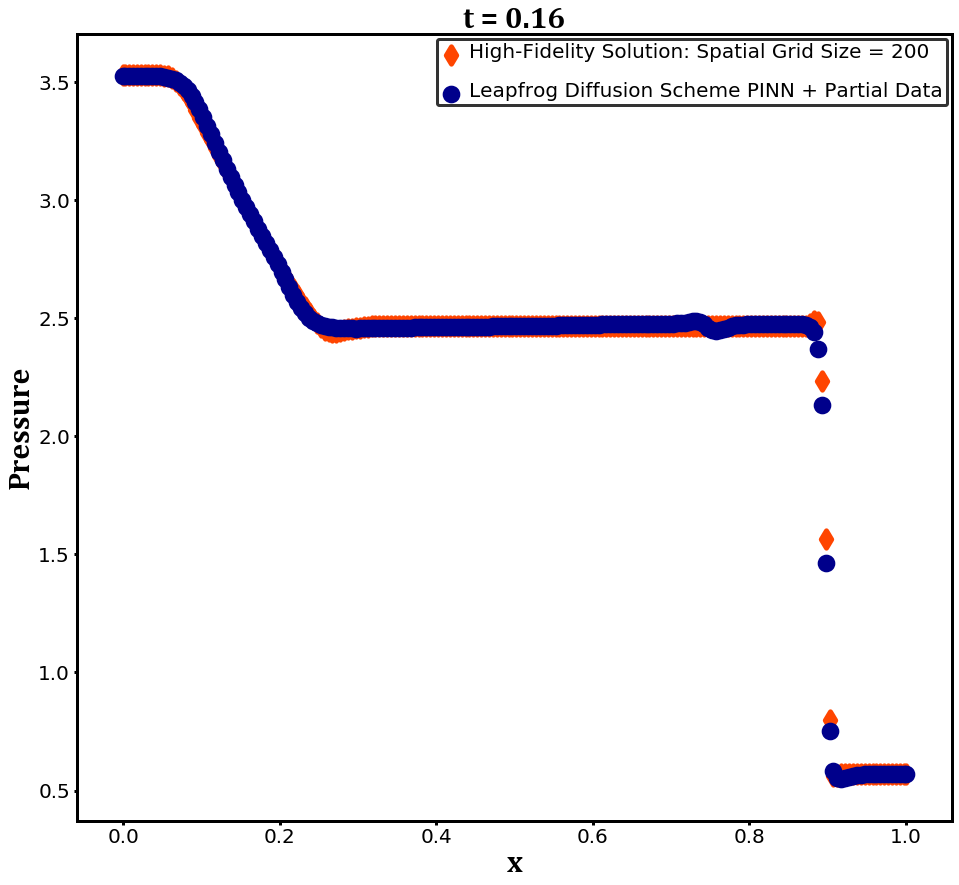}
  
\end{subfigure}
\begin{subfigure}[b]{.48\textwidth}
  \centering
  \includegraphics[width=1.0\linewidth]{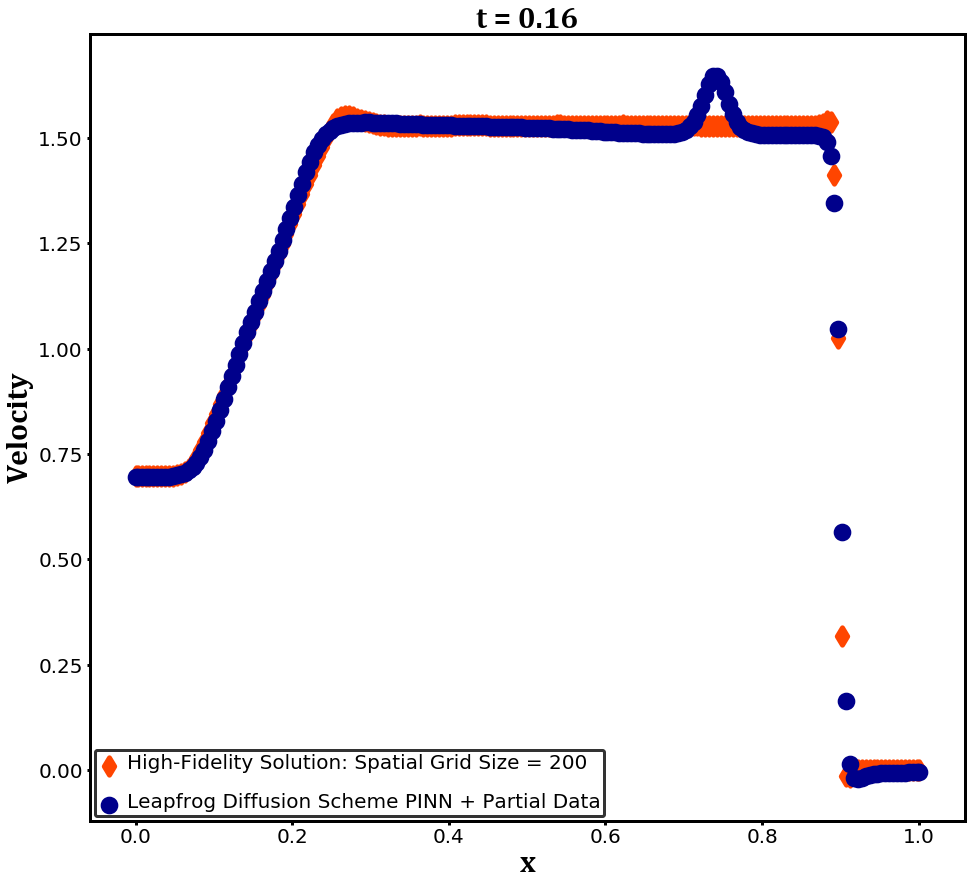}
\end{subfigure}
\caption{Final time prediction $(t=0.16)$ of the Lax problem  by N-PINN based on Leapfrog and diffusion with additional data (blue) in comparison with the reference solution (red).}
\label{Leapfrog Diffusion scheme with dat: Final time of Lax problem}
\end{figure}

\begin{figure}[H]\centering
\begin{subfigure}[b]{.48\textwidth}
  \centering
  \includegraphics[width=1.0\linewidth]{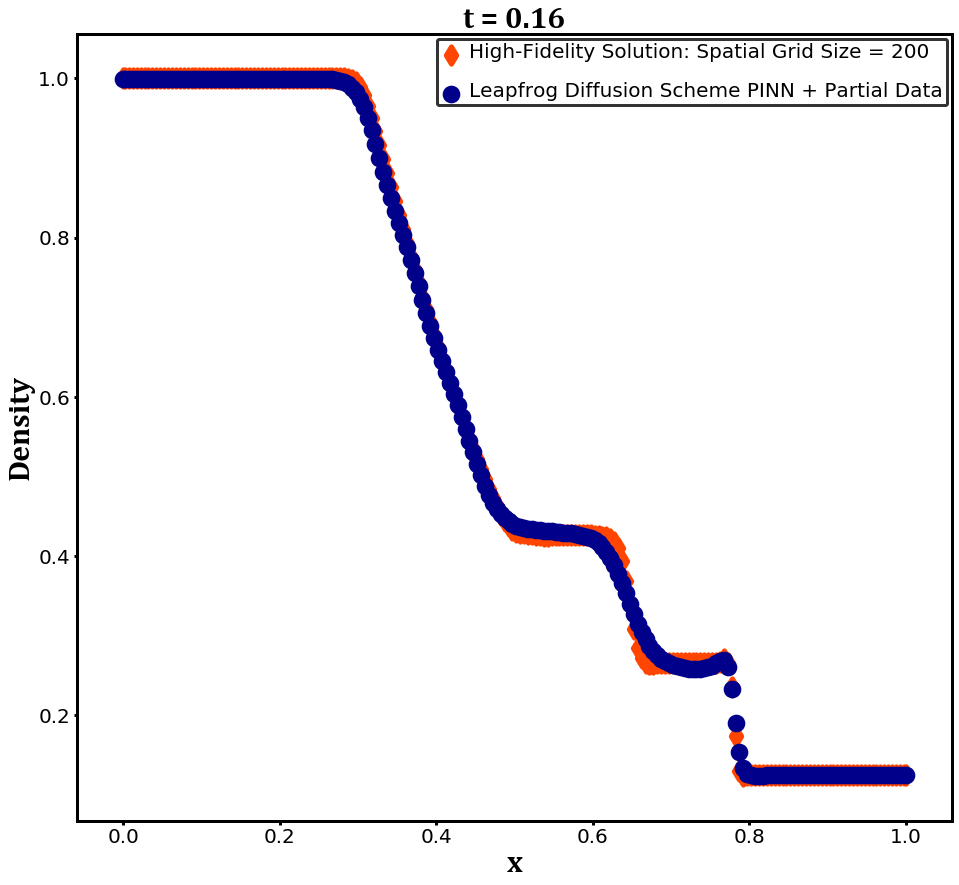}
  
\end{subfigure}
\begin{subfigure}[b]{.48\textwidth}
  \centering
  \includegraphics[width=1.0\linewidth]{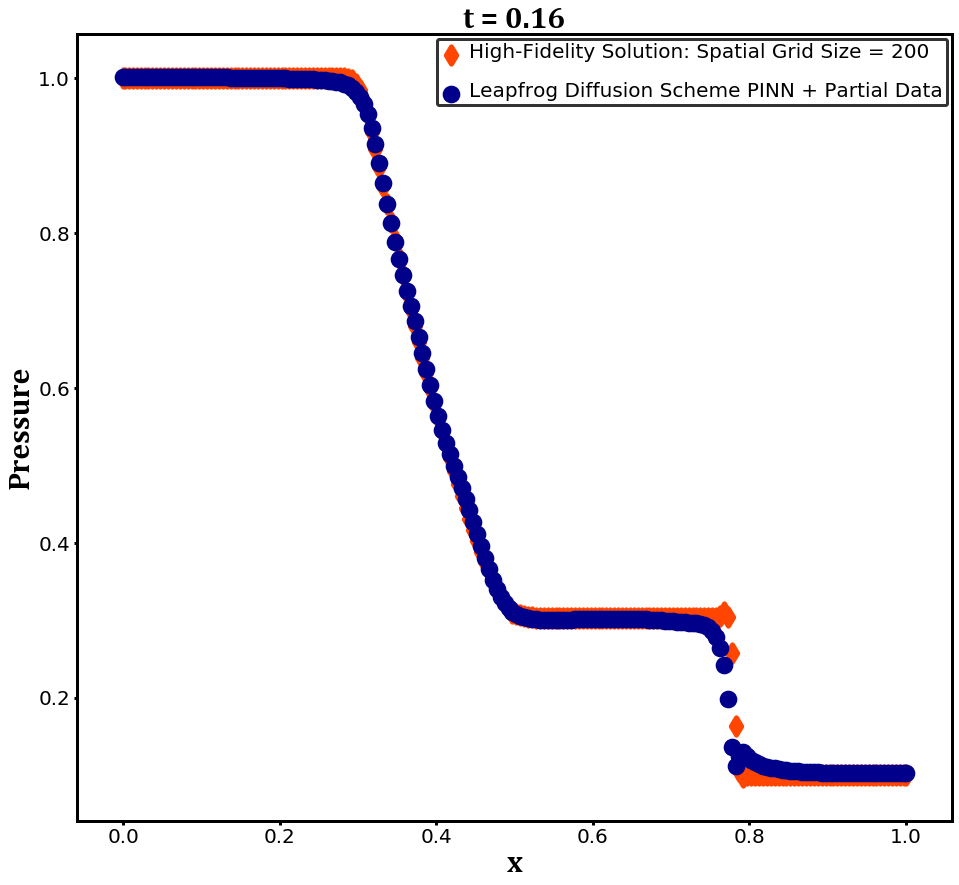}
\end{subfigure}
\begin{subfigure}[b]{.48\textwidth}
  \centering
  \includegraphics[width=1.0\linewidth]{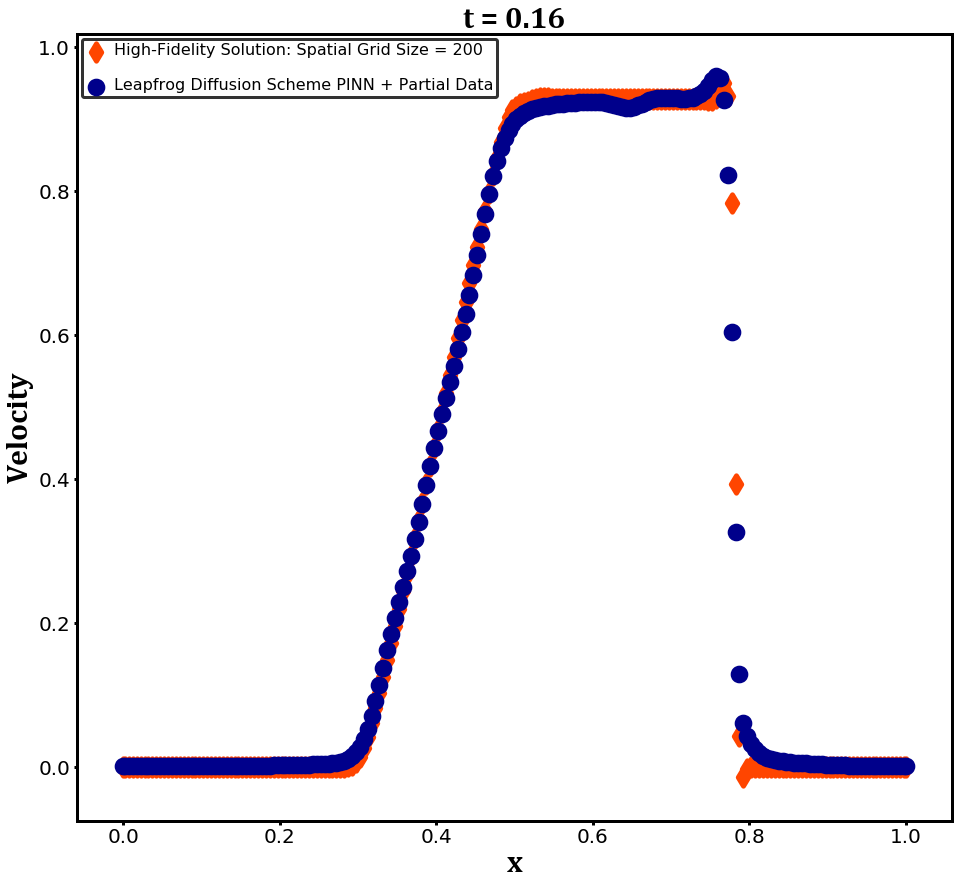}
  
\end{subfigure}
\caption{Final time prediction $(t=0.16)$ of the Sod problem by N-PINN based on Leapfrog and diffusion with additional data (blue) in comparison with the reference solution (red).}
\label{Leapfrog Diffusion scheme with dat: Final time of Sod problem}
\end{figure}

We present the results of N-PINN based on Leapfrog and diffusion with additional training data for the Lax and Sod problem. Because of the additional training data, the number of neurons in the associated DNN increases to $240$ each layer for the Lax problem while the number of layers remains to be $6$. For the Sod problem, the associated DNN remains the same as before ($6$ layers and $180$ neurons in each layer.) Meanwhile, the weights for the loss function are adjusted. We assign the weight $0.1$ for the loss from initial value, $0.1$ for the loss from boundary conditions, $0.5$ for the loss from residue error of the underlying scheme,  and $0.3$ for the loss from additional training data.
Fig.~\ref{Leapfrog Diffusion scheme with dat: Final time of Lax problem} shows the final-time profiles of density, pressure and velocity of the N-PINN prediction of the Lax problem and reference solution. Except for the numerical artifact around the contact for the velocity field, the overall resolution is much improved from Fig.~\ref{Leapfrog Diffusion scheme: Final time of Lax problem} which has similar setup in N-PINN but without the additional training data.  Fig.~\ref{Leapfrog Diffusion scheme with dat: Final time of Sod problem},
shows the final-time profiles of density, pressure and velocity of the N-PINN prediction of the Sod problem and reference solution. The quality of the prediction is much better than those in Fig.~\ref{Rusanov scheme: Final time of Sod problem} which has similar setup in N-PINN but without the additional training data. 


\section{2-Coarse-Grid Neural Network }
It is an efficient strategy to use a low cost numerical solution of a PDE as the input of a neural network to predict a more accurate solution of the PDE. In \cite{MFNN_liu19}, the numerical solution in a coarse mesh is used in a neural network to predict the solution in a fine mesh. In \cite{Nguyen2020ASP},
the solution of a second order wave equation and its gradient computed by a first order scheme is used to predict the solution computed by a higher order scheme. In this study, we introduce a novel neural network, the 2-Coarse-Grid neural network (2CGNN) to accurately compute solutions of systems of conservation laws which may contain shocks and contacts. Instead of taking a spatiotemporal point $(x,t)$ as input, 2CGNN first needs two different coarse grid solutions to be computed. The input of the neural network is designed to incorporate a small domain of dependence of the solution at certain space-time location as well as an approaching sequence of approximations. Once trained, the neural network is able to predict a fine discontinuous solution efficiently based on some low-cost computations.


\subsection{Input, Output and Loss Functional}
Consider the scalar conservation law (\ref{cons-law}) on an interval $[a,b]$ which is partitioned with the coarsest uniform grid $a=x_0<x_1<\dots <x_M=b$ having spatial grid size $\Delta x=x_1-x_0$ and time step size $\Delta t$. Refine the grid to obtain a finer uniform grid with spatial grid size $\frac12 \Delta x$ and time step size $\frac12 \Delta t$. Let $L$ be
a low-cost scheme used to compute (\ref{cons-law}) on both grids.
Given a grid point $x_{i'}$ at the time $t_{n'}$ (on the coarsest uniform grid) where the solution is to be predicted by DNN, we choose the coarsest grid solution (computed by $L$) at $3$ points $x_{i'-1}$, $x_{i'}$ and $x_{i'+1}$ at the time level $t_{n'-1}$ and also at the point $(x_{i'}, t_{n'})$ as the first part of the input, and the finer grid solution (also computed by $L$) at the same spatiotemporal locations as the second part of input. Note that, the chosen $4$ spatiotemporal locations of either grid enclose a (space-time) domain of dependence of the exact solution at $(x_{i'}, t_{n'})$ (with $\Delta t$ satisfying the CFL restriction.)
Since the two parts of the input solution have different levels of approximation to the exact or reference solution, DNN utilizes the information to extrapolate a prediction of the exact or reference solution.




Denote the first part of the input as
$$u^{n'-1}_{i'-1}, u^{n'-1}_{i'}, u^{n'-1}_{i'+1}, u^{n'}_{i'}~,$$
and 
the second part of input as
$$u^{n''-2}_{i''-2}, u^{n''-2}_{i''}, u^{n''-2}_{i''+2}, u^{n''}_{i''}~.$$
Note that the space-time index $(i',n')$ in the coarsest grid refers to the same location as $(i'',n'')$ does in the finer grid, $(i'-1, n'-1)$ refers to the same location as $(i''-2, n''-2)$ does, and so on.

\begin{figure}[H]\centering
\begin{subfigure}[b]{.8\textwidth}
  \centering
  \includegraphics[width=1.0\linewidth]{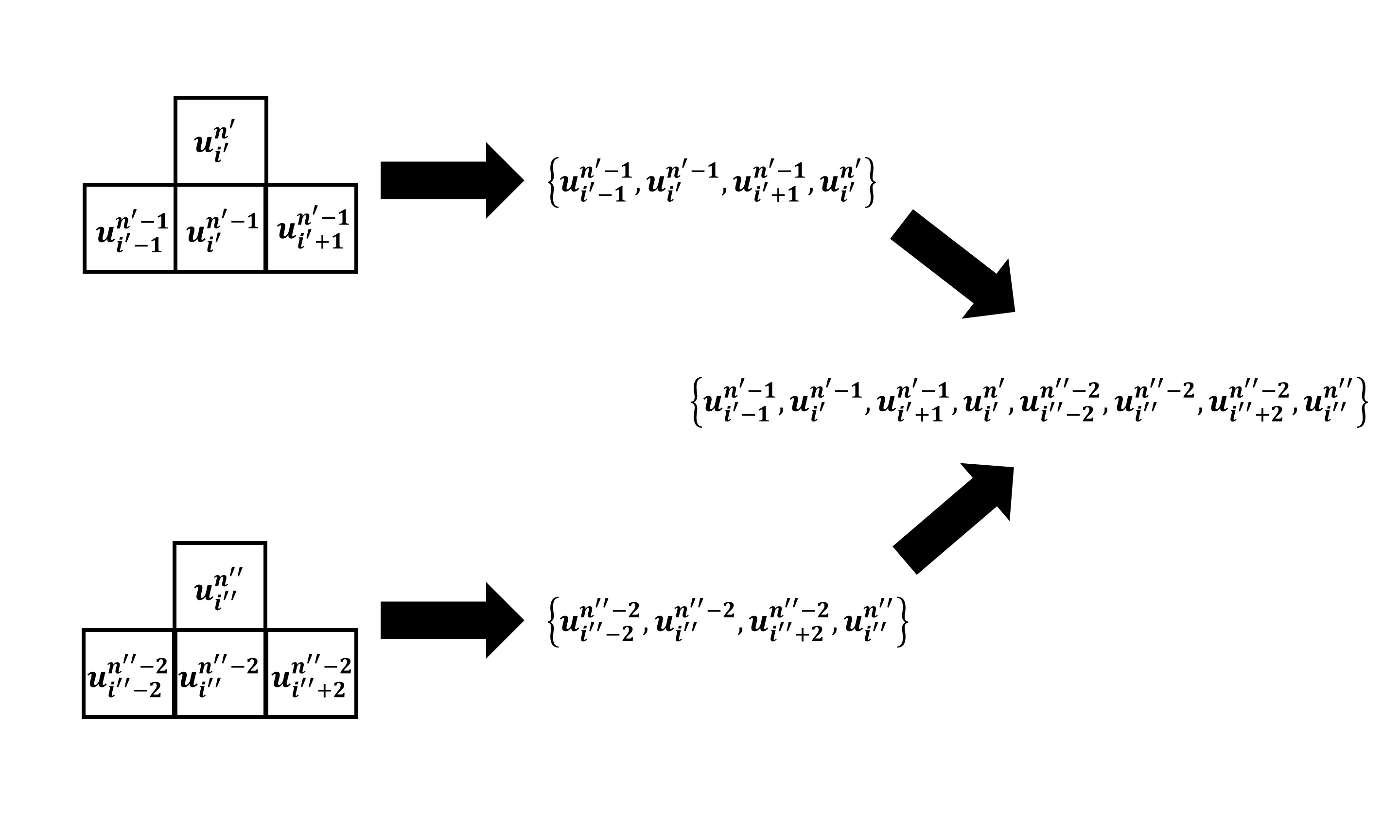}
  \caption{}
\end{subfigure}
\hfill
\begin{subfigure}[b]{.8\textwidth}
  \centering
  \includegraphics[width=1.0\linewidth]{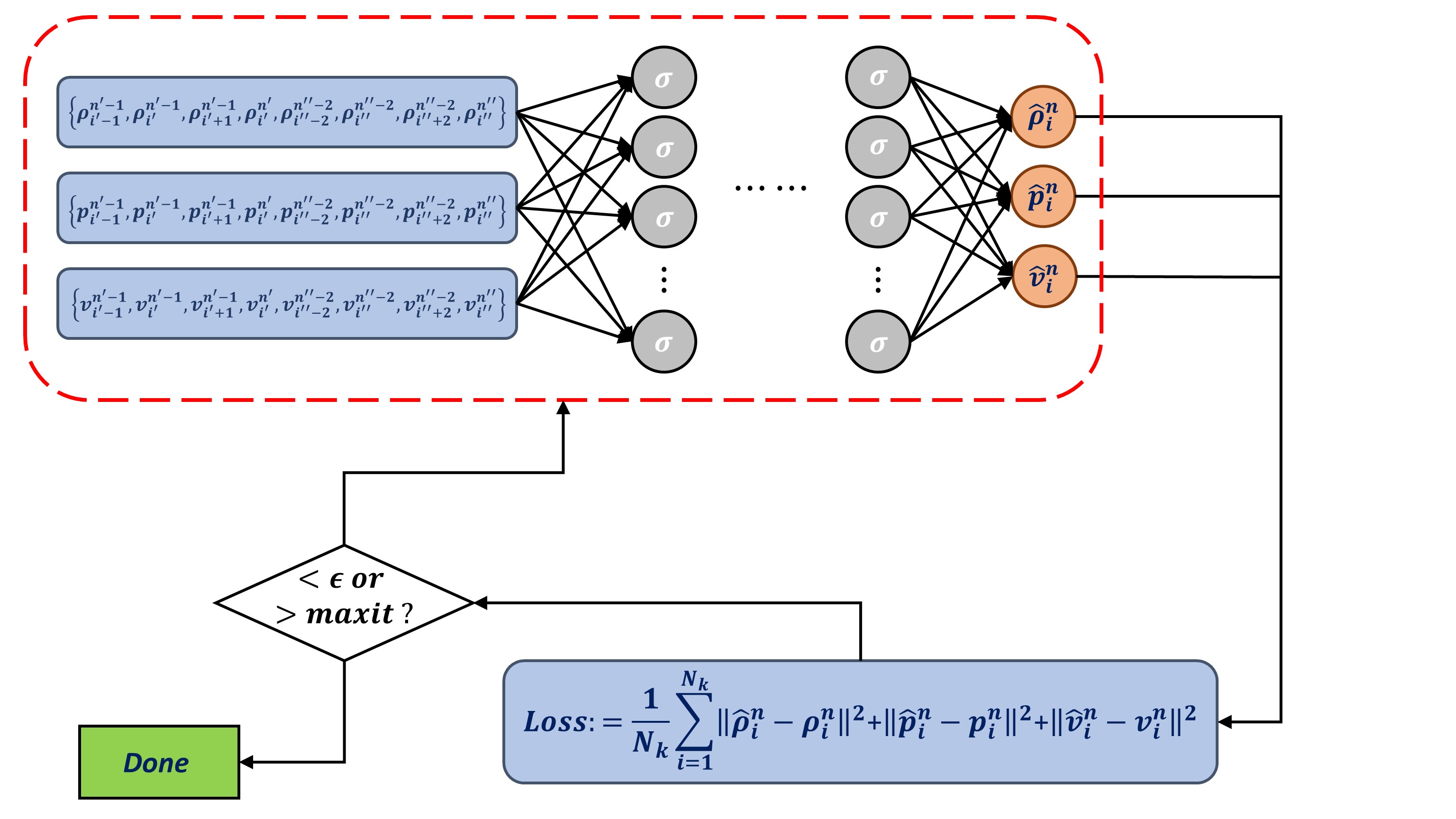}
  \caption{}
\end{subfigure}
\caption{Schematic of 2-Coarse-Grid neural network: (a) Procedure of transferring 2 coarse grids into input format; (b) Training procedure of 2-Coarse-Grid neural network.}
\label{2-Coarse-Grid NN}
\end{figure}

Suppose we are interested in the predicted solution at $(x,t)$ which is referred to as $(i',n')$ in the coarsest grid,
the input of 2CGNN is 
\begin{equation}
\label{standard-input}
    \{u^{n'-1}_{i'-1}, u^{n'-1}_{i'}, u^{n'-1}_{i'+1}, u^{n'}_{i'}, u^{n''-2}_{i''-2}, u^{n''-2}_{i''}, u^{n''-2}_{i''+2}, u^{n''}_{i''}\}~,
\end{equation}
called ``input of $u$'', and the corresponding output of 2CGNN
is the predicted solution at $(x,t)$.
See Fig.~\ref{2-Coarse-Grid NN} for an illustration of 2CGNN.
For the Euler system, the input and output of 2CGNN are made up of corresponding ones for each prime variable. For example, the input can be
the vector 
$$\{ {\rm input}\; {\rm of}\; \rho,\;
{\rm input}\; {\rm of}\; v,\;
{\rm input}\; {\rm of}\; p\}
$$
with $8\times 3=24$ elements.
And the corresponding output will be 
$$\{\rho, v, p\}
$$
with $3$ elements.

The loss function measures the difference between the output and the reference solution corresponding to the input, and is defined as follows.
$$
\begin{array}{ccc}
{\rm Loss}&=&\sum_{k}\|({\rm output} \; {\rm corresponding}\; {\rm to}\; k^{th}\; {\rm set}\; {\rm of}\; {\rm input})- \\
&& ({\rm reference}\; {\rm solution}\; {\rm corresponding}\; {\rm to}\; k^{th}\; {\rm set}\; {\rm of}\; {\rm input}) \|_2^2~,
\end{array}
$$

where $\|\cdot\|_2$ is the $2$-norm, and the summation goes through every set of input in the training data.

\subsection{Generation of Input and Training Data}
The associated DNN for 2CGNN consists of $8$ hidden layers and each layer has $300$ neurons.
During the training process, DNN minimizes the difference between outputs of the neural network and a reference solution by using L-BFGS optimizer in TensorFLow (with the number of iterations of optimization procedure under $50000$.) After the training is done, the neural network is used to predict a solution (different from the training data) given an input computed by the same low-cost scheme(s) and grids which are used to compute inputs of the training data.  We try to use the trained neural network to predict final solutions for $5$ initial values of the Euler system which include the original initial value of the Lax (or Sod) problem, $\pm3\%$ and $\pm5\%$ perturbations of the initial value. 

In order to generate the training data, we use a first order scheme on
the coarsest uniform grid ($50$ cells) and the finer uniform grid ($100$ cells) to compute the input data from several initial values of the Euler system which include $\pm2\%$, $\pm4\%$, $\pm6\%$, $\pm8\%$ and $\pm10\%$ perturbations of the initial value of the Lax (or Sod) problem.
The high resolution reference solutions of the training data are computed on a uniform grid
with $200$ cells (note that solution values at grid points need to be interpolated from cell averages.) The first order scheme used for generating the inputs is a stable version of the Leapfrog and diffusion scheme (\ref{leapfrog-diffusion-stabilized}), the Rusanov scheme or the Leapfrog and diffusion splitting scheme (\ref{leapfrog-diffusion-splitting}). The time step size for the first order schemes is fixed to a constant, e.g., $\Delta t$ for the $50$-cell grid and $\frac12 \Delta t$ for the $100$-cell grid, bounded by the CFL restriction based on the largest characteristic speed estimated throughout the computational domain. Using different low-cost schemes in 2CGNN helps demonstrate the robustness of the method. Prediction results based on schemes (\ref{leapfrog-diffusion-stabilized}) and (\ref{leapfrog-diffusion-splitting}) are shown in Appendix A and B respectively. 

\begin{equation}
\label{leapfrog-diffusion-stabilized}
    \frac{U^{n+1}_i - U^{n-1}_i}{2\Delta t} +\frac{{f(U)}|^{n}_{i+1}-{f(U)}|^{n}_{i-1}}{2\Delta x} - \alpha\frac{U^{n-1}_{i+1} - 2\cdot U^{n-1}_i + U^{n-1}_{i-1}}{\Delta x^2}=0~,
\end{equation}
\\
where $\alpha = \Delta x$.

\begin{equation}
\label{leapfrog-diffusion-splitting}
\left \{
\begin{array}{l}
    \frac{\tilde{U}_i - U^{n-1}_i}{2\Delta t} +\frac{{f(U)}|^{n}_{i+1}-{f(U)}|^{n}_{i-1}}{2\Delta x} =0~, \\
    \frac{U^{n+1}_i - \tilde{U}_i}{\Delta t}- \alpha\frac{\tilde{U}_{i+1} - 2\cdot \tilde{U}_i + \tilde{U}_{i-1}}{\Delta x^2}=0~,
\end{array}
\right .
\end{equation}
\\
where $\alpha = \Delta x$.


\subsection{Results and Discussions for The Lax Problem}
We first use the Rusanov Scheme to compute inputs of 2CGNN. The prediction of the final-time solution of the Lax problem is shown in Fig.~\ref{2CGNN Rusanov input: Final time of lax problem, original}. 
It's clear that 2CGNN captures the shock and contact very sharp, and the prediction is comparable or even better than the reference solutions used in training (and as the ``exact'' solution in the figure.)
Fig.~\ref{2CGNN Rusanov input: Final time of lax problem, +3} (\ref{2CGNN Rusanov input: Final time of lax problem, -3}) shows the prediction of the final-time solution of the Euler system with its initial value being $+3\%$  ( $-3\%$)
perturbation of that of the Lax problem.
Fig.~\ref{2CGNN Rusanov input: Final time of lax problem, +5} (\ref{2CGNN Rusanov input: Final time of lax problem, -5}) shows the prediction of the final-time solution of the Euler system with its initial value being $+5\%$  ( $-5\%$)
perturbation of that of the Lax problem.
 Since the coarsest grid for input has $50$ grid cells, the spatial grid for the predicted solution also has $50$ grid cells. The predictions do not depict smeared solutions like their low-cost input solutions, and are significant improvements over those by N-PINN. Note that the cases for predictions are not included in the training data. 
 
 

\begin{figure}[H]\centering
\begin{subfigure}[b]{.48\textwidth}
  \centering
  \includegraphics[width=1.0\linewidth]{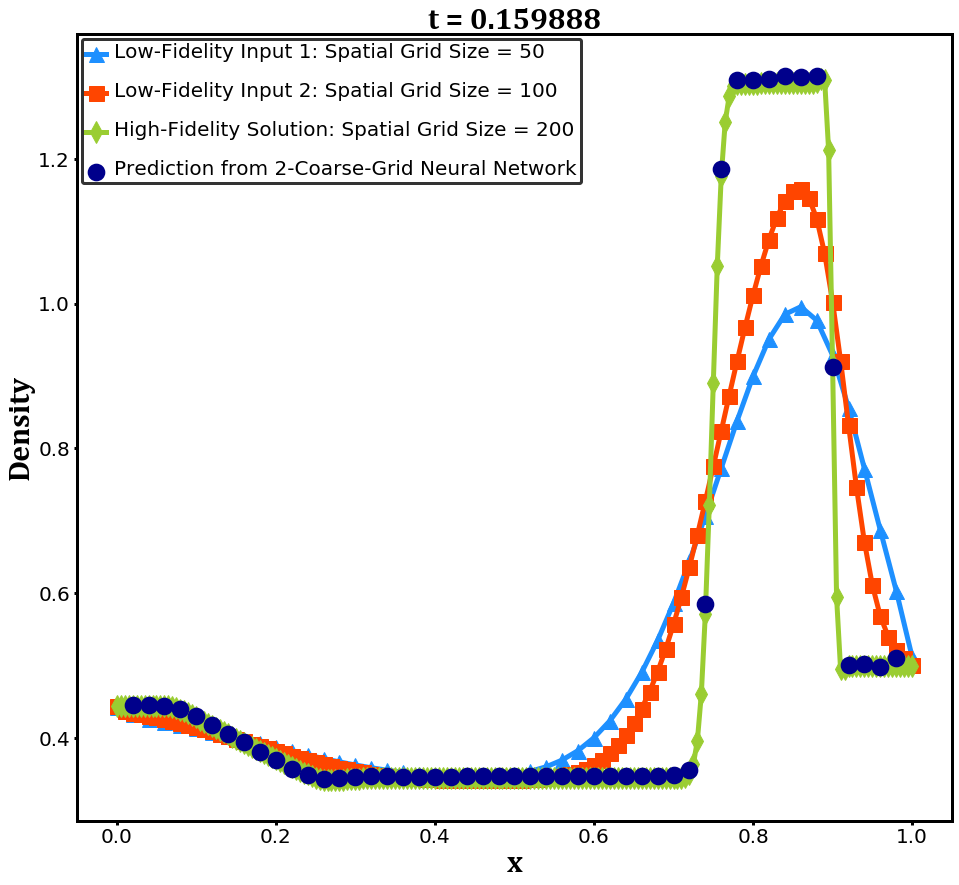}
\end{subfigure}
\begin{subfigure}[b]{.48\textwidth}
  \centering
  \includegraphics[width=1.0\linewidth]{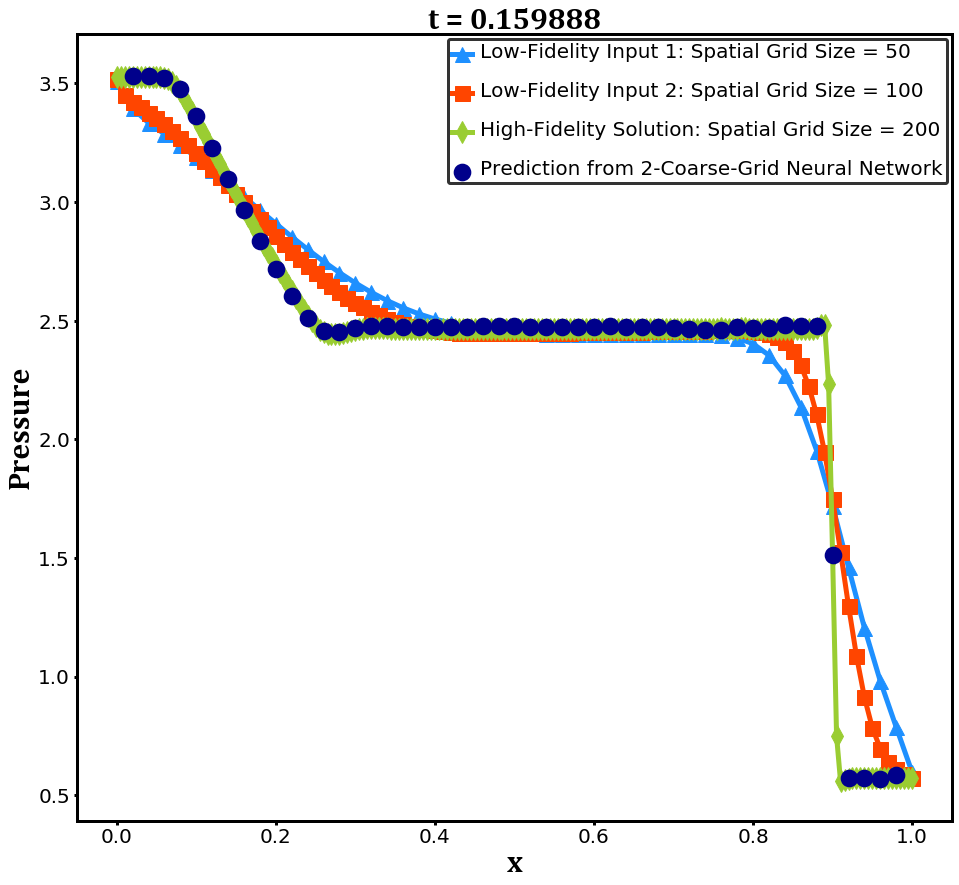}
\end{subfigure}
\begin{subfigure}[b]{.48\textwidth}
  \centering
  \includegraphics[width=1.0\linewidth]{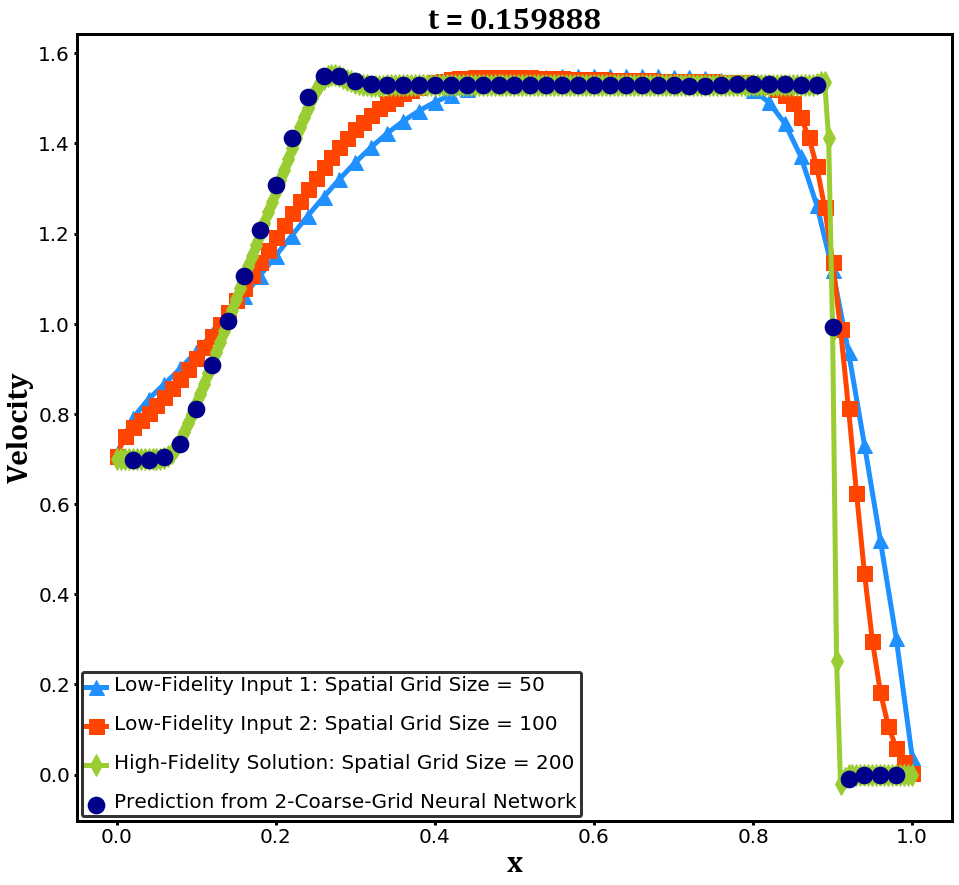}
\end{subfigure}
\caption{2CGNN prediction of final-time $(t=0.159888)$ solution of {\bf Lax problem}  (dark blue), low-fidelity input solutions (blue and red) by Rusanov scheme on $2$ different grids (with $50$ and $100$ cells resp.), and ``exact'' (reference) solution (green).}
\label{2CGNN Rusanov input: Final time of lax problem, original}
\end{figure}

\begin{figure}[H]\centering
\begin{subfigure}[b]{.48\textwidth}
  \centering
  \includegraphics[width=1.0\linewidth]{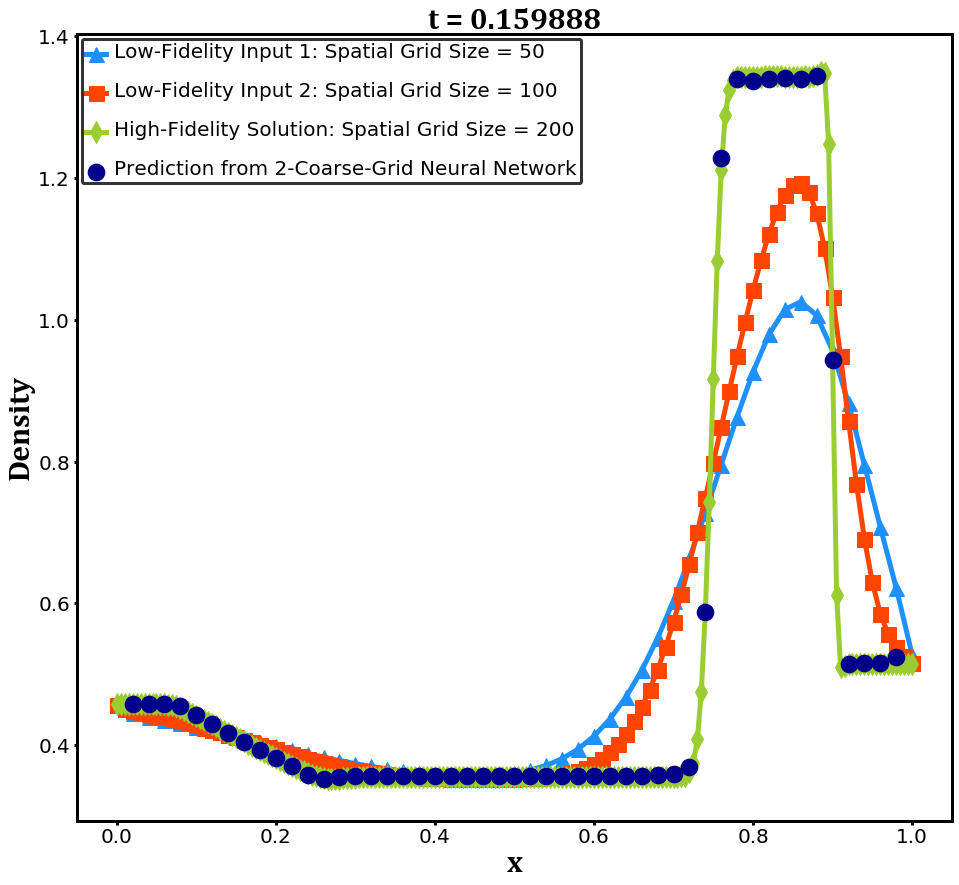}
\end{subfigure}
\begin{subfigure}[b]{.48\textwidth}
  \centering
  \includegraphics[width=1.0\linewidth]{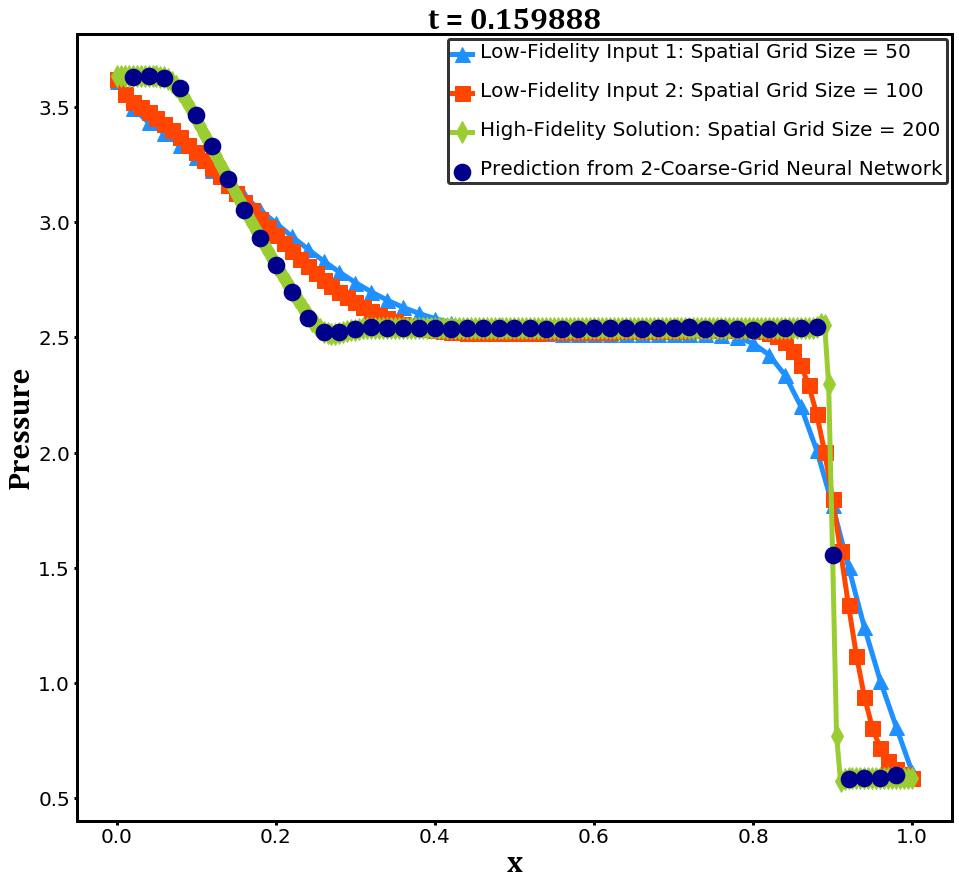}
\end{subfigure}
\begin{subfigure}[b]{.48\textwidth}
  \centering
  \includegraphics[width=1.0\linewidth]{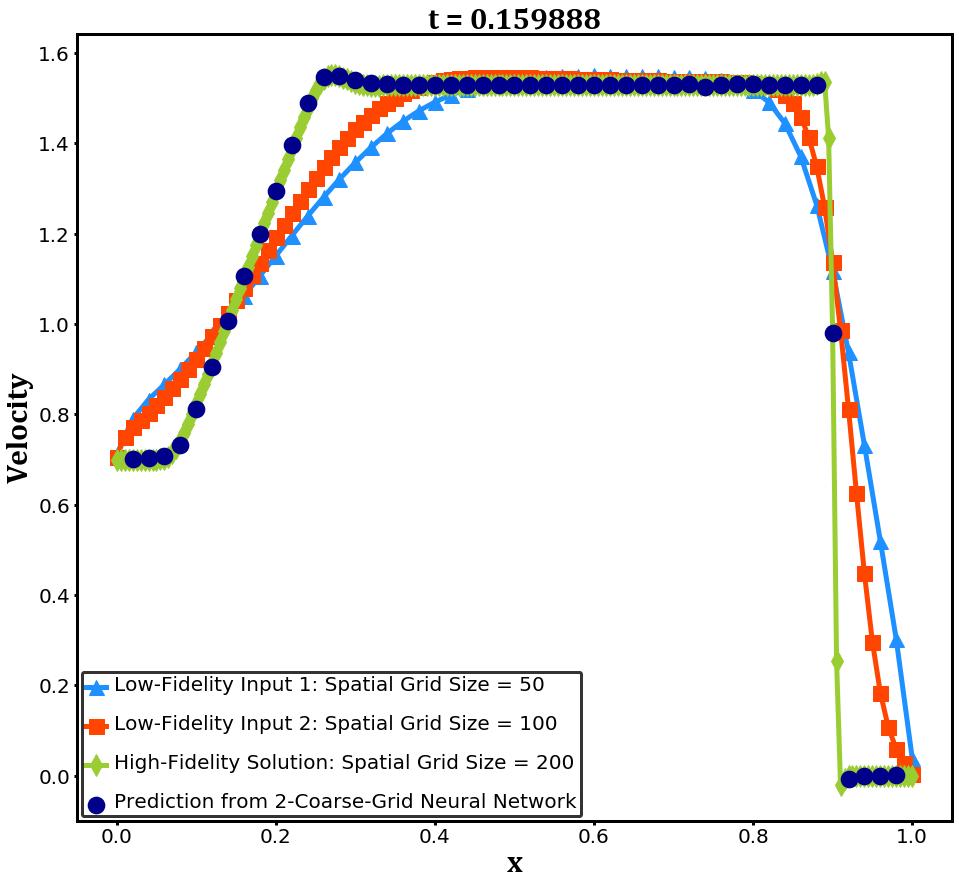}
\end{subfigure}
\caption{2CGNN prediction of final-time $(t=0.159888)$ solution (dark blue) of the Euler system with its {\bf initial value being $+3\%$ perturbation of that of the Lax problem}, low-fidelity input solutions (blue and red) by Rusanov scheme on $2$ different grids (with $50$ and $100$ cells resp.), and ``exact'' (reference) solution (green).}
\label{2CGNN Rusanov input: Final time of lax problem, +3}
\end{figure}

\begin{figure}[H]\centering
\begin{subfigure}[b]{.48\textwidth}
  \centering
  \includegraphics[width=1.0\linewidth]{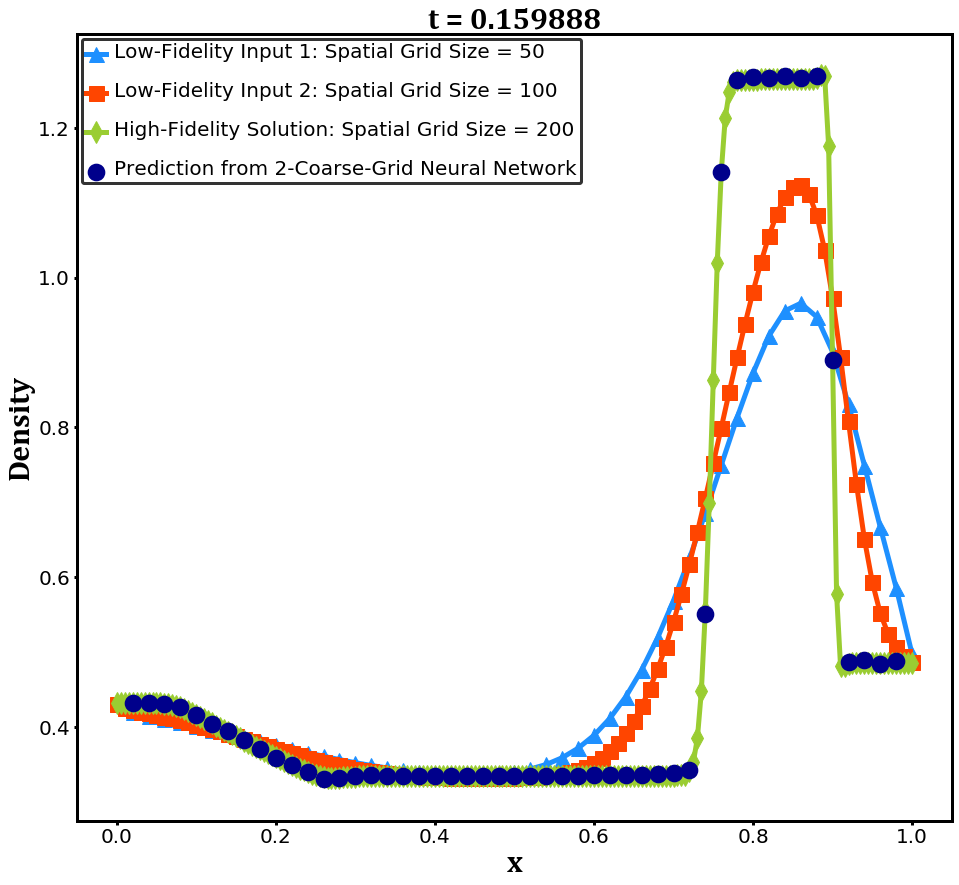}
\end{subfigure}
\begin{subfigure}[b]{.48\textwidth}
  \centering
  \includegraphics[width=1.0\linewidth]{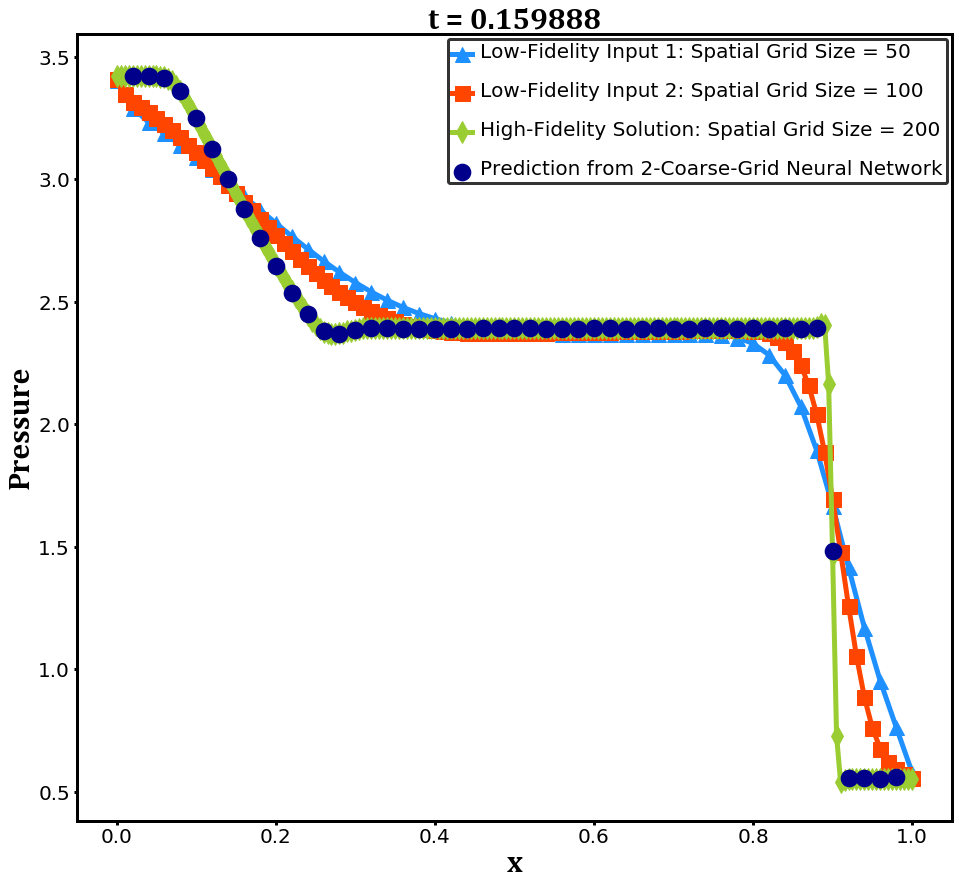}
\end{subfigure}
\begin{subfigure}[b]{.48\textwidth}
  \centering
  \includegraphics[width=1.0\linewidth]{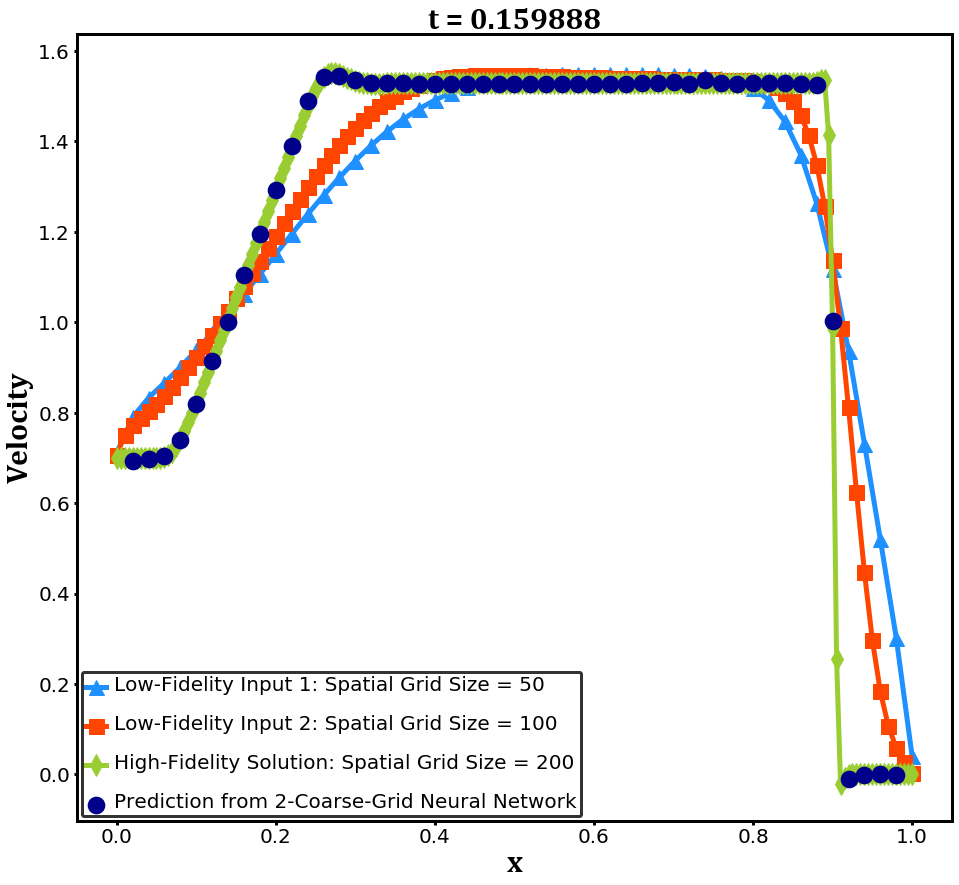}
\end{subfigure}
\caption{2CGNN prediction of final-time $(t=0.159888)$ solution (dark blue) of the Euler system with its {\bf initial value being $-3\%$ perturbation of that of the Lax problem}, low-fidelity input solutions (blue and red) by Rusanov scheme on $2$ different grids (with $50$ and $100$ cells resp.), and ``exact'' (reference) solution (green).}
\label{2CGNN Rusanov input: Final time of lax problem, -3}
\end{figure}

\begin{figure}[H]\centering
\begin{subfigure}[b]{.48\textwidth}
  \centering
  \includegraphics[width=1.0\linewidth]{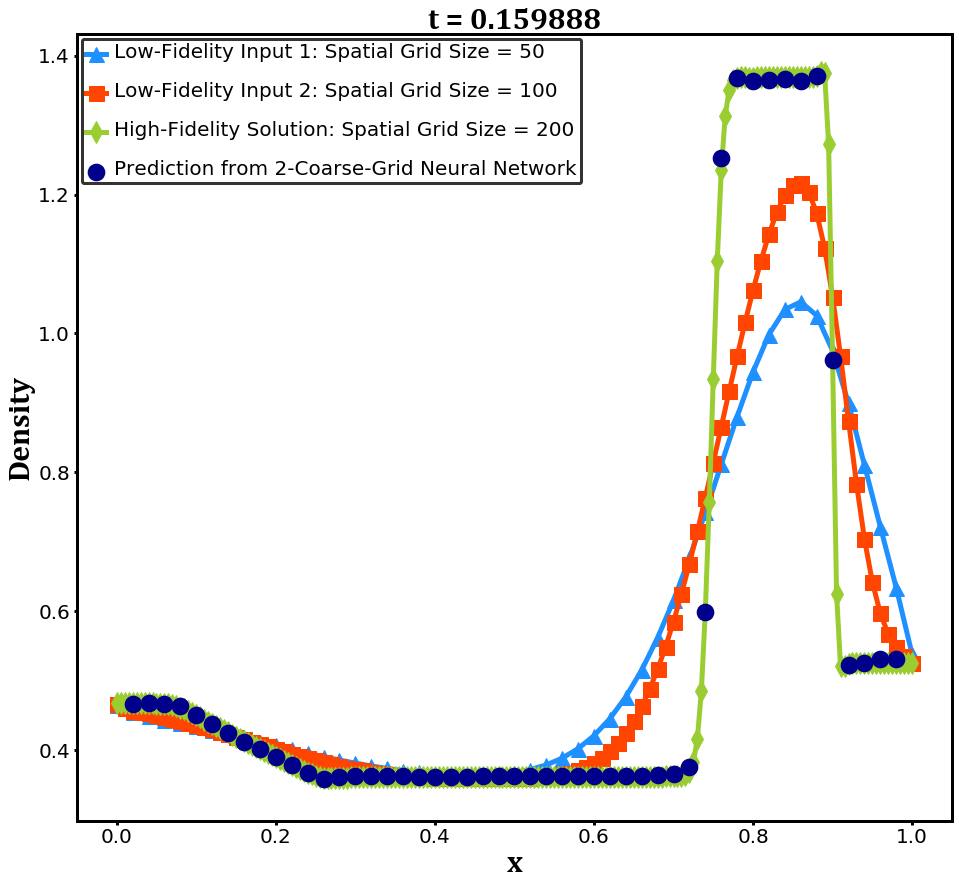}
\end{subfigure}
\begin{subfigure}[b]{.48\textwidth}
  \centering
  \includegraphics[width=1.0\linewidth]{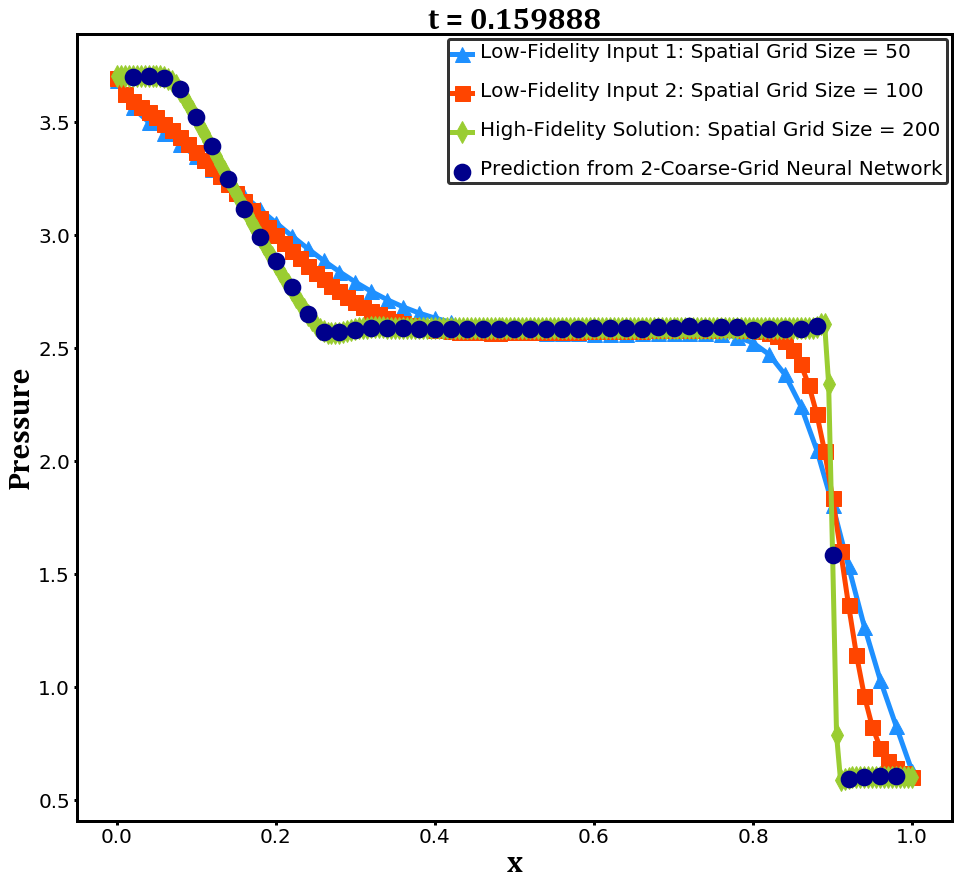}
\end{subfigure}
\begin{subfigure}[b]{.48\textwidth}
  \centering
  \includegraphics[width=1.0\linewidth]{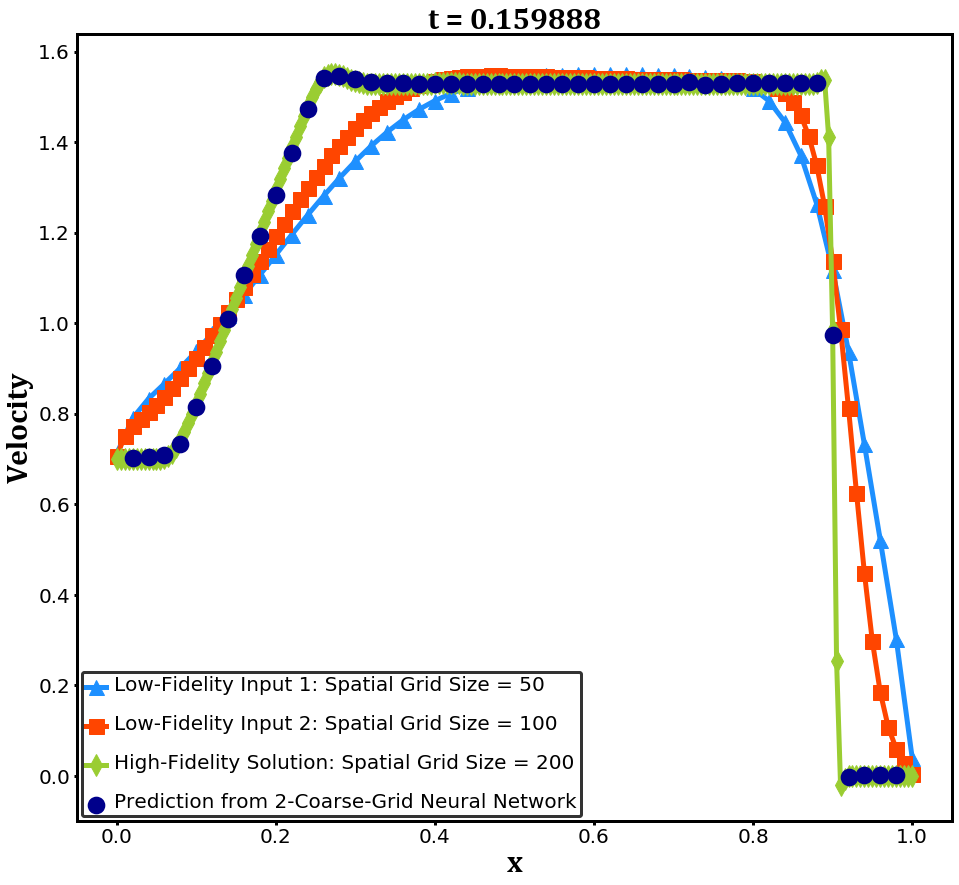}
\end{subfigure}
\caption{2CGNN prediction of final-time $(t=0.159888)$ solution (dark blue) of the Euler system with its {\bf initial value being $+5\%$ perturbation of that of the Lax problem}, low-fidelity input solutions (blue and red) by Rusanov scheme on $2$ different grids (with $50$ and $100$ cells resp.), and ``exact'' (reference) solution (green).}
\label{2CGNN Rusanov input: Final time of lax problem, +5}
\end{figure}

\begin{figure}[H]\centering
\begin{subfigure}[b]{.48\textwidth}
  \centering
  \includegraphics[width=1.0\linewidth]{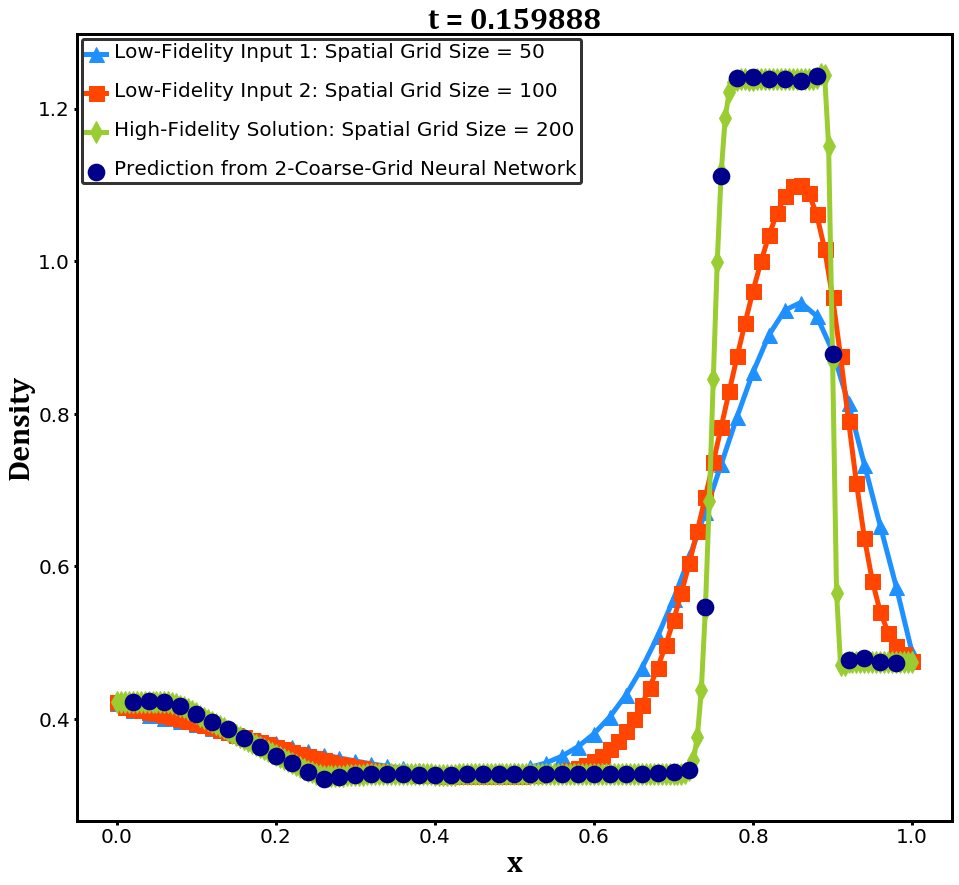}
\end{subfigure}
\begin{subfigure}[b]{.48\textwidth}
  \centering
  \includegraphics[width=1.0\linewidth]{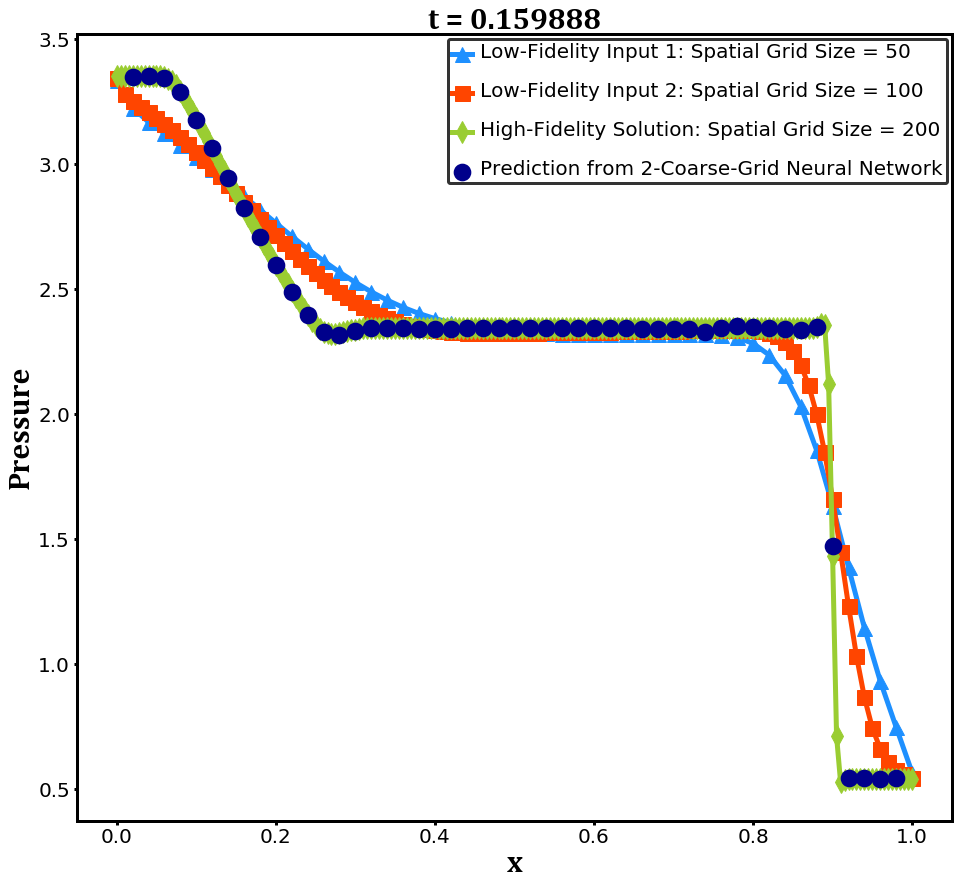}
\end{subfigure}
\begin{subfigure}[b]{.48\textwidth}
  \centering
  \includegraphics[width=1.0\linewidth]{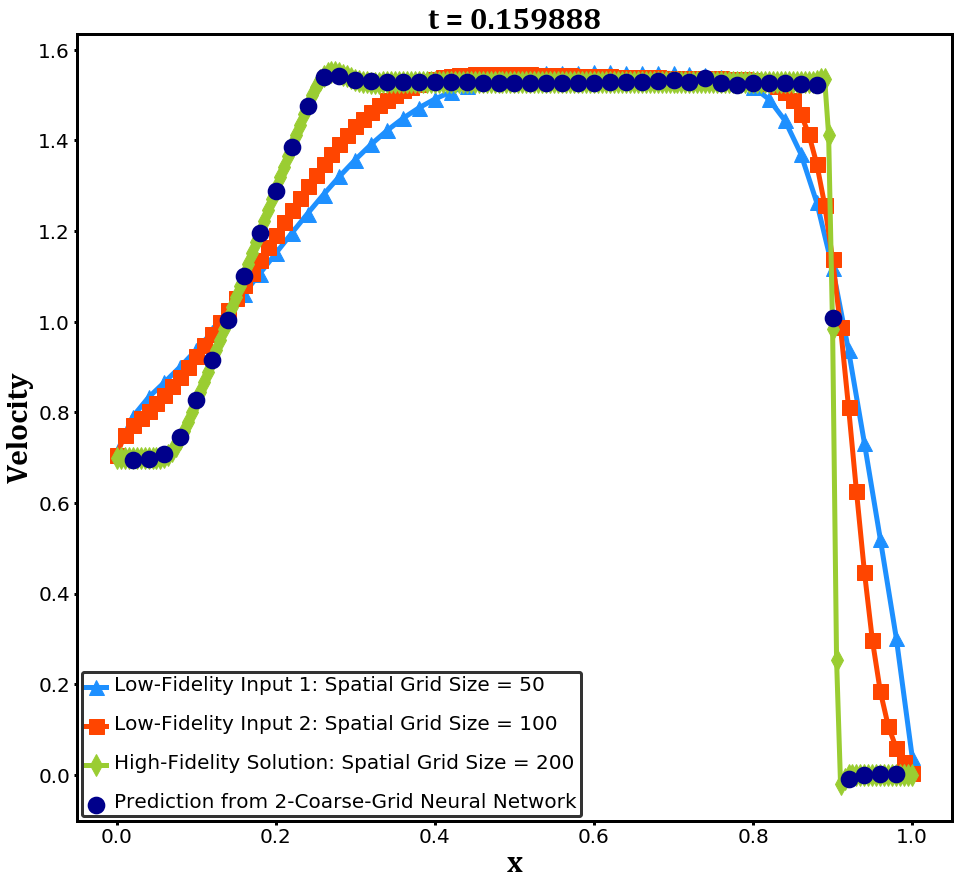}
\end{subfigure}
\caption{2CGNN prediction of final-time $(t=0.159888)$ solution (dark blue) of the Euler system with its {\bf initial value being $-5\%$ perturbation of that of the Lax problem}, low-fidelity input solutions (blue and red) by Rusanov scheme on $2$ different grids (with $50$ and $100$ cells resp.), and ``exact'' (reference) solution (green).}
\label{2CGNN Rusanov input: Final time of lax problem, -5}
\end{figure}

\subsection{Results for The Sod Problem}
The 2CGNN prediction of the final-time solution of the Sod problem is shown in Fig.~\ref{2CGNN Rusanov input: Final time of sod problem, original}. 
It's clear that 2CGNN captures the shock and contact really sharp, even without the small overshoot of the reference solution (note that similar reference solutions are used in the training of 2CGNN.)
Fig.~\ref{2CGNN Rusanov input: Final time of sod problem, +5} shows the prediction of the final-time solution of the Euler system with its initial value being $+5\%$ perturbation of that of the Sod problem. The 2CGNN predictions for
the final-time solutions of the Euler system with its initial value being $\pm  3\%$ and $-5\%$ perturbations of that of the Sod problem perform similarly well and are not shown here. 
Compared to N-PINN and PINN with additional data, 2CGNN once trained is much more efficient and accurate to predict the solution of the Euler system as long as it doesn't vary too much from the training data. 

\begin{figure}[H]\centering
\begin{subfigure}[b]{.48\textwidth}
  \centering
  \includegraphics[width=1.0\linewidth]{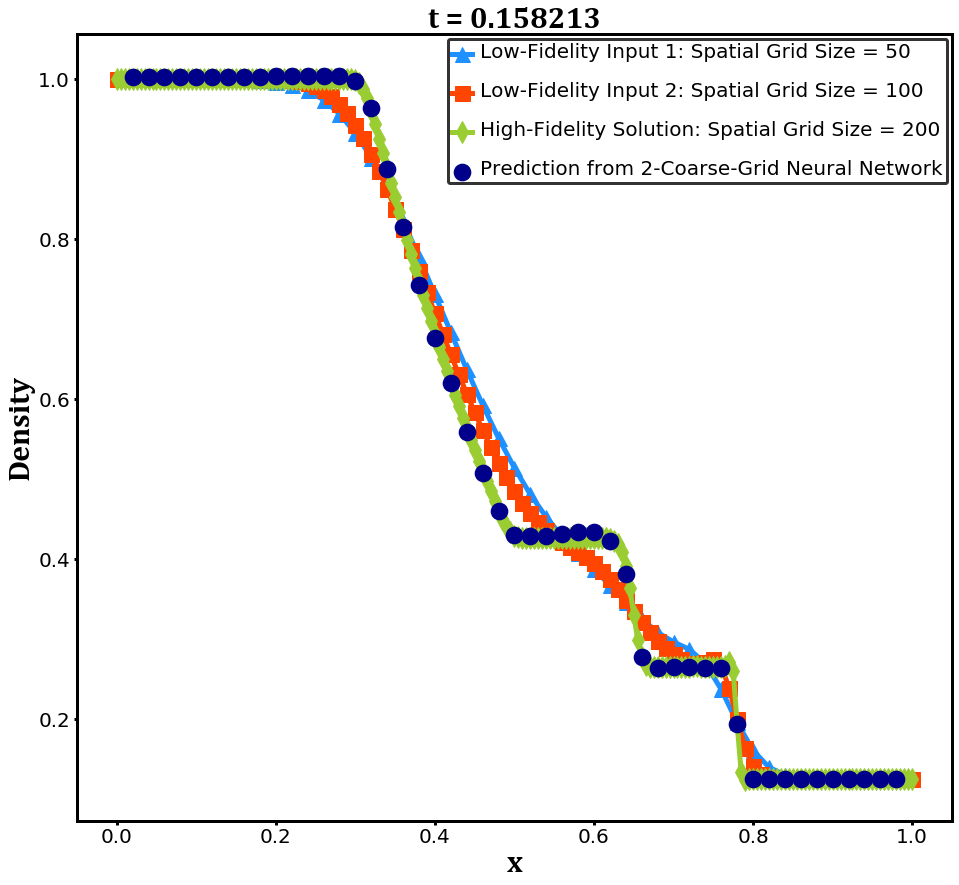}
\end{subfigure}
\begin{subfigure}[b]{.48\textwidth}
  \centering
  \includegraphics[width=1.0\linewidth]{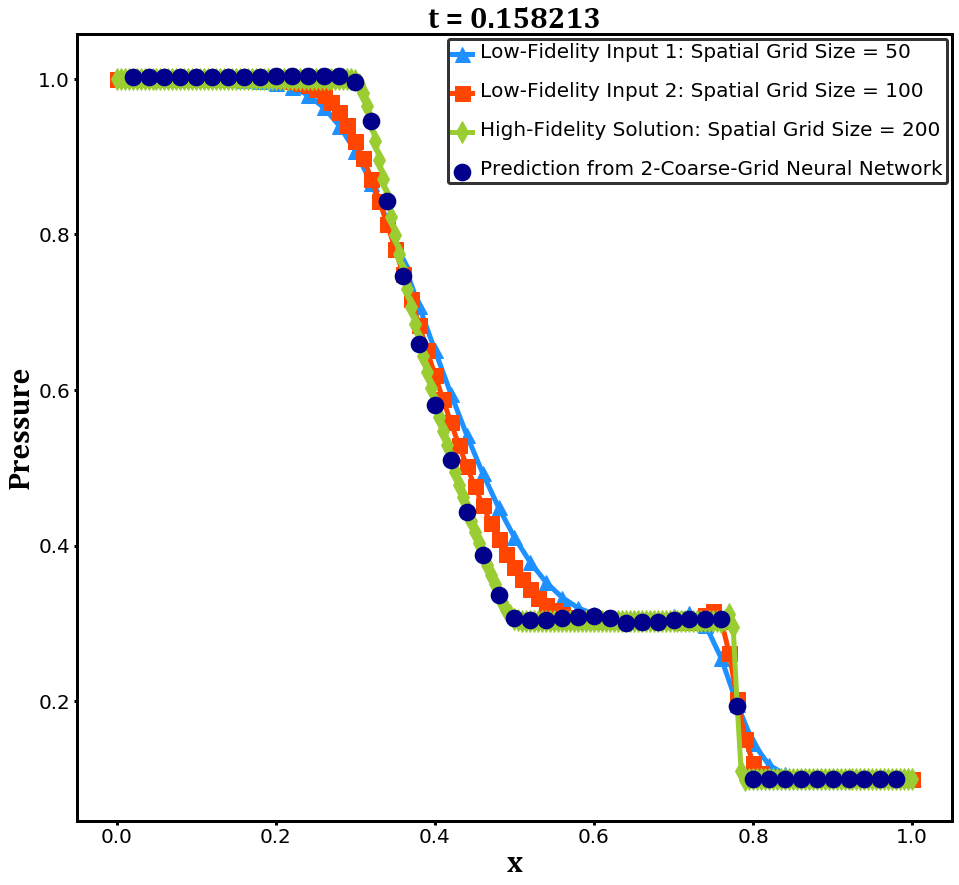}
\end{subfigure}
\begin{subfigure}[b]{.48\textwidth}
  \centering
  \includegraphics[width=1.0\linewidth]{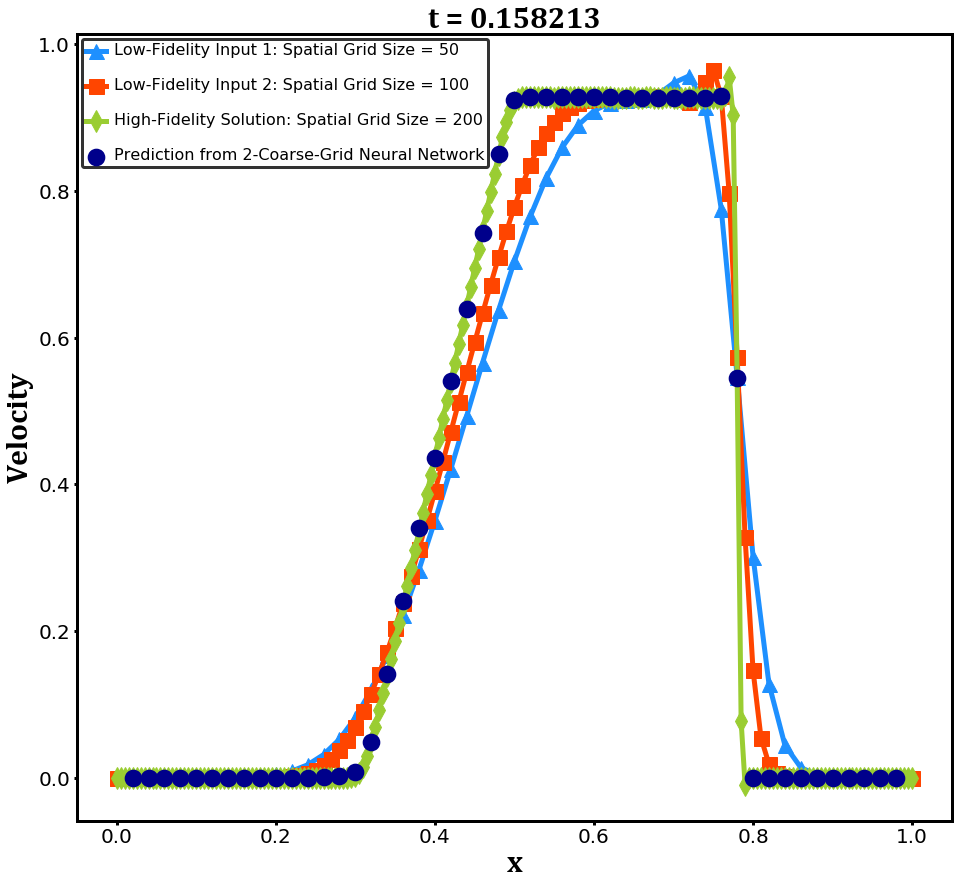}
\end{subfigure}
\caption{2CGNN prediction of final-time $(t=0.158213)$ solution (dark blue) of the {\bf Sod problem}, low-fidelity input solutions (blue and red) by Rusanov scheme on $2$ different grids (with $50$ and $100$ cells resp.), and ``exact'' (reference) solution (green).}
\label{2CGNN Rusanov input: Final time of sod problem, original}
\end{figure}

\begin{figure}[H]\centering
\begin{subfigure}[b]{.48\textwidth}
  \centering
  \includegraphics[width=1.0\linewidth]{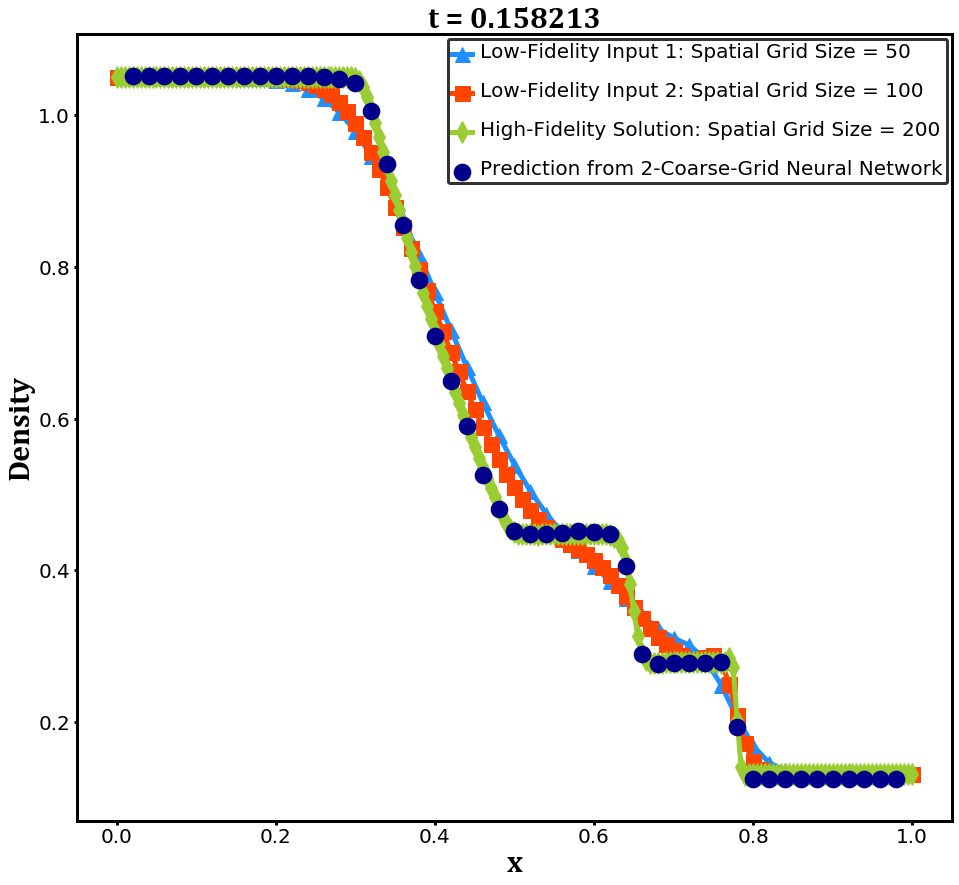}
\end{subfigure}
\begin{subfigure}[b]{.48\textwidth}
  \centering
  \includegraphics[width=1.0\linewidth]{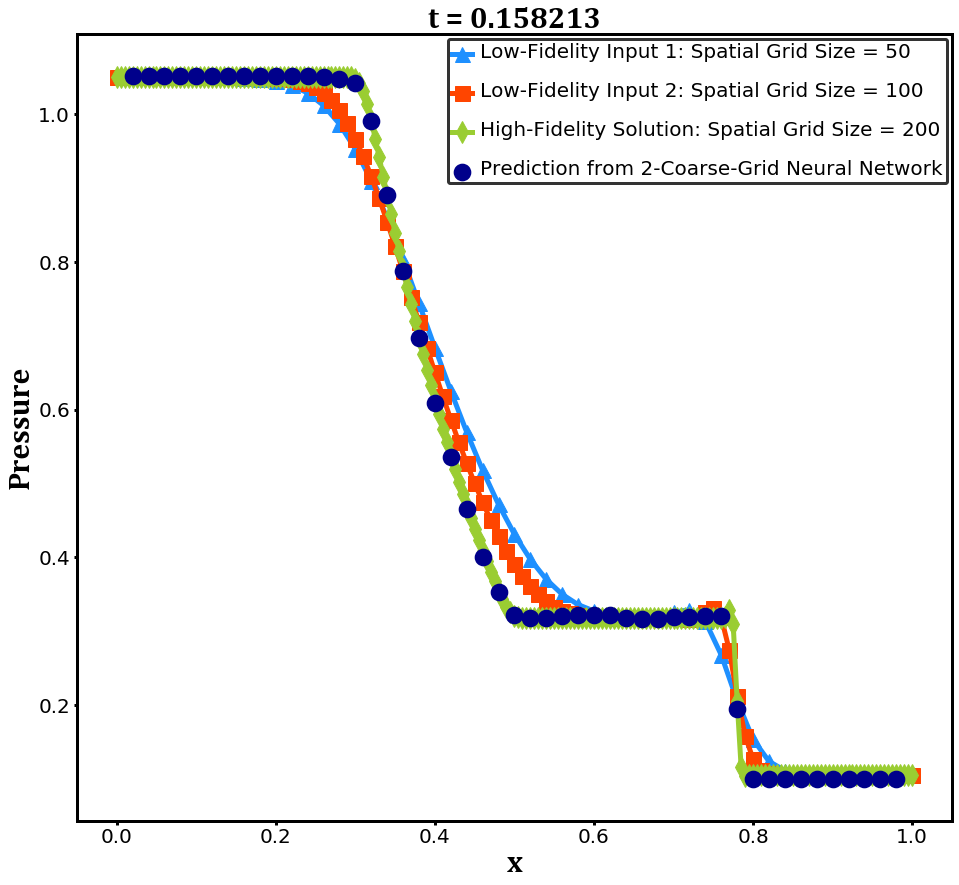}
\end{subfigure}
\begin{subfigure}[b]{.48\textwidth}
  \centering
  \includegraphics[width=1.0\linewidth]{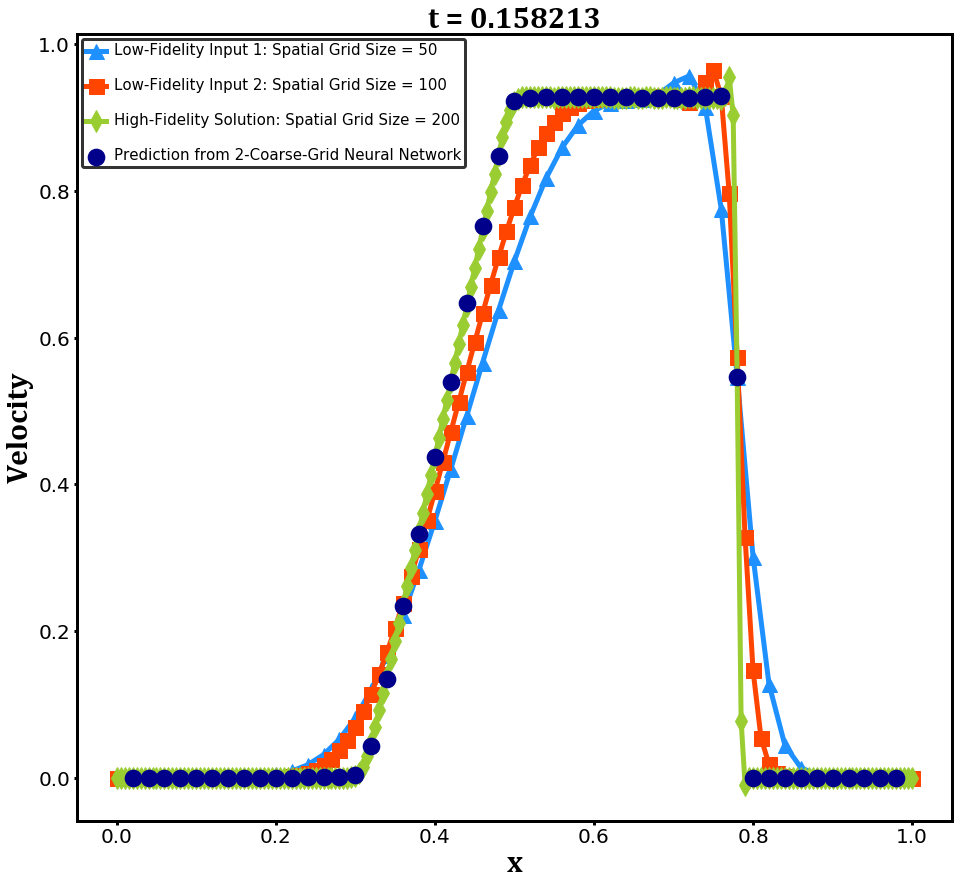}
\end{subfigure}
\caption{2CGNN prediction of final-time $(t=0.158213)$ solution (dark blue) of the Euler system with its {\bf initial value being $+5\%$ perturbation of that of the Sod problem}, low-fidelity input solutions (blue and red) by Rusanov scheme on $2$ different grids (with $50$ and $100$ cells resp.), and ``exact'' (reference) solution (green).}
\label{2CGNN Rusanov input: Final time of sod problem, +5}
\end{figure}

\subsection{Variations of the Input}
There could be alternative input formats. For example, one could use
as many data points as possible from the finer grid, and the input 
for a scalar conservation law can be changed from (\ref{standard-input}) to
\begin{equation}
\label{non-standard-input-1}
    \{u^{n'-1}_{i'-1}, u^{n'-1}_{i'}, u^{n'-1}_{i'+1}, u^{n'}_{i'}, u^{n''-2}_{i''-2}, u^{n''-2}_{i''-1}, u^{n''-2}_{i''}, u^{n''-2}_{i''+1}, u^{n''-2}_{i''+2}, u^{n''}_{i''}\}~,
\end{equation}
similarly for the Euler system. Prediction results from this type of input is similar to (slightly better than) those from input (\ref{standard-input}). 
Another variation is to change input (\ref{non-standard-input-1}) to
\begin{equation}
\label{non-standard-input-2}
    \{u^{n'-1}_{i'-1}, u^{n'-1}_{i'}, u^{n'-1}_{i'+1}, u^{n'}_{i'}, u^{n''-2}_{i''-1}, u^{n''-2}_{i''}, u^{n''-2}_{i''+1}, u^{n''}_{i''}\}~,
\end{equation}
similarly for the Euler system. Prediction results from this type of input is similar to (slightly worse than) those from input (\ref{standard-input}).

Input (\ref{non-standard-input-2}) computed by the Leapfrog diffusion and splitting scheme (\ref{leapfrog-diffusion-splitting}) generates some artifacts for predictions of the Sod problem, see Fig.~\ref{2CGNN non-standard2 Leapfrog and diffusion splitting input: Final time of sod problem, original}.

\begin{figure}[H]\centering
\begin{subfigure}[b]{.48\textwidth}
    \centering
    \includegraphics[width=1.0\linewidth]{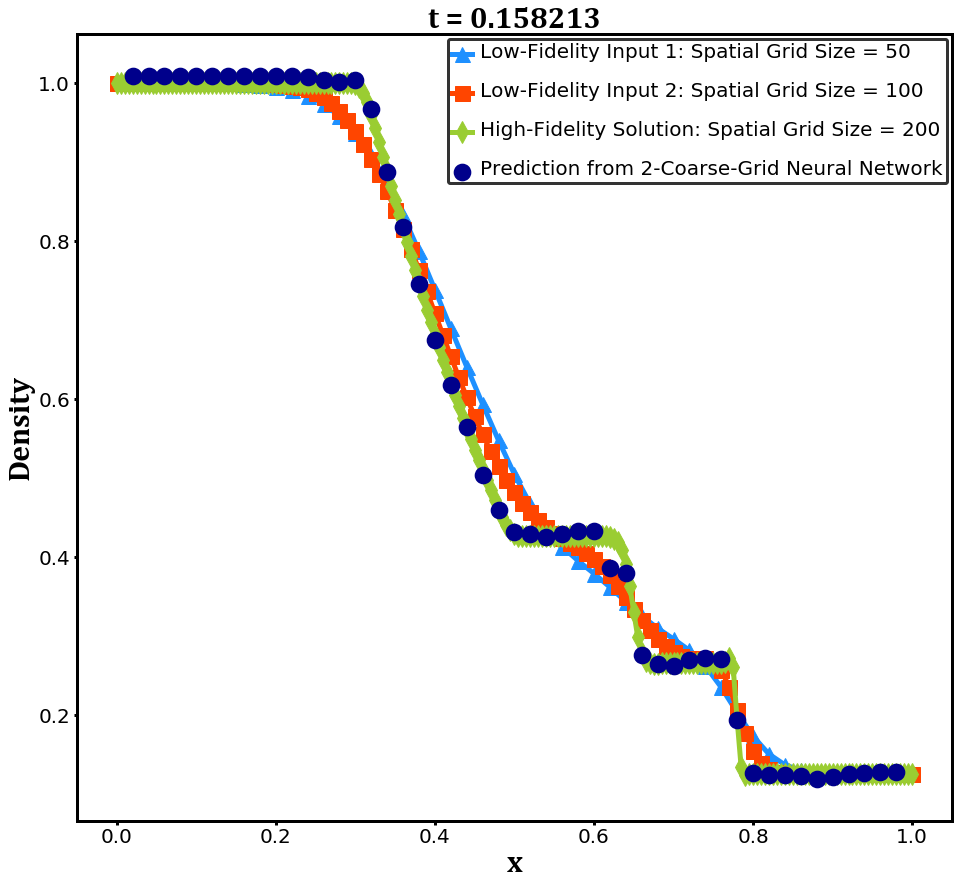}
\end{subfigure}
\begin{subfigure}[b]{.48\textwidth}
  \centering
  \includegraphics[width=1.0\linewidth]{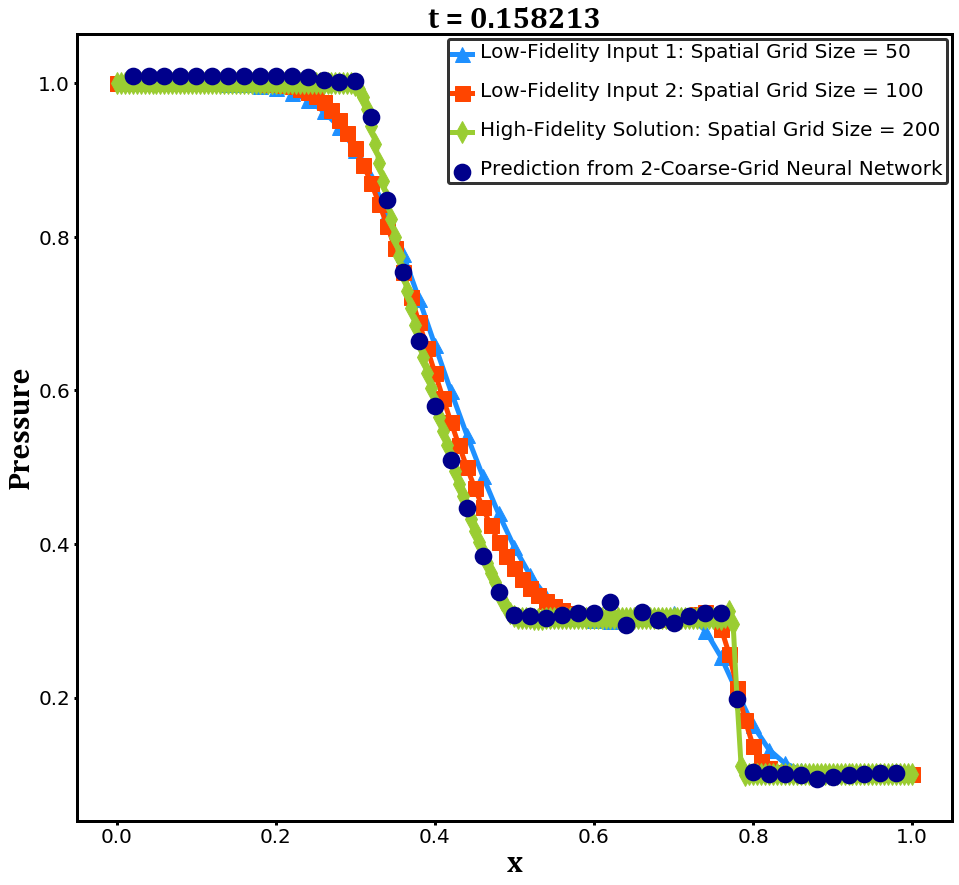}
\end{subfigure}
\caption{2CGNN (with non-standard input (\ref{non-standard-input-2})) prediction of final-time $(t=0.158213)$ solution (dark blue) of the {\bf Sod problem}, low-fidelity input solutions (blue and red) by Leapfrog and diffusion splitting scheme (\ref{leapfrog-diffusion-splitting}) on $2$ different grids (with $50$ and $100$ cells resp.), and ``exact'' (reference) solution (green).}
\label{2CGNN non-standard2 Leapfrog and diffusion splitting input: Final time of sod problem, original}
\end{figure}

We use $9$ layers and $100$ neurons on each layer for 2CGNN (with non-standard input (\ref{non-standard-input-1})) to predict the Sod problem while there is no change on neurons and layers for other cases. Table~\ref{2CGNN-1 Leapfrog and diffusion splitting input results} summarizes the prediction results of input (\ref{non-standard-input-1}) computed by the leapfrog diffusion splitting scheme (\ref{leapfrog-diffusion-splitting}).  

\begin{table}[H]
\centering
\caption{Relative $l_2$ errors of 2CGNN (with non-standard input (\ref{non-standard-input-1})) predictions of solutions on coarsest uniform spatial grid (50 cells) of the Euler system (with initial values perturbed from those of the Lax or Sod problem) with low-fidelity input solutions computed by the Leapfrog and diffusion splitting scheme (\ref{leapfrog-diffusion-splitting}) on $2$ different grids (with $50$ and $100$ cells.)}
\begin{tabular}{|c|c|c|} \hline
&\textbf{Lax Problem} &\textbf{Sod Problem} \\ \hline
Initial Value&Relative $l_2$ Error &Relative $l_2$ Error \\ \hline
Original &$8.17e-03$ &$8.35e-03$ \\ \hline
$+3\%$ &$4.06e-03$ &$2.81e-02$ \\ \hline
$-3\%$ &$5.96e-03$ &$6.63e-03$ \\ \hline
$+5\%$ &$5.27e-03$ &$9.37e-03$ \\ \hline
$-5\%$ &$1.06e-02$ &$9.47e-03$ \\ \hline
$+7\%$ &$8.51e-03$ &$8.27e-03$ \\ \hline
$-7\%$ &$1.49e-02$ &$2.46e-02$ \\ \hline
\end{tabular}
\label{2CGNN-1 Leapfrog and diffusion splitting input results}
\end{table}

\section{2-Diffusion-Coefficient Neural Network }
Instead of computing for the input on two grids, one could generate the input on
a single grid but for two perturbed equations with two different diffusion coefficients. This is referred to as the 2-Diffusion-Coefficient Neural Network (2DCNN.) The perturbed scalar conservation law can be written as
\begin{equation}
\label{cons-law+diffusion}
   \frac{\partial U}{\partial t} + \frac{\partial f(U)}{\partial x} = \alpha \frac{\partial^2 U}{\partial x^2}, \; x\in\Omega\subset \mathcal{R},\; t\in[0, T]~.
\end{equation}
To fill the two parts of an input, one could compute equation (\ref{cons-law+diffusion}) with two different values of $\alpha$ using a low-cost scheme on a single grid. 
For example, we can approximate (\ref{cons-law+diffusion}) using the Leapfrog and diffusion splitting scheme (\ref{leapfrog-diffusion-splitting}) with $\alpha=\Delta x$ and $c\Delta x$ as follows
\begin{equation}
\label{leapfrog-diffusion-splitting-diff-coef-1}
\left \{
\begin{array}{l}
    \frac{\tilde{U}_i - U^{n-1}_i}{2\Delta t} +\frac{{f(U)}|^{n}_{i+1}-{f(U)}|^{n}_{i-1}}{2\Delta x} =0~, \\
    \frac{U^{n+1}_i - \tilde{U}_i}{\Delta t}- \Delta x \frac{\tilde{U}_{i+1} - 2\cdot \tilde{U}_i + \tilde{U}_{i-1}}{\Delta x^2}=0~,
\end{array}
\right .
\end{equation}
and 
\begin{equation}
\label{leapfrog-diffusion-splitting-diff-coef-2}
\left \{
\begin{array}{l}
    \frac{\tilde{V}_i - V^{n-1}_i}{2\Delta t} +\frac{{f(V)}|^{n}_{i+1}-{f(V)}|^{n}_{i-1}}{2\Delta x} =0~, \\
    \frac{V^{n+1}_i - \tilde{V}_i}{\Delta t}- c\Delta x \frac{\tilde{V}_{i+1} - 2\cdot \tilde{V}_i + \tilde{V}_{i-1}}{\Delta x^2}=0~.
\end{array}
\right .
\end{equation}
Then the input computed on a single uniform grid can be written as
\begin{equation}
\label{2DCNN-input}
    \{U^{n-1}_{i-1}, U^{n-1}_{i}, U^{n-1}_{i+1}, U^{n}_{i}, V^{n-1}_{i-1}, V^{n-1}_{i}, V^{n-1}_{i+1}, V^{n}_{i}\}~.
\end{equation}
similarly for the Euler system.
See Fig.~\ref{2DCNN Leapfrog and diffusion splitting input: Final time of Lax problem, dx and 4dx} and \ref{2DCNN Leapfrog and diffusion splitting input: Final time of Lax problem +5, dx and 4dx} for the predictions by 2DCNN, with the input computed on a uniform grid with $100$ cells and $c=4$. Note that for the second equation of
(\ref{leapfrog-diffusion-splitting-diff-coef-2}), the larger diffusion coefficient may require a time step size smaller than that of the first equation of (\ref{leapfrog-diffusion-splitting-diff-coef-2}). Therefore one may need to compute the second equation of
(\ref{leapfrog-diffusion-splitting-diff-coef-2}) in two
time steps, with the time step size $\frac12\Delta t$ for each time step.

\begin{figure}[H]\centering
\begin{subfigure}[b]{.48\textwidth}
  \centering
  \includegraphics[width=1.0\linewidth]{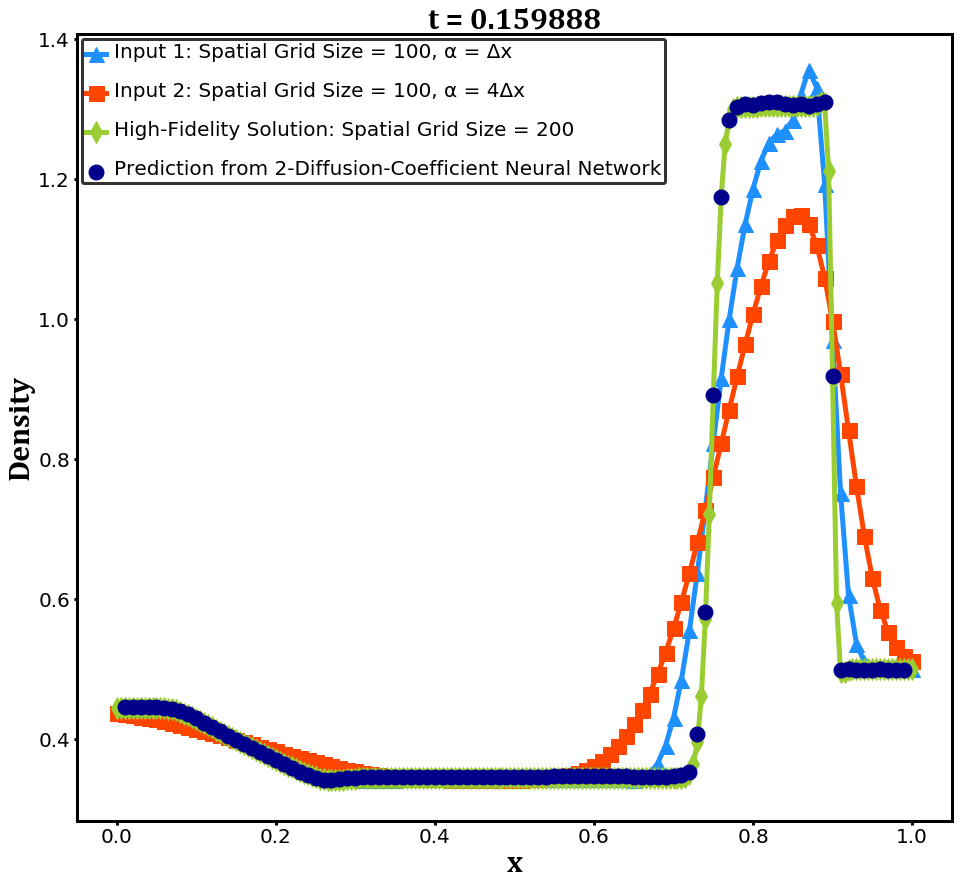}
\end{subfigure}
\begin{subfigure}[b]{.48\textwidth}
  \centering
  \includegraphics[width=1.0\linewidth]{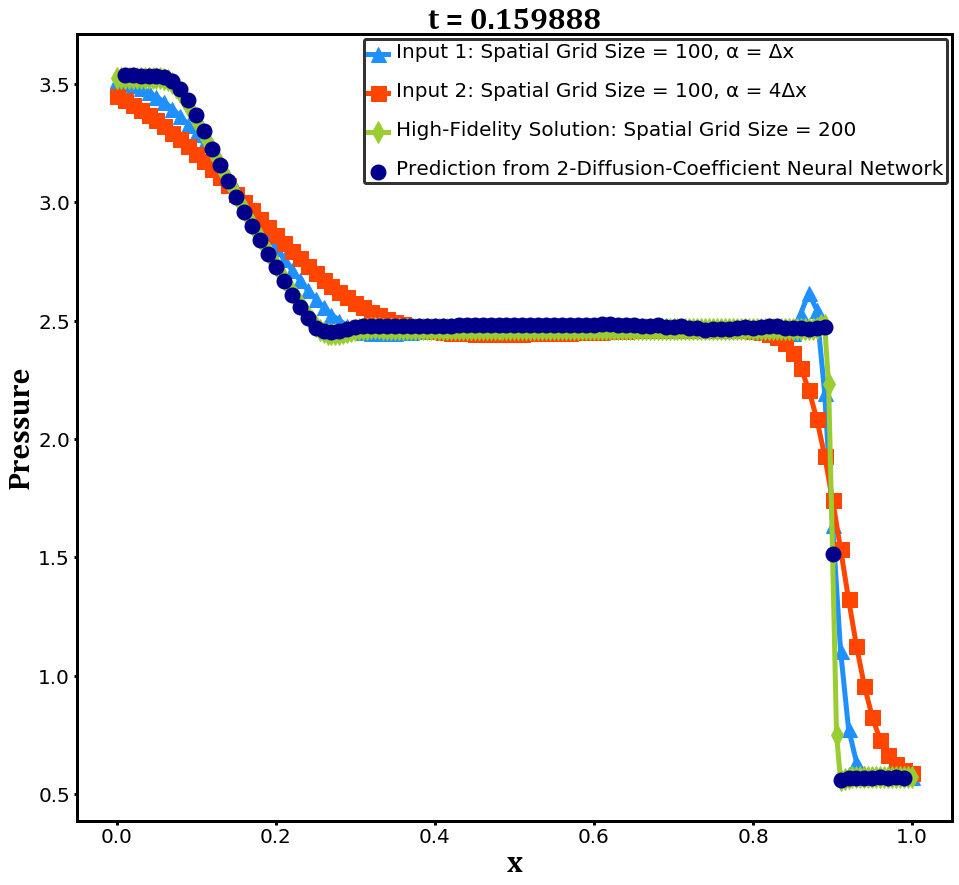}
\end{subfigure}
\begin{subfigure}[b]{.48\textwidth}
  \centering
  \includegraphics[width=1.0\linewidth]{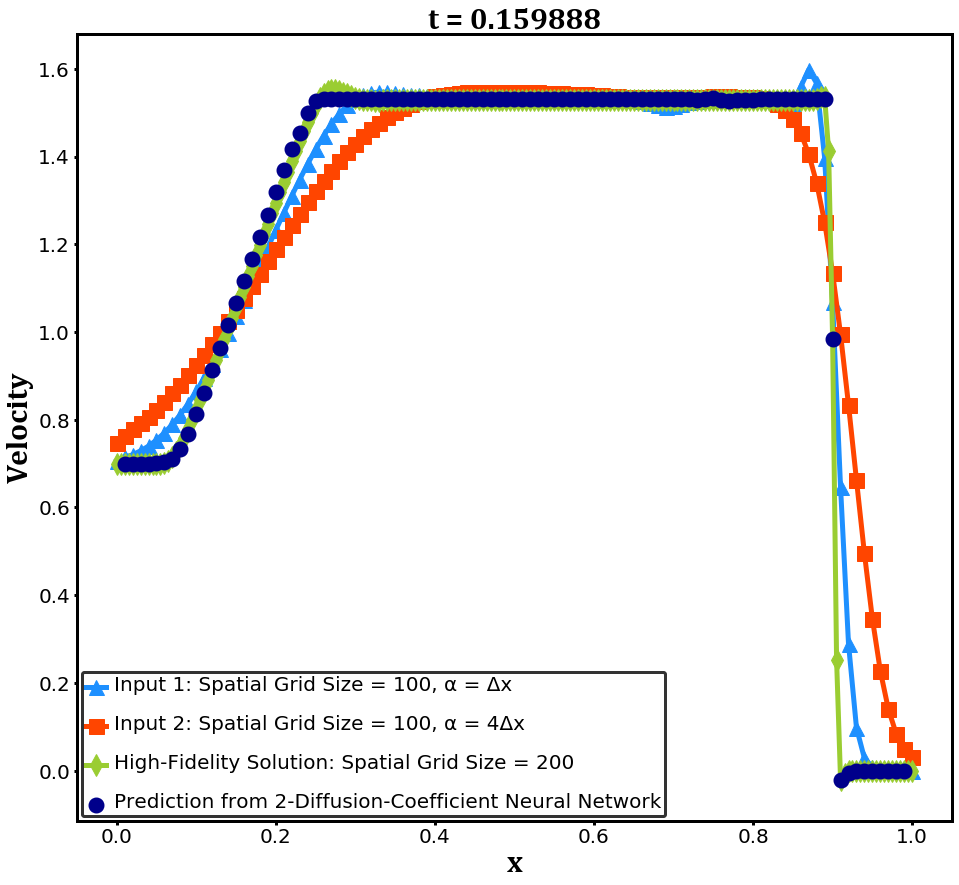}
\end{subfigure}
\caption{2DCNN prediction of final-time $(t=0.159888)$ solution of {\bf Lax problem}  (dark blue), low-fidelity input solutions (blue and red) by the Leapfrog and diffusion splitting scheme (\ref{leapfrog-diffusion-splitting}) on the uniform grid with $100$ cells, with diffusion coefficient $\alpha = \Delta x$ and $4\Delta x$ resp., and ``exact'' (reference) solution (green).
}
\label{2DCNN Leapfrog and diffusion splitting input: Final time of Lax problem, dx and 4dx}
\end{figure}

\begin{figure}[H]\centering
\begin{subfigure}[b]{.48\textwidth}
  \centering
  \includegraphics[width=1.0\linewidth]{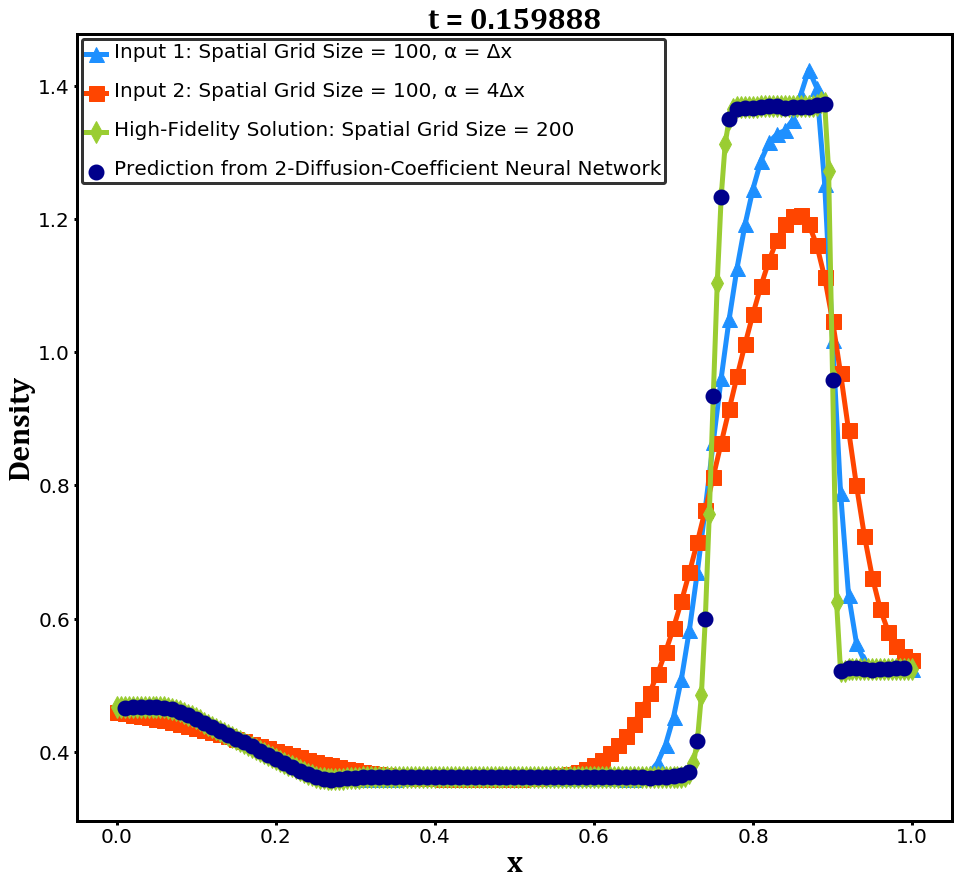}
\end{subfigure}
\begin{subfigure}[b]{.48\textwidth}
  \centering
  \includegraphics[width=1.0\linewidth]{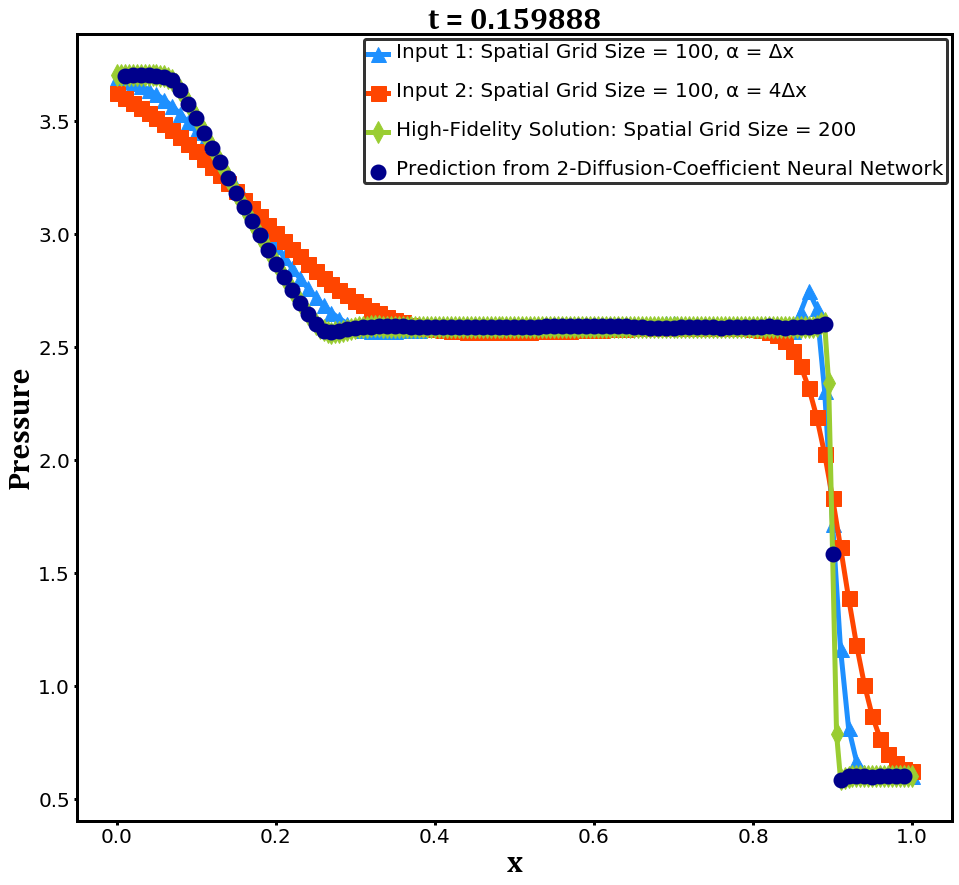}
\end{subfigure}
\begin{subfigure}[b]{.48\textwidth}
  \centering
  \includegraphics[width=1.0\linewidth]{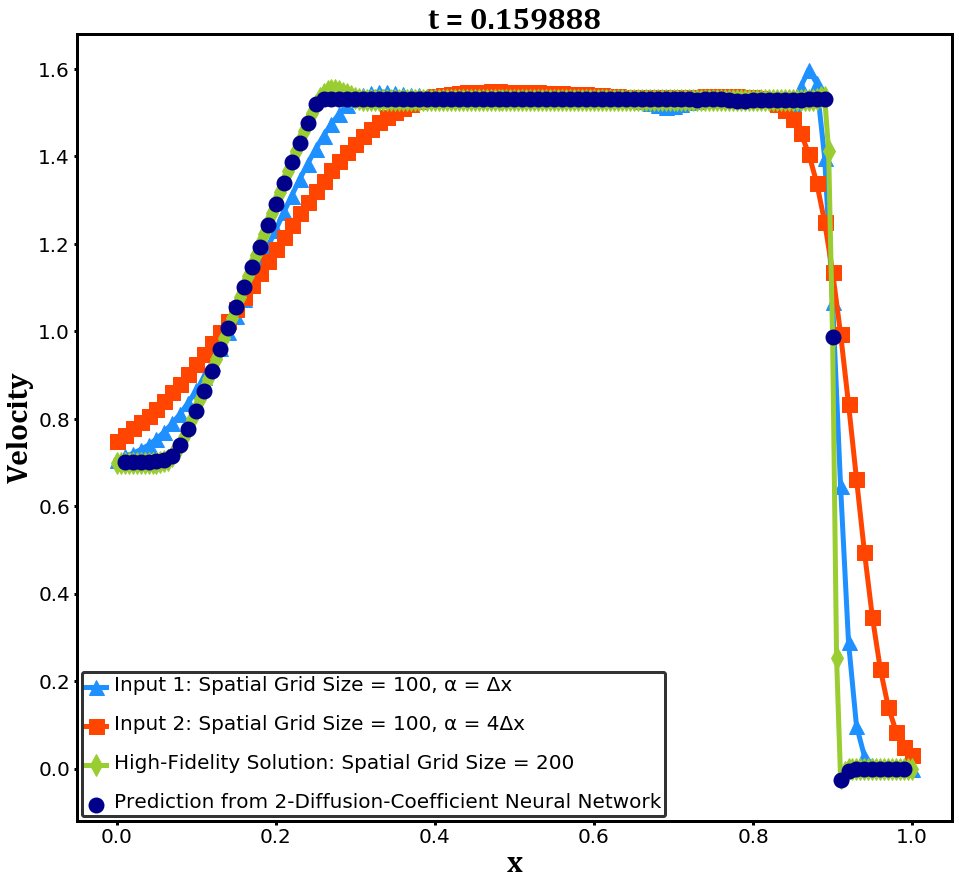}
\end{subfigure}
\caption{2DCNN prediction of final-time $(t=0.159888)$ solution (dark blue) of the Euler system with its {\bf initial value being $+5\%$ perturbation of that of the Lax problem}, low-fidelity input solutions (blue and red) by the Leapfrog and diffusion splitting scheme (\ref{leapfrog-diffusion-splitting}) on the uniform grid with $100$ cells, with diffusion coefficient $\alpha = \Delta x$ and $4\Delta x$ resp., and ``exact'' (reference) solution (green).}
\label{2DCNN Leapfrog and diffusion splitting input: Final time of Lax problem +5, dx and 4dx}
\end{figure}


\section{Summary and Comparison of All Methods}
The tables below summarize the results of all methods we have tested in aforementioned sections. 
\begin{table}[H]
\centering
\caption{Comparison of N-PINN: Relative $l_2$ norm errors between the predicted and ``exact'' (reference) solutions for the Lax problem, Sod problem and Burgers’ equation.}\label{tb: Summary of N-PINN errors }
\begin{tabular}{|c|c|c|c|c|c|c|c|} \hline
\textbf{N-PINN with } & \multicolumn{3}{|c|}{\textbf{Lax Problem}} & \multicolumn{3}{|c|}{\textbf{Sod Problem}} & \textbf{Burgers'}\\ & \multicolumn{3}{|c|}{ } & \multicolumn{3}{|c|}{ } & \textbf{Equation}\\ \hline
& Density & Pressure & Velocity & Density & Pressure & Velocity &  \\ \hline
Rusanov Scheme & $1.43e-01$ & $3.85e-02$ & $8.63e-02$ & $6.06e-02$ & $5.63e-02$ & $2.76e-01$ & $1.33e-02$\\ \hline
Leapfrog and Diffusion & $1.67e-01$ & $4.63e-02$ & $9.07e-02$ & $4.75e-02$ & $4.94e-02$ & $2.25e-01$ & $1.11e-02$\\ \hline
Leapfrog and Diffusion  & & & & & & &\\ w/ Additional Data  & $4.33e-02$ & $7.94e-03$ & $2.79e-02$ & $1.26e-02$ & $1.56e-02$ & $1.81e-02$ & N/A\\ 
\hline\end{tabular}
\label{N-PINN Results}
\end{table}

\begin{table}[H]
\centering
\caption{Relative $l_2$ errors of 2CGNN predictions of solutions on coarsest uniform spatial grid (50 cells) of the Euler system (with initial values perturbed from those of the Lax or Sod problem) with low-fidelity input solutions computed by the Rusanov scheme on $2$ different grids (with $50$ and $100$ cells.)}\label{tb: Summary of total errors from predictions, Rusanov scheme}
\begin{tabular}{|c|c|c|} \hline
&\textbf{Lax Problem} &\textbf{Sod Problem} \\ \hline
Initial Value&Relative $l_2$ Error &Relative $l_2$ Error \\ \hline
Original &$5.97e-03$ &$6.20e-03$ \\ \hline
$+3\%$ &$7.29e-03$ &$5.66e-03$ \\ \hline
$-3\%$ &$9.10e-03$ &$9.05e-03$ \\ \hline
$+5\%$ &$8.61e-03$ &$7.84e-03$ \\ \hline
$-5\%$ &$1.32e-02$ &$1.29e-02$ \\ \hline
$+7\%$ &$9.96e-03$ &$9.72e-03$ \\ \hline
$-7\%$ &$1.75e-02$ &$1.45e-02$ \\ \hline
\end{tabular}
\label{2CGNN Rusanov input results}
\end{table}

\begin{table}[H]
\centering
\caption{Relative $l_2$ errors of 2CGNN predictions of solutions on coarsest uniform spatial grid (50 cells) of the Euler system (with initial values perturbed from those of the Lax or Sod problem) with low-fidelity input solutions computed by the stabilized Leapfrog and diffusion scheme (\ref{leapfrog-diffusion-stabilized}) on $2$ different grids (with $50$ and $100$ cells.)}\label{tb: Summary of total errors from predictions, Leapfrog Diffusion scheme}
\begin{tabular}{|c|c|c|} \hline
&\textbf{Lax Problem} &\textbf{Sod Problem} \\ \hline
Initial Value&Relative $l_2$ Error &Relative $l_2$ Error \\ \hline
Original &$8.73e-03$ &$6.22e-03$ \\ \hline
$+3\%$ &$5.47e-03$ &$7.68e-03$ \\ \hline
$-3\%$ &$1.14e-02$ &$1.13e-02$ \\ \hline
$+5\%$ &$1.01e-02$ &$6.47e-03$ \\ \hline
$-5\%$ &$1.29e-02$ &$1.11e-02$ \\ \hline
$+7\%$ &$1.42e-02$ &$8.37e-03$ \\ \hline
$-7\%$ &$1.53e-02$ &$1.77e-02$ \\ \hline
\end{tabular}
\label{2CGNN Leapfrog Diff. input results}
\end{table}

\begin{table}[H]
\centering
\caption{Relative $l_2$ errors of 2CGNN predictions of solutions on coarsest uniform spatial grid (50 cells) of the Euler system (with initial values perturbed from those of the Lax or Sod problem) with low-fidelity input solutions computed by the Leapfrog and diffusion splitting scheme (\ref{leapfrog-diffusion-splitting}) on $2$ different grids (with $50$ and $100$ cells.)
}\label{tb: Summary of total errors from predictions, Leapfrog Diffusion scheme with operator splitting}
\begin{tabular}{|c|c|c|} \hline
&\textbf{Lax Problem} &\textbf{Sod Problem} \\ \hline
Initial Value&Relative $l_2$ Error &Relative $l_2$ Error\\ \hline
Original &$8.56e-03$ &$5.72e-03$ \\ \hline
$+3\%$ &$4.80e-03$ &$8.03e-03$ \\ \hline
$-3\%$ &$6.86e-03$ &$8.26e-03$ \\ \hline
$+5\%$ &$4.96e-03$ &$6.60e-03$ \\ \hline
$-5\%$ &$7.05e-03$ &$1.84e-02$ \\ \hline
$+7\%$ &$6.22e-03$ &$8.55e-03$ \\ \hline
$-7\%$ &$8.28e-03$ &$1.05e-02$ \\ \hline
\end{tabular}
\label{2CGNN Leapfrog Diff. operator split input results}
\end{table}

\begin{table}[H]
\centering
\caption{Relative $l_2$ errors of 2DCNN predictions of solutions on finer uniform spatial grid (100 cells) of the Euler system (with initial values perturbed from those of the Lax or Sod problem) with low-fidelity input solutions computed by the Leapfrog and diffusion splitting scheme (\ref{leapfrog-diffusion-splitting}) with $2$ different diffusion coefficients ($\alpha=\Delta x$ and $4\Delta x$) on a uniform grid with $100$ cells.
}\label{tb: Summary of total errors from 2DCNN predictions, Leapfrog Diffusion scheme with operator splitting}
\begin{tabular}{|c|c|c|} \hline
&\textbf{Lax Problem} &\textbf{Sod Problem} \\ \hline
Initial Value&Relative $l_2$ Error &Relative $l_2$ Error\\ \hline
Original &$5.24e-03$ &$6.47e-03$ \\ \hline
$+3\%$ &$3.43e-03$ &$6.80e-03$ \\ \hline
$-3\%$ &$5.34e-03$ &$8.57e-03$ \\ \hline
$+5\%$ &$7.13e-03$ &$7.52e-03$ \\ \hline
$-5\%$ &$6.59e-03$ &$1.07e-02$ \\ \hline
$+7\%$ &$8.19e-03$ &$8.87e-03$ \\ \hline
$-7\%$ &$7.51e-03$ &$1.30e-02$ \\ \hline
\end{tabular}
\label{2DCNN Leapfrog Diff. operator split input results}
\end{table}

\section{Conclusion}
We have studied a few neural networks based on first order schemes to solve
the Burgers' equation and one dimensional Euler system. An N-PINN is able to
generate an approximate solution comparable to those obtained from first order schemes. Given additional data of the solution, N-PINN can improve its approximation conveniently. We also introduce 2CGNN and 2DCNN which significantly improve from N-PINN and are efficient to use after they have been trained. In fact, they are able to output very sharp shocks and contacts out of smeared solution profiles in the input. And they don't seem to be very sensitive to
the low-cost schemes used for computing inputs either, see Appendix A and B for additional prediction results with some alternative schemes used for computing inputs.
We will extend 2CGNN and 2DCNN to multidimensional nonlinear hyperbolic systems in the future.

\section*{Acknowledgements}

The authors thank Prof. Wen Shen, Penn. State U. for a helpful discussion.

\section*{Appendix A. Stabilized Leapfrog and Diffusion Scheme (\ref{leapfrog-diffusion-stabilized}) for Computing Inputs for 2CGNN}

\begin{figure}[H]\centering
\begin{subfigure}[b]{.48\textwidth}
  \centering
  \includegraphics[width=1.0\linewidth]{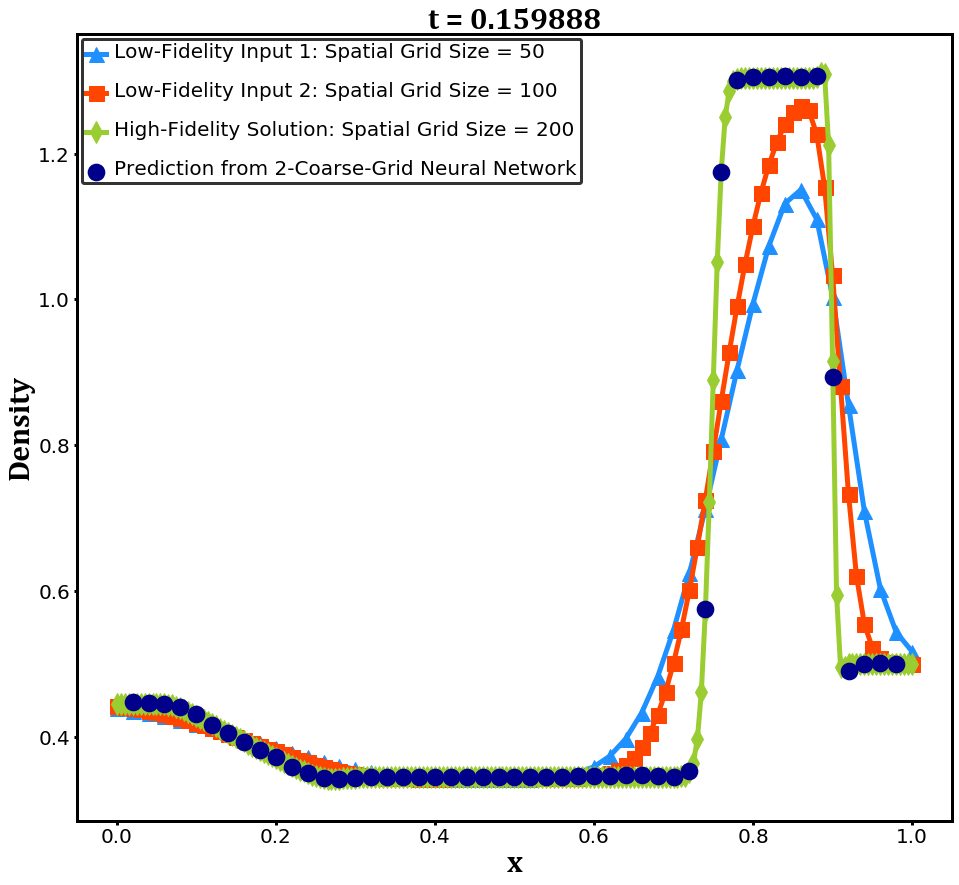}
\end{subfigure}
\begin{subfigure}[b]{.48\textwidth}
  \centering
  \includegraphics[width=1.0\linewidth]{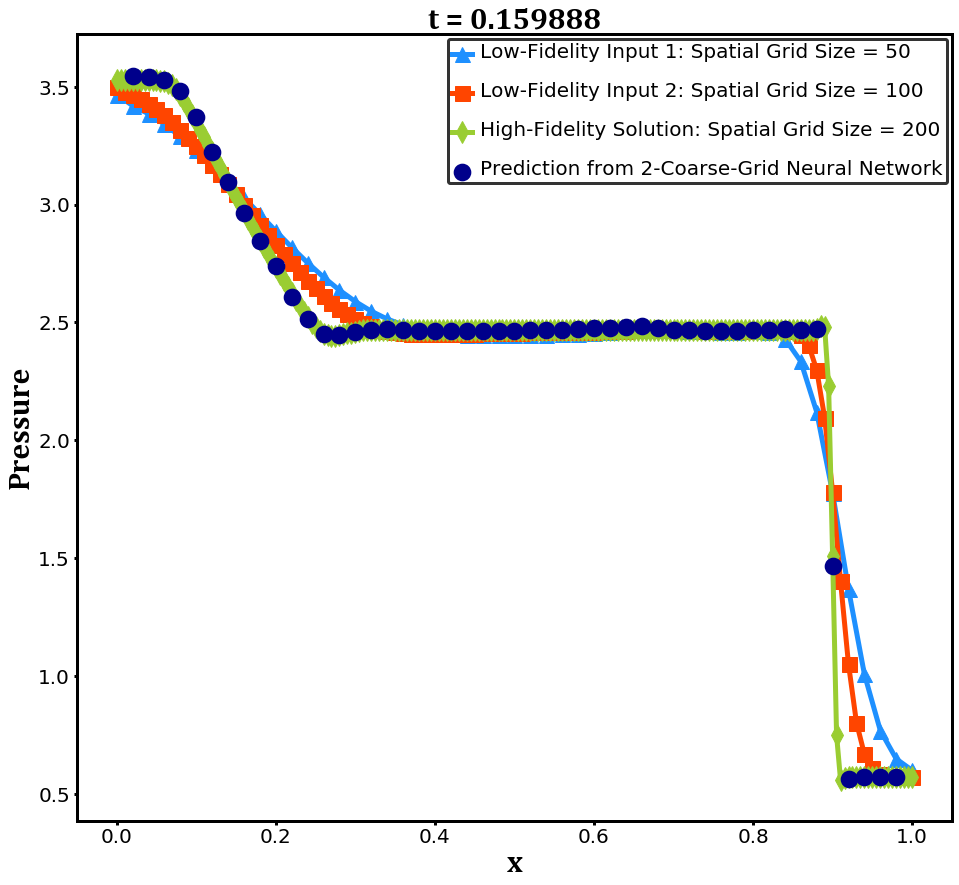}
\end{subfigure}
\begin{subfigure}[b]{.48\textwidth}
  \centering
  \includegraphics[width=1.0\linewidth]{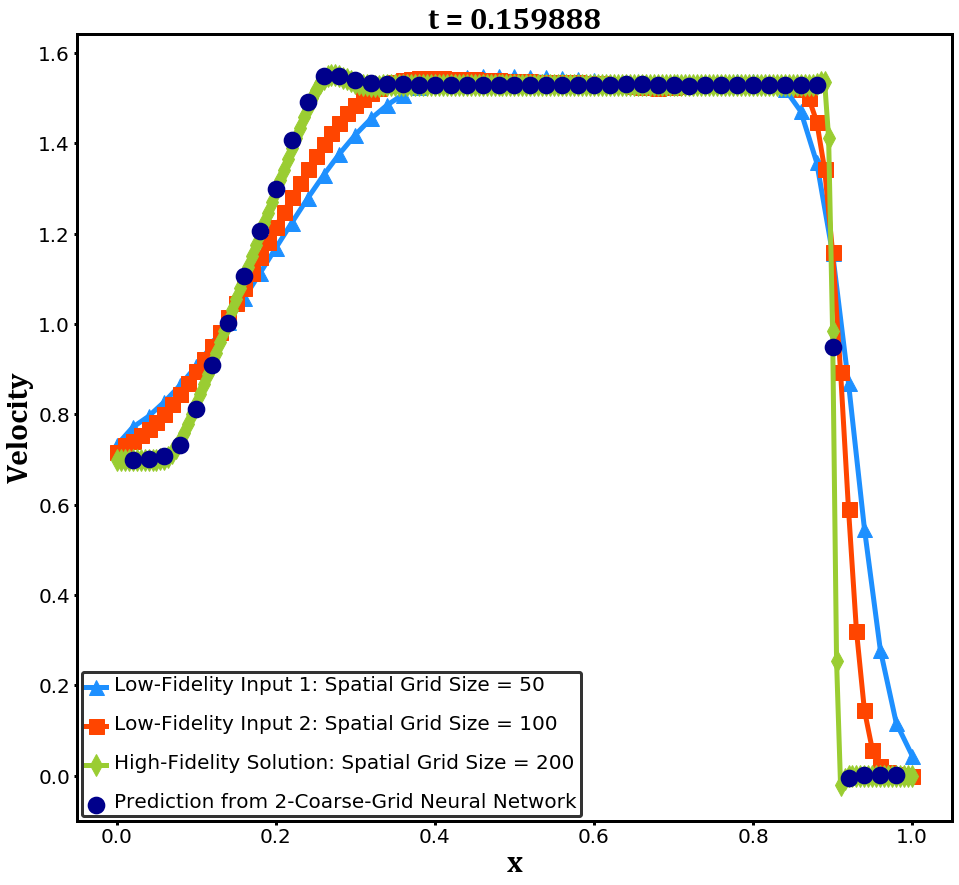}
\end{subfigure}
\captionsetup{labelformat=empty}
\caption{2CGNN prediction of final-time $(t=0.159888)$ solution (dark blue) of the {\bf Lax problem}, low-fidelity input solutions (blue and red) by stabilized Leapfrog and diffusion scheme (\ref{leapfrog-diffusion-stabilized}) on $2$ different grids (with $50$ and $100$ cells resp.), and ``exact'' (reference) solution (green).
}
\label{2CGNN leapfrog diffusion input: Final time of lax problem, original}
\end{figure}

\begin{figure}[H]\centering
\begin{subfigure}[b]{.48\textwidth}
  \centering
  \includegraphics[width=1.0\linewidth]{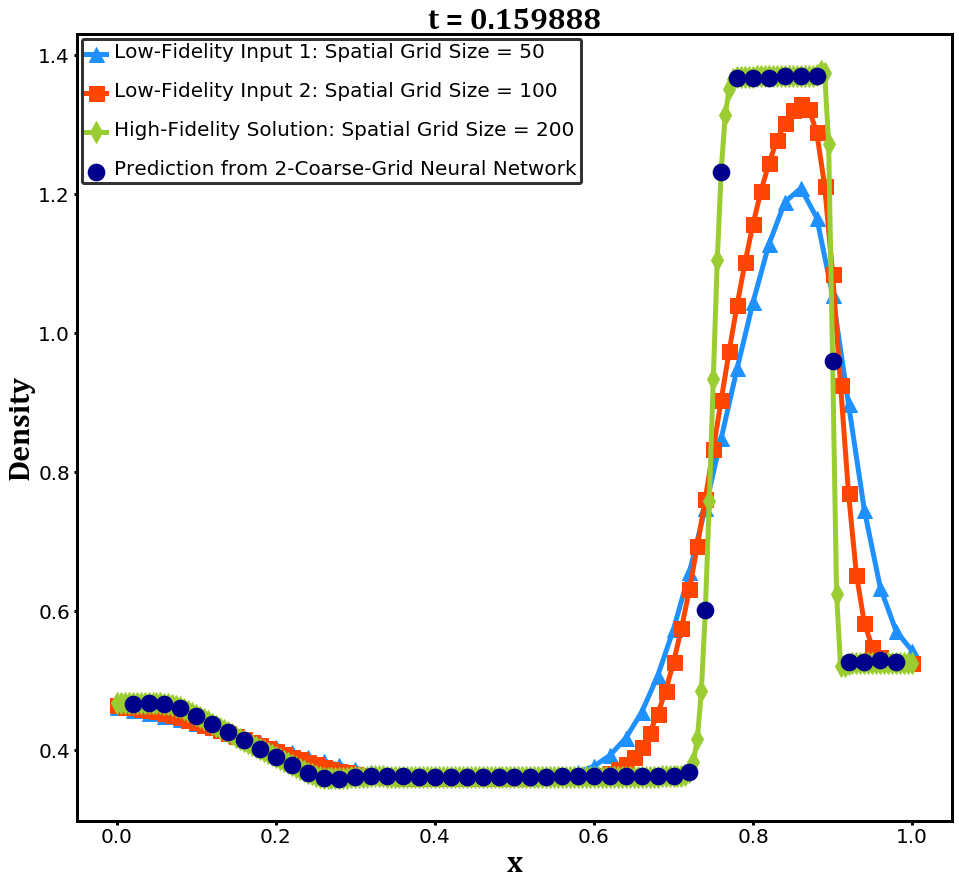}
\end{subfigure}
\begin{subfigure}[b]{.48\textwidth}
  \centering
  \includegraphics[width=1.0\linewidth]{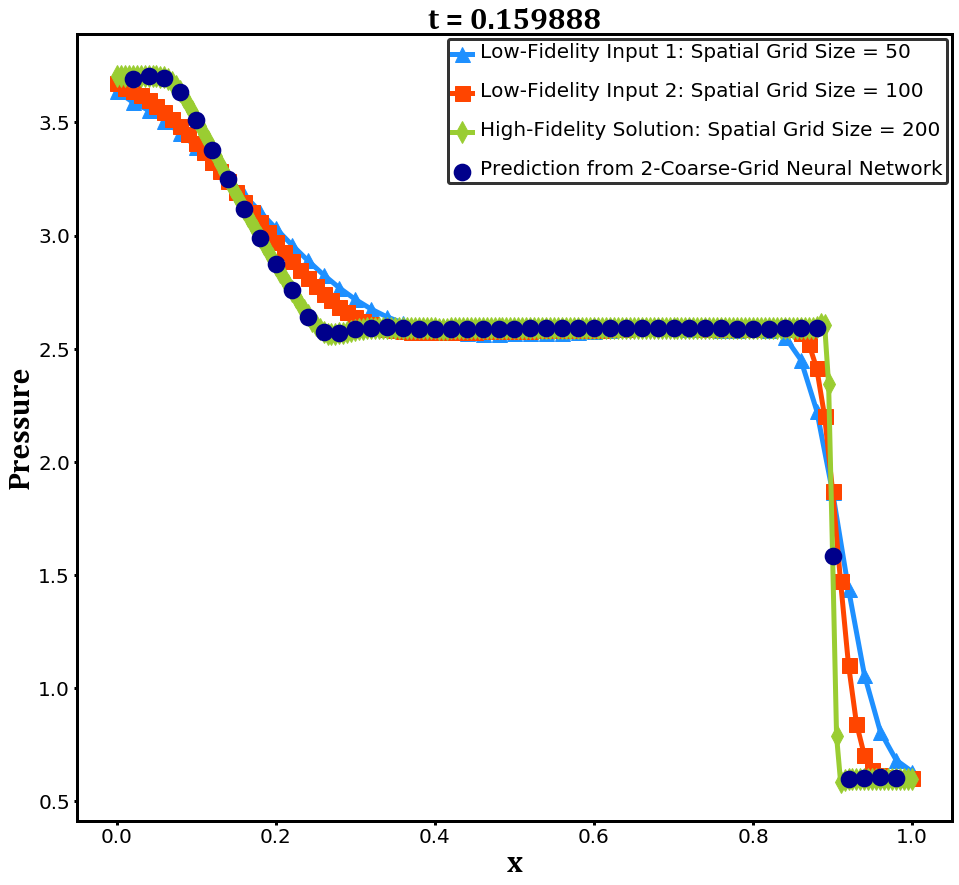}
\end{subfigure}
\begin{subfigure}[b]{.48\textwidth}
  \centering
  \includegraphics[width=1.0\linewidth]{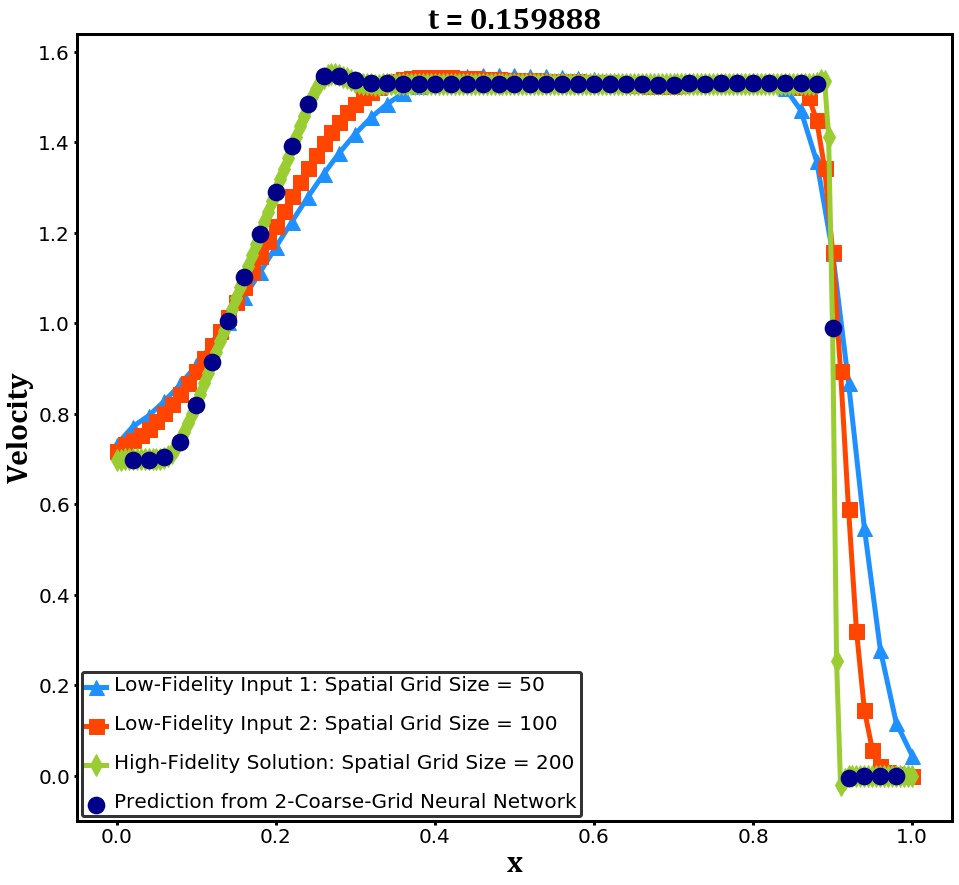}
\end{subfigure}
\captionsetup{labelformat=empty}
\caption{2CGNN prediction of final-time $(t=0.159888)$ solution (dark blue) of the Euler system with its {\bf initial value being $+5\%$ perturbation of that of the Lax problem}, low-fidelity input solutions (blue and red) by stabilized Leapfrog and diffusion scheme (\ref{leapfrog-diffusion-stabilized}) on $2$ different grids (with $50$ and $100$ cells resp.), and ``exact'' (reference) solution (green).
}
\label{2CGNN leapfrog diffusion input: Final time of lax problem, p5}
\end{figure}

\begin{figure}[H]\centering
\begin{subfigure}[b]{.48\textwidth}
  \centering
  \includegraphics[width=1.0\linewidth]{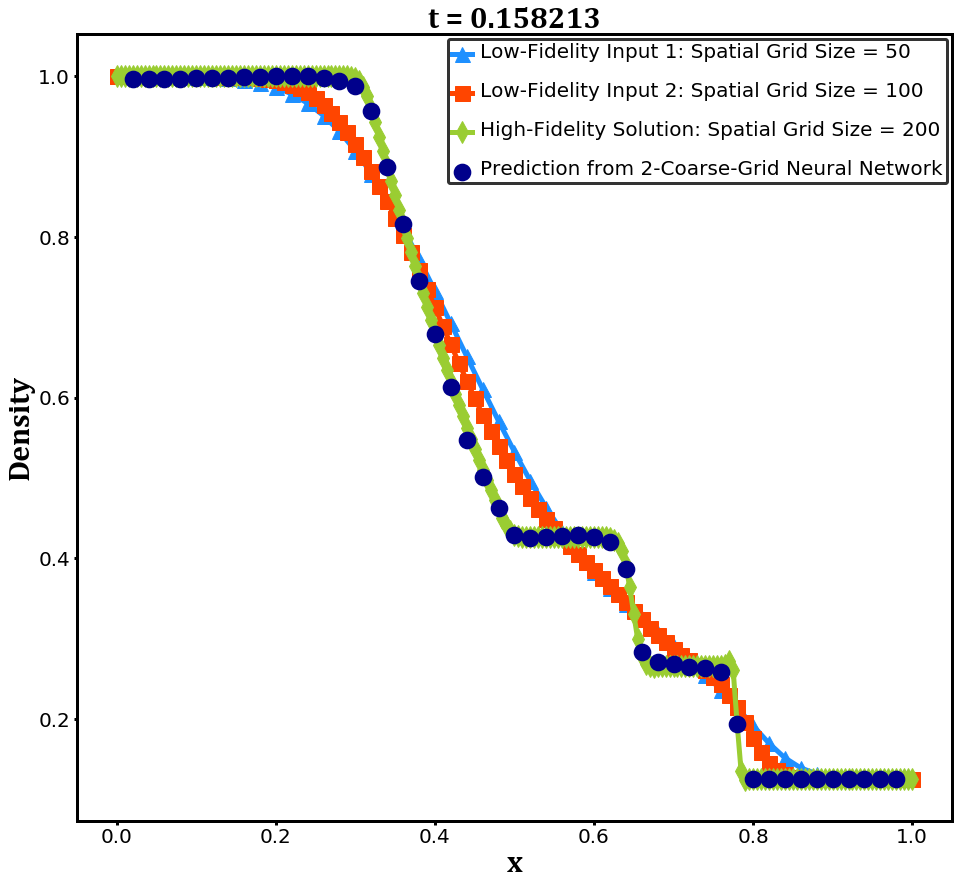}
\end{subfigure}
\begin{subfigure}[b]{.48\textwidth}
  \centering
  \includegraphics[width=1.0\linewidth]{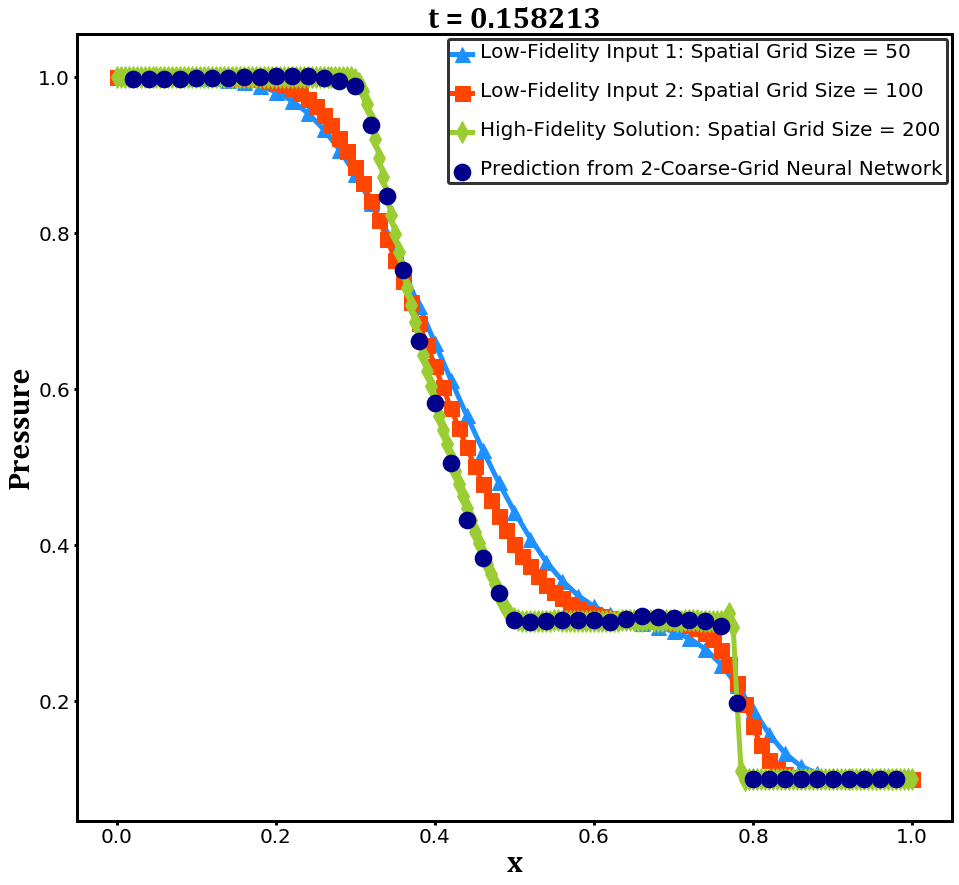}
\end{subfigure}
\begin{subfigure}[b]{.48\textwidth}
  \centering
  \includegraphics[width=1.0\linewidth]{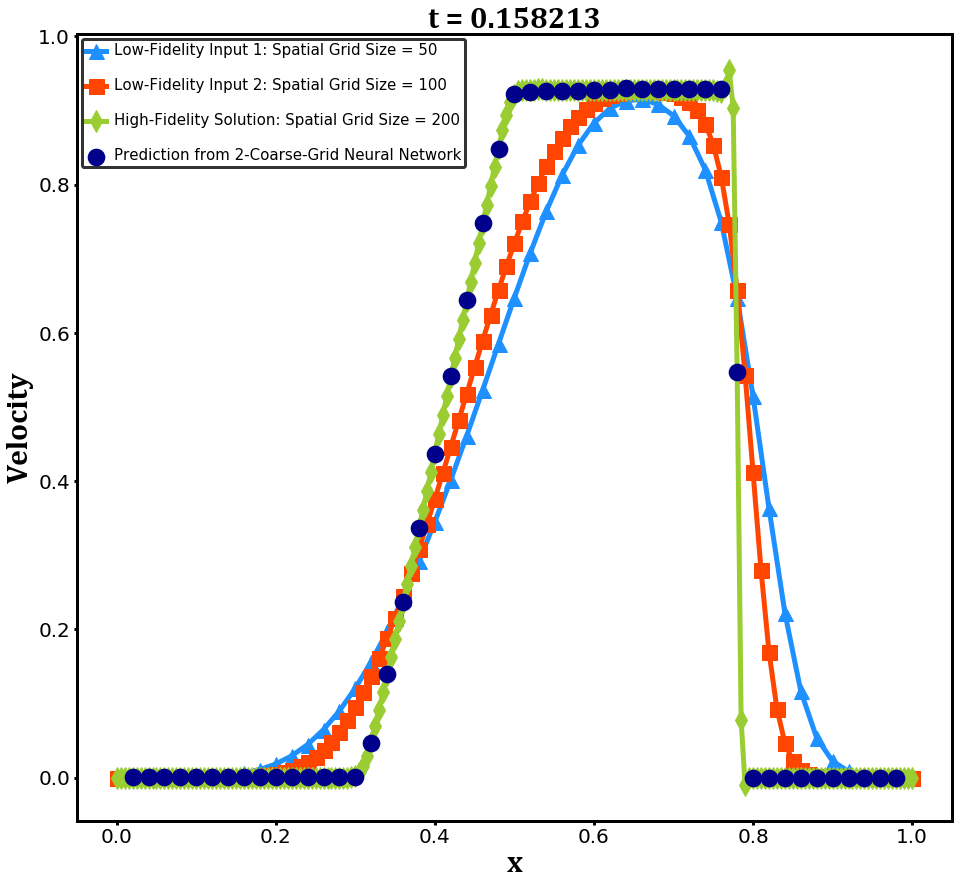}
\end{subfigure}
\captionsetup{labelformat=empty}
\caption{2CGNN prediction of final-time $(t=0.158213)$ solution (dark blue) of the {\bf Sod problem}, low-fidelity input solutions (blue and red) by stabilized Leapfrog and diffusion scheme (\ref{leapfrog-diffusion-stabilized}) on $2$ different grids (with $50$ and $100$ cells resp.), and ``exact'' (reference) solution (green).
}
\label{2CGNN leapfrog diffusion input: Final time of sod problem, original}
\end{figure}

\begin{figure}[H]\centering
\begin{subfigure}[b]{.48\textwidth}
  \centering
  \includegraphics[width=1.0\linewidth]{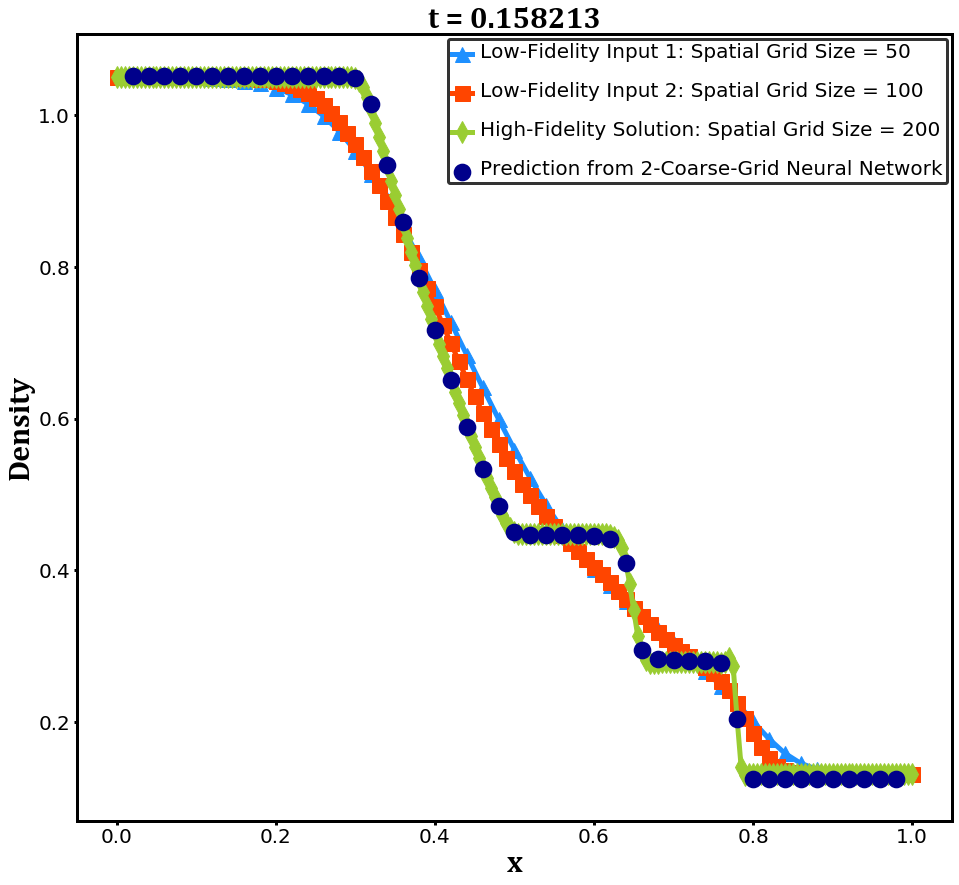}
\end{subfigure}
\begin{subfigure}[b]{.48\textwidth}
  \centering
  \includegraphics[width=1.0\linewidth]{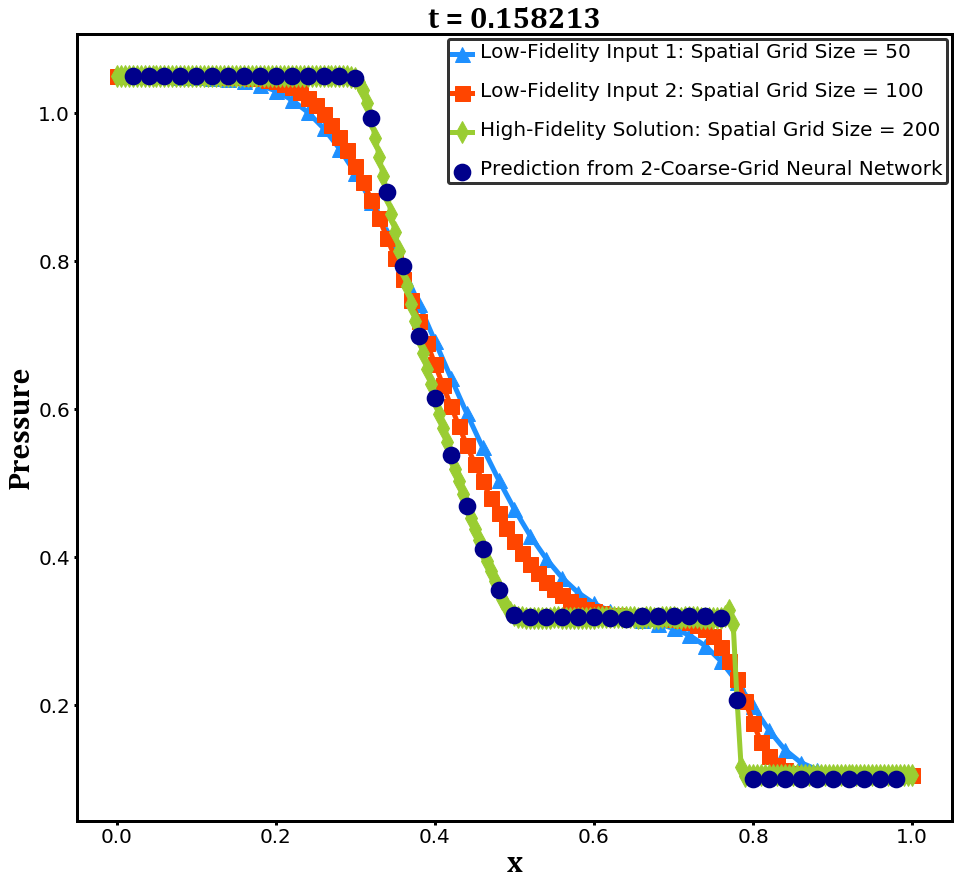}
\end{subfigure}
\begin{subfigure}[b]{.48\textwidth}
  \centering
  \includegraphics[width=1.0\linewidth]{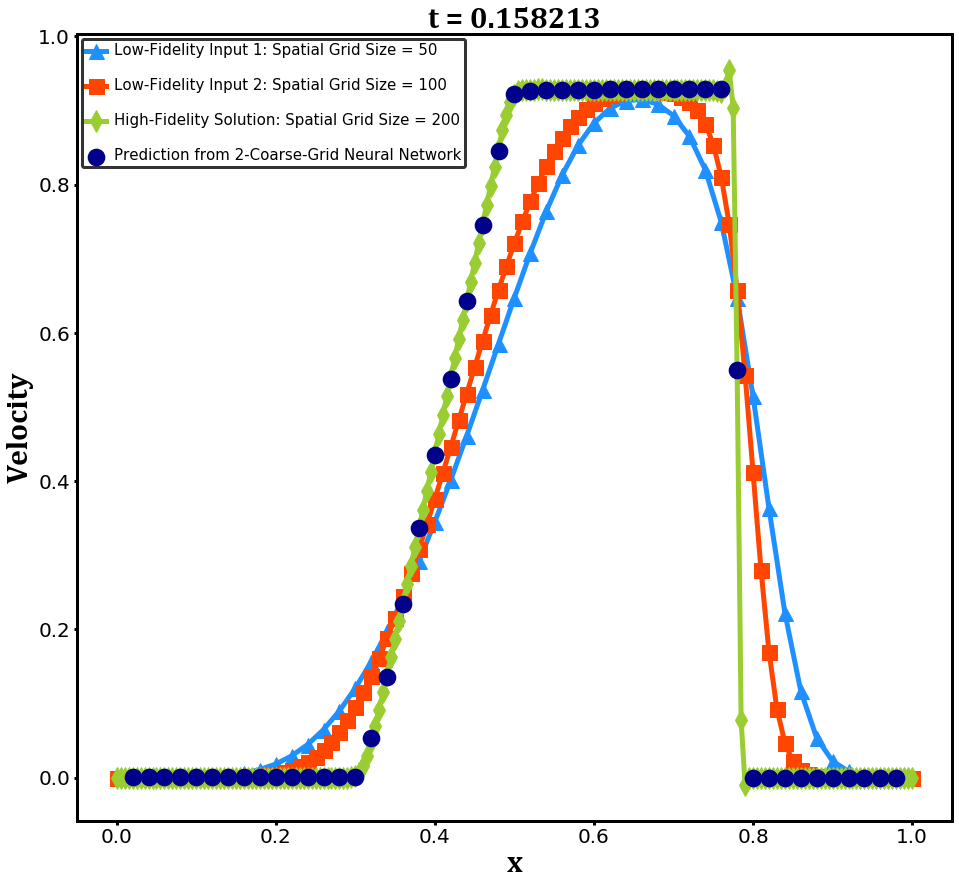}
\end{subfigure}
\captionsetup{labelformat=empty}
\caption{2CGNN prediction of final-time $(t=0.158213)$ solution (dark blue) of the Euler system with its {\bf initial value being $+5\%$ perturbation of that of the Sod problem}, low-fidelity input solutions (blue and red) by stabilized Leapfrog and diffusion scheme (\ref{leapfrog-diffusion-stabilized}) on $2$ different grids (with $50$ and $100$ cells resp.), and ``exact'' (reference) solution (green).
}
\label{2CGNN leapfrog diffusion input: Final time of sod problem, p5}
\end{figure}

\section*{Appendix B. Leapfrog and Diffusion Splitting Scheme (\ref{leapfrog-diffusion-splitting}) for Computing Inputs for 2CGNN}

\begin{figure}[H]\centering
\begin{subfigure}[b]{.48\textwidth}
  \centering
  \includegraphics[width=1.0\linewidth]{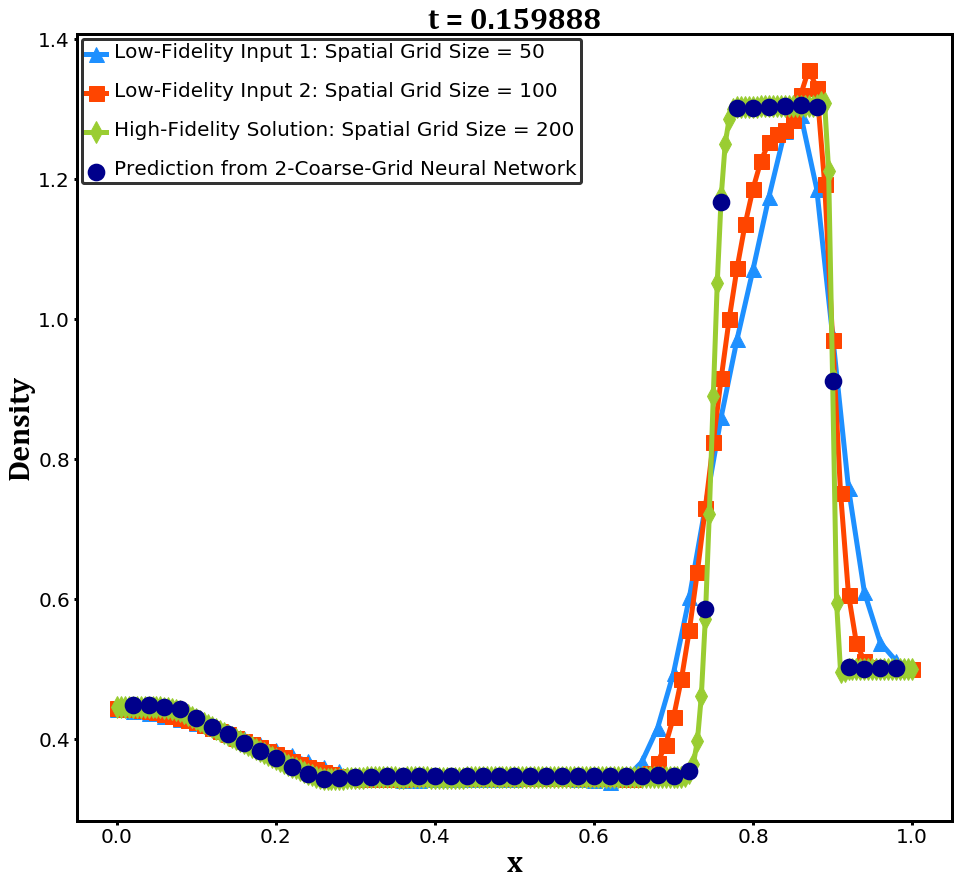}
\end{subfigure}
\begin{subfigure}[b]{.48\textwidth}
  \centering
  \includegraphics[width=1.0\linewidth]{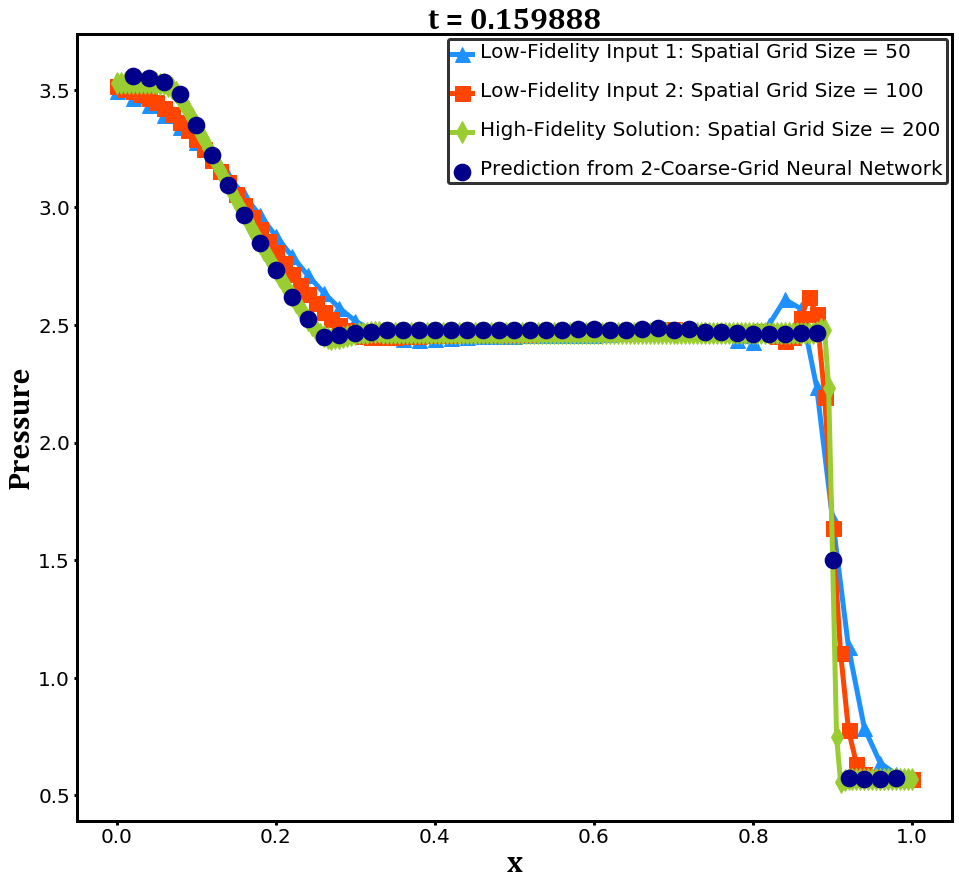}
\end{subfigure}
\begin{subfigure}[b]{.48\textwidth}
  \centering
  \includegraphics[width=1.0\linewidth]{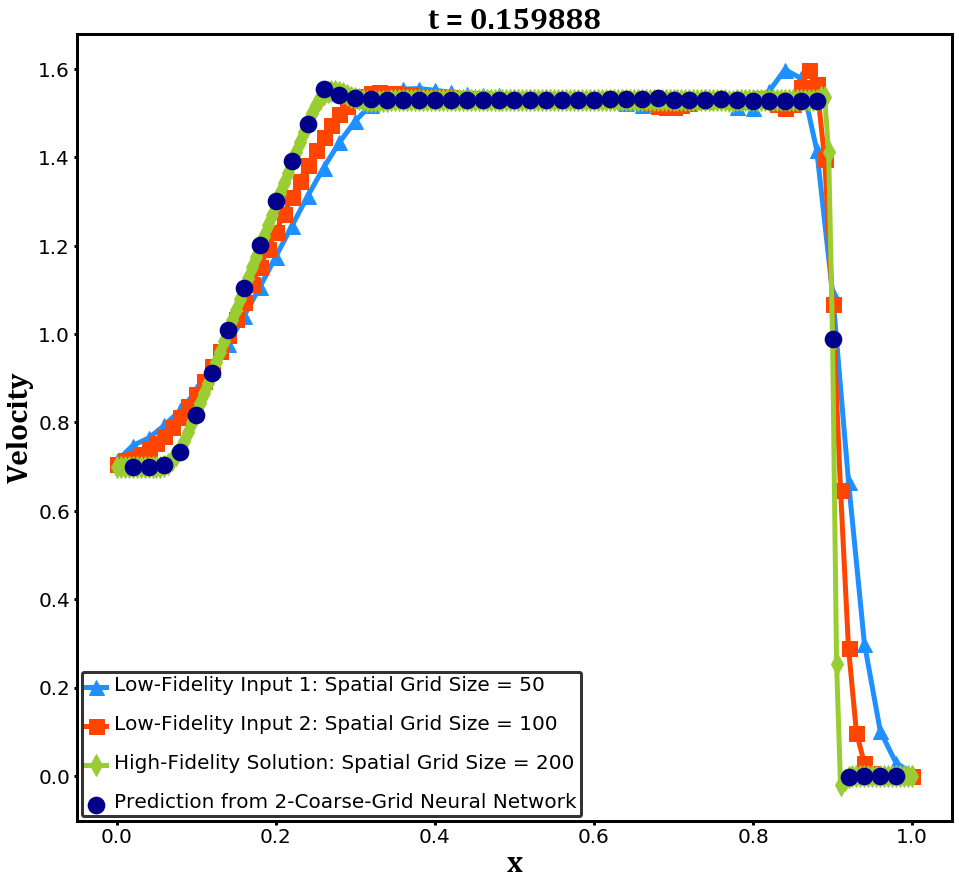}
\end{subfigure}
\captionsetup{labelformat=empty}
\caption{2CGNN prediction of final-time $(t=0.159888)$ solution (dark blue) of the {\bf Lax problem}, low-fidelity input solutions (blue and red) by Leapfrog and diffusion splitting scheme (\ref{leapfrog-diffusion-splitting}) on $2$ different grids (with $50$ and $100$ cells resp.), and ``exact'' (reference) solution (green).
}
\label{2CGNN leapfrog diff w/ split input: Final time of lax problem, original}
\end{figure}

\begin{figure}[H]\centering
\begin{subfigure}[b]{.48\textwidth}
  \centering
  \includegraphics[width=1.0\linewidth]{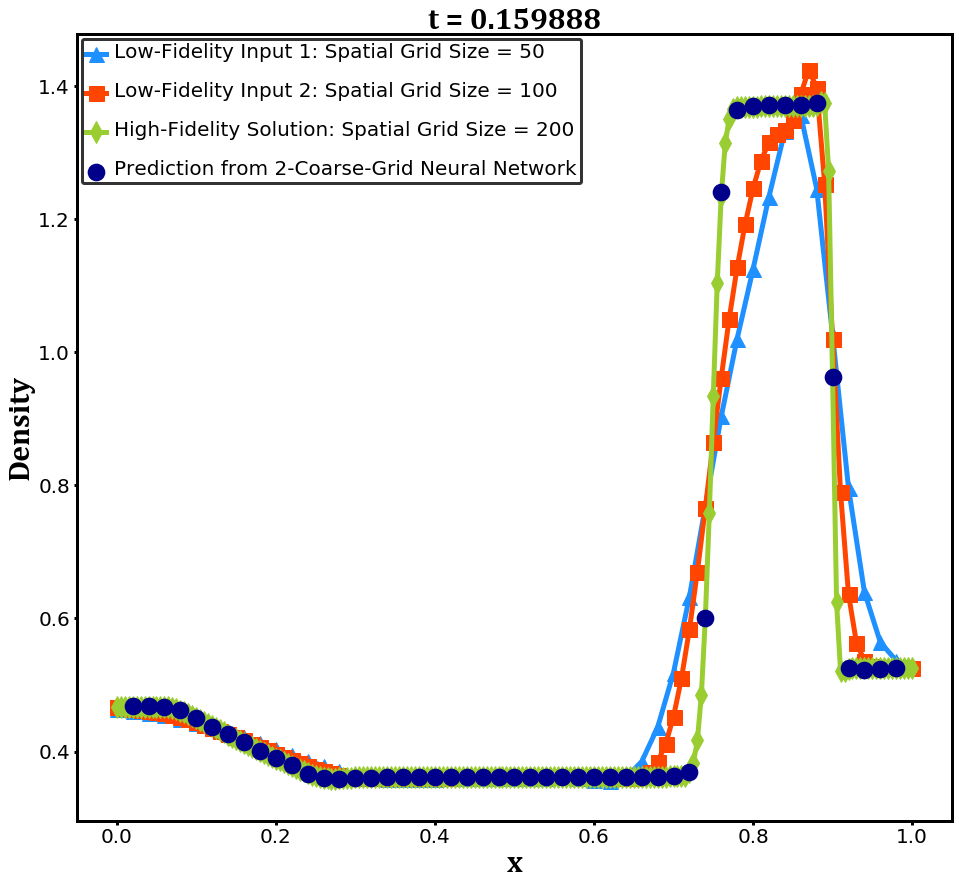}
\end{subfigure}
\begin{subfigure}[b]{.48\textwidth}
  \centering
  \includegraphics[width=1.0\linewidth]{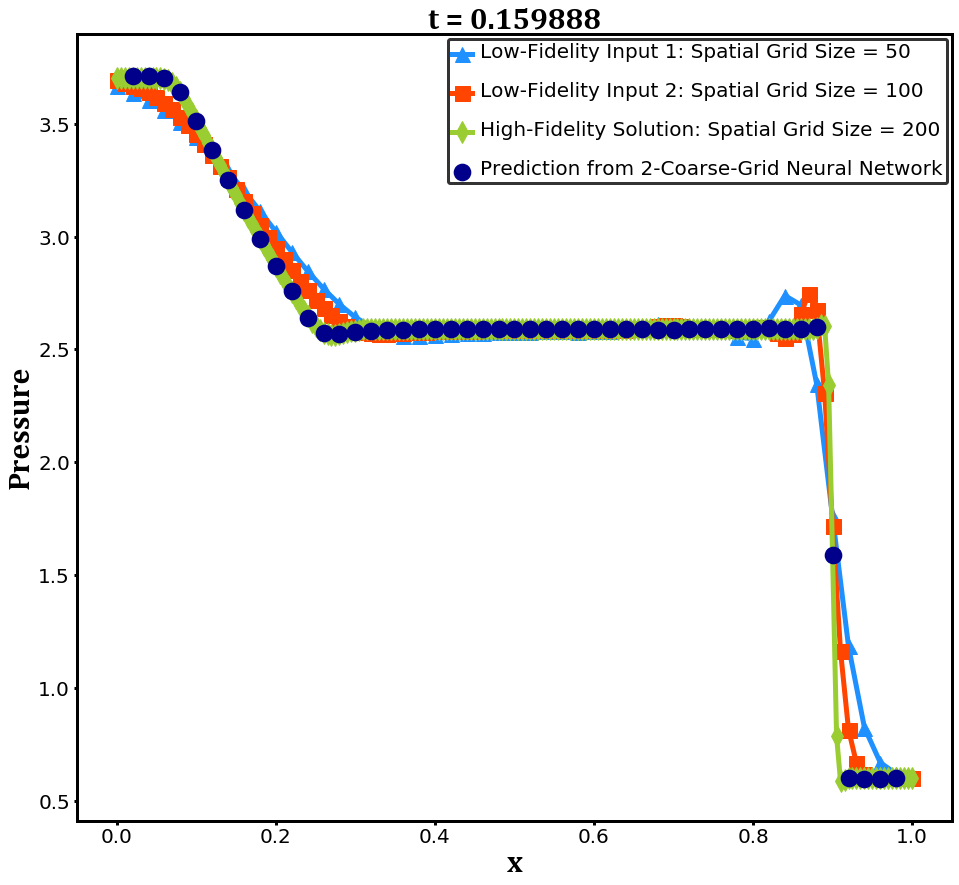}
\end{subfigure}
\begin{subfigure}[b]{.48\textwidth}
  \centering
  \includegraphics[width=1.0\linewidth]{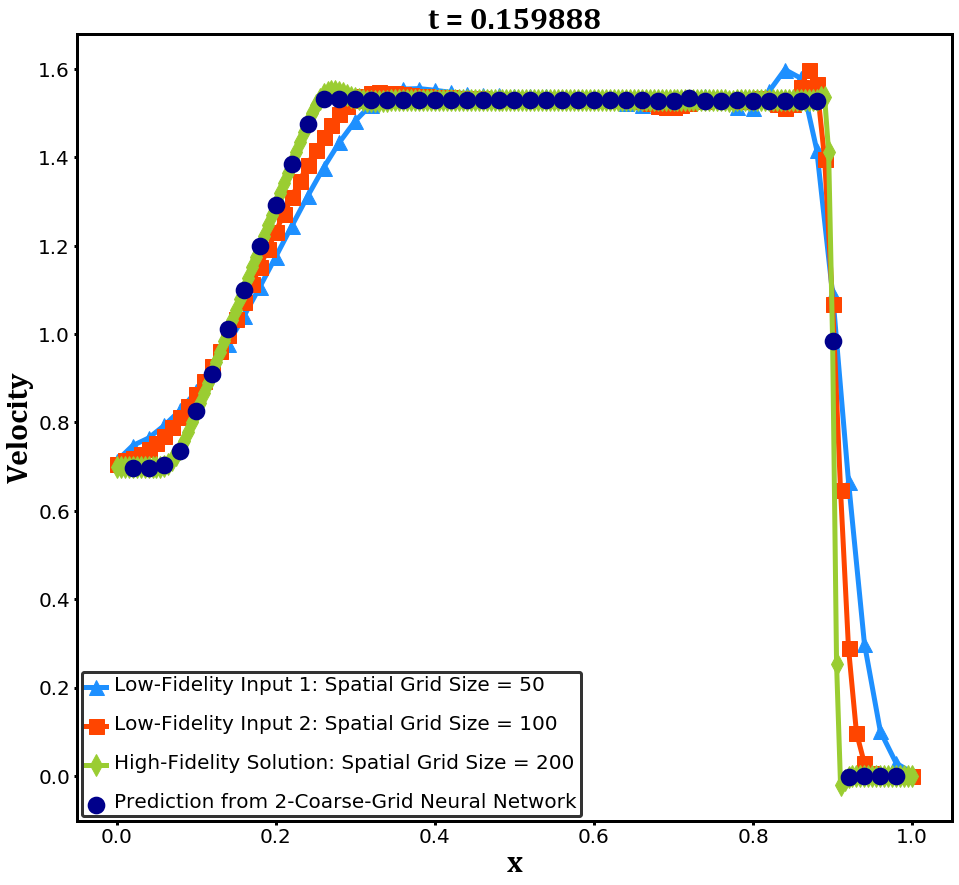}
\end{subfigure}
\captionsetup{labelformat=empty}
\caption{2CGNN prediction of final-time $(t=0.159888)$ solution (dark blue) of the Euler system with its {\bf initial value being $+5\%$ perturbation of that of the Lax problem}, low-fidelity input solutions (blue and red) by Leapfrog and diffusion splitting scheme (\ref{leapfrog-diffusion-splitting}) on $2$ different grids (with $50$ and $100$ cells resp.), and ``exact'' (reference) solution (green).
}
\label{2CGNN leapfrog diff w/ split input: Final time of lax problem, p5}
\end{figure}

\begin{figure}[H]\centering
\begin{subfigure}[b]{.48\textwidth}
  \centering
  \includegraphics[width=1.0\linewidth]{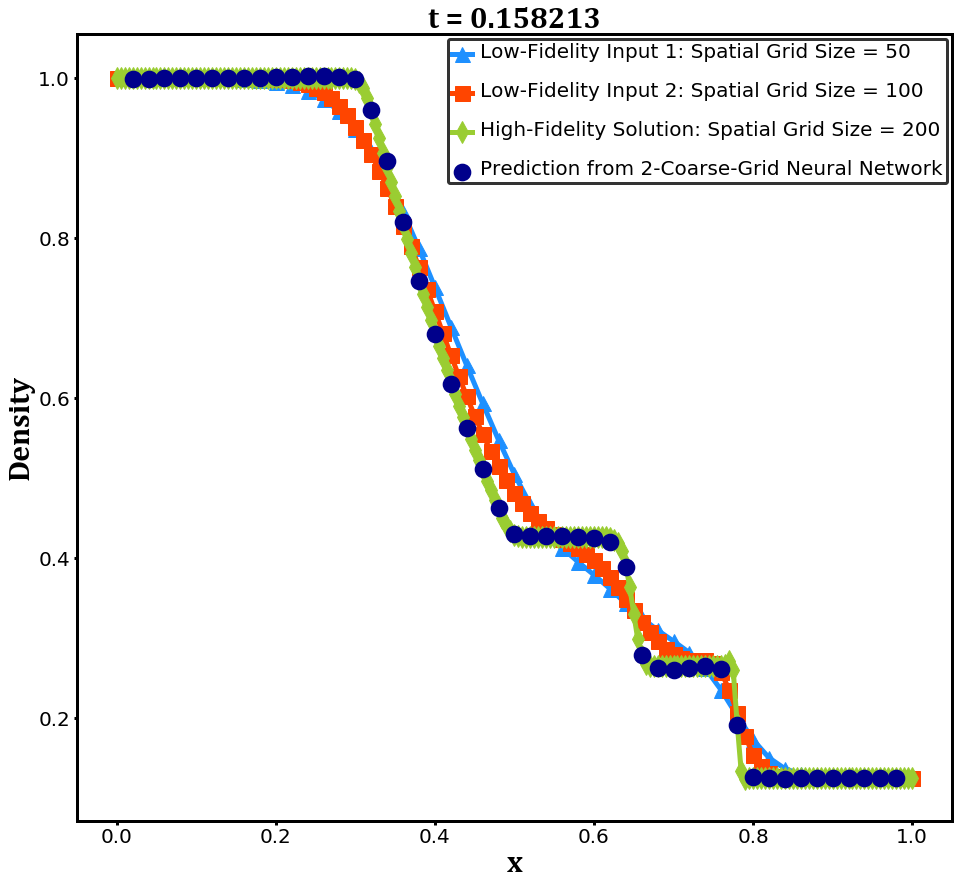}
\end{subfigure}
\begin{subfigure}[b]{.48\textwidth}
  \centering
  \includegraphics[width=1.0\linewidth]{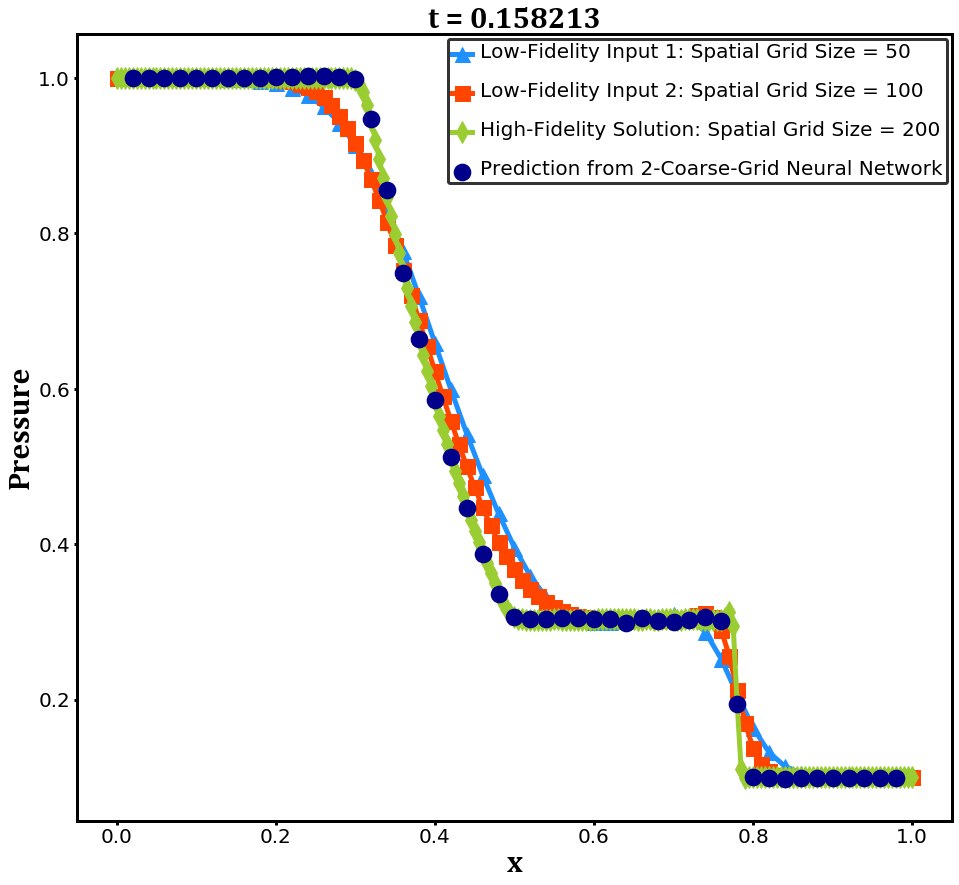}
\end{subfigure}
\begin{subfigure}[b]{.48\textwidth}
  \centering
  \includegraphics[width=1.0\linewidth]{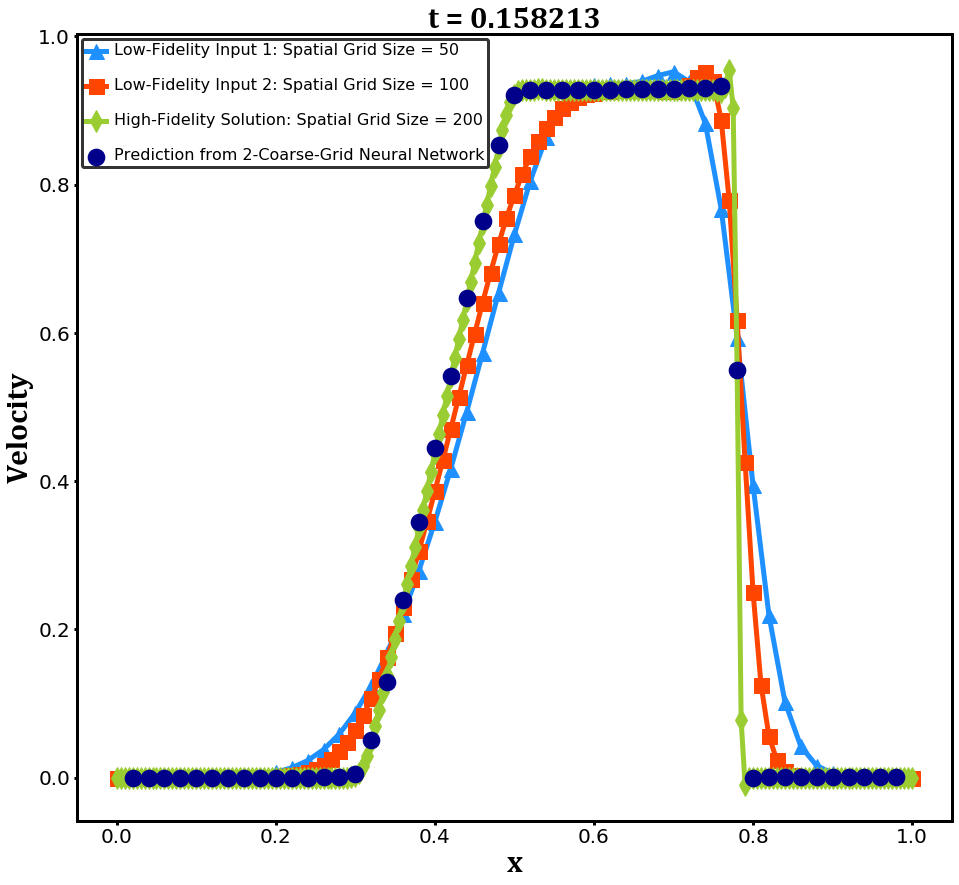}
\end{subfigure}
\captionsetup{labelformat=empty}
\caption{2CGNN prediction of final-time $(t=0.158213)$ solution (dark blue) of the {\bf Sod problem}, low-fidelity input solutions (blue and red) by Leapfrog and diffusion splitting scheme (\ref{leapfrog-diffusion-splitting}) on $2$ different grids (with $50$ and $100$ cells resp.), and ``exact'' (reference) solution (green).
}
\label{2CGNN leapfrog diff w/ split input: Final time of sod problem, original}
\end{figure}

\begin{figure}[H]\centering
\begin{subfigure}[b]{.48\textwidth}
  \centering
  \includegraphics[width=1.0\linewidth]{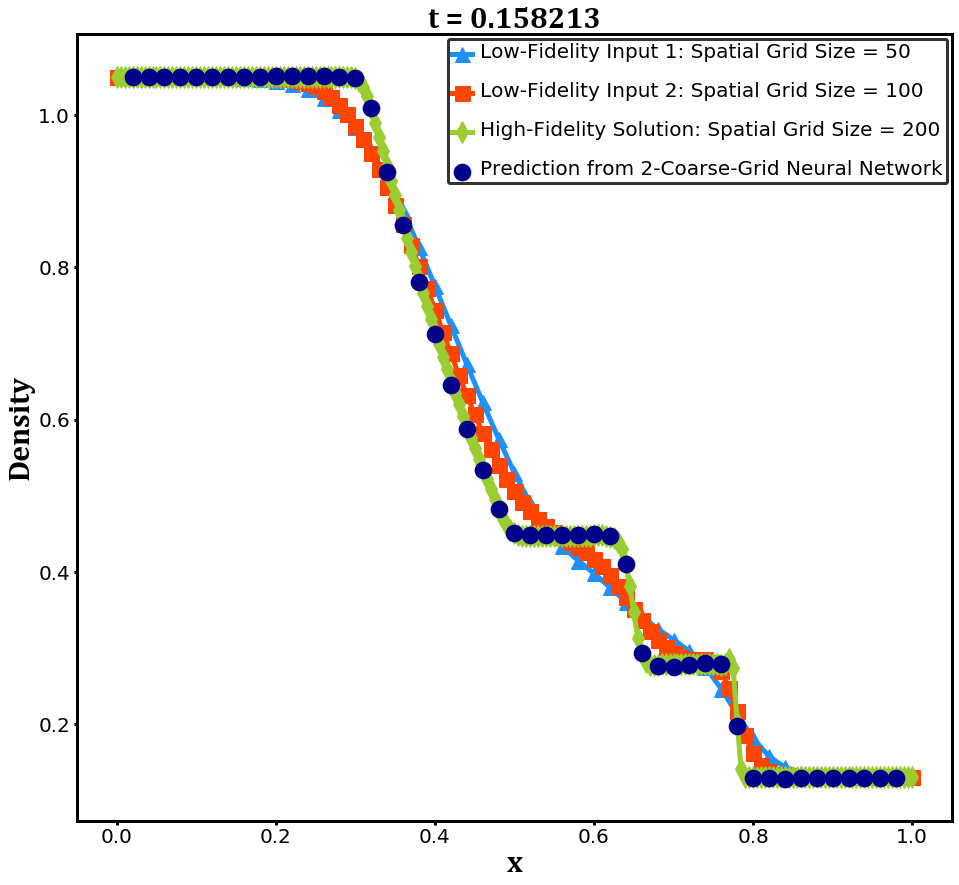}
\end{subfigure}
\begin{subfigure}[b]{.48\textwidth}
  \centering
  \includegraphics[width=1.0\linewidth]{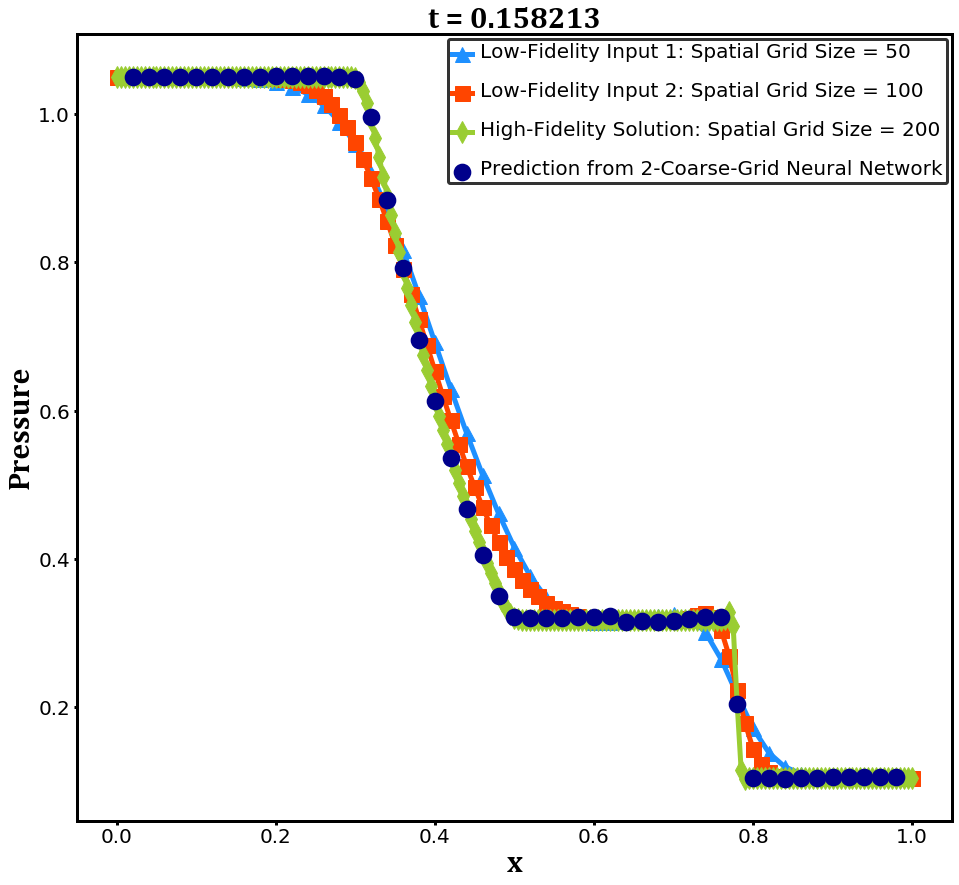}
\end{subfigure}
\begin{subfigure}[b]{.48\textwidth}
  \centering
  \includegraphics[width=1.0\linewidth]{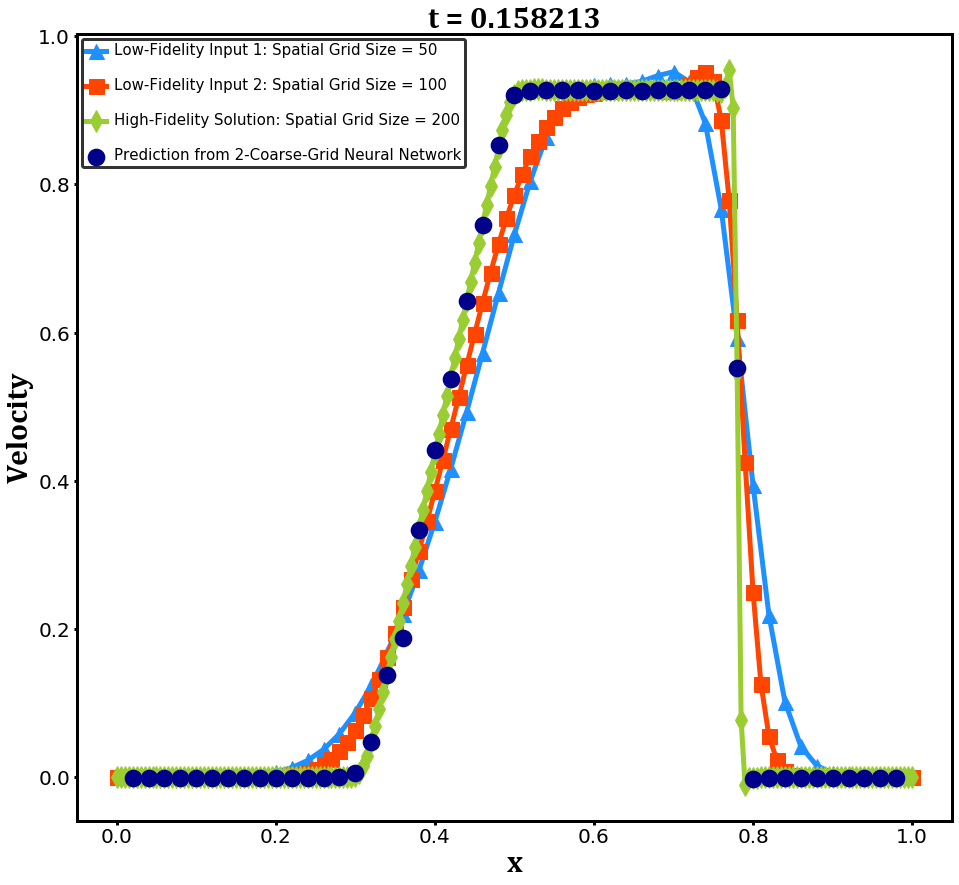}
\end{subfigure}
\captionsetup{labelformat=empty}
\caption{2CGNN prediction of final-time $(t=0.158213)$ solution (dark blue) of the Euler system with its {\bf initial value being $+5\%$ perturbation of that of the Sod problem}, low-fidelity input solutions (blue and red) by Leapfrog and diffusion splitting scheme (\ref{leapfrog-diffusion-splitting}) on $2$ different grids (with $50$ and $100$ cells resp.), and ``exact'' (reference) solution (green).
}
\label{2CGNN leapfrog diff w/ split input: Final time of sod problem, p5}
\end{figure}

\bibliography{mybibfile}

\end{document}